\documentclass[12pt]{article}

\usepackage[a4paper,margin=2cm]{geometry}

\usepackage[T1]{fontenc}
\usepackage[utf8]{inputenc}
\usepackage[style=numeric,citestyle=numeric,url=false]{biblatex}
\usepackage{bbm}
\usepackage{float}
\usepackage{amsmath}
\usepackage{amsthm}
\usepackage{amsfonts}
\usepackage{amssymb}
\usepackage{graphicx}
\usepackage{icomma}
\usepackage{latexsym}
\usepackage{textcomp}
\usepackage[all]{xy}
\usepackage{cancel}
\usepackage{dsfont}
\usepackage{caption}
\usepackage{listings}
\usepackage{lmodern}
\usepackage{xcolor}
\usepackage{moreverb}
\usepackage{mathrsfs}
\usepackage{stmaryrd}
\usepackage{thmtools}
\usepackage{hyperref}
\usepackage{cleveref}

\newcommand{\suml}[0]{\sum\limits}

\newcommand{\mbf}[1]{\mathbf{#1}}
\newcommand{\mscr}[1]{\mathscr{#1}}
\newcommand{\mcal}[1]{\mathcal{#1}}
\newcommand{\mbb}[1]{\mathbb{#1}}

\newcommand{\x}[0]{\mbf{x}}

\newcommand{\etal}[0]{\textit{et al. }}
\newcommand{\PX}[0]{\mscr{P}_X}
\newcommand{\ds}[0]{\displaystyle}

\newcommand{\slfrac}[2]{\left.#1\middle/#2\right.}

\newtheorem{lemme}{Lemma}

\newtheorem{hypothese}{Assumption}
\newtheorem{prop}{Proposition}
\newtheorem{theoreme}{Theorem}
\newtheorem{coro}{Corollary}

\title{Robust estimation for ergodic Markovian processes}
\author{Alexandre Lecestre\thanks{This project has received funding from the European Union's Horizon 2020 research and innovation programme under grant agreement N° 811017}}

\begin{document}
\maketitle
\begin{abstract}
We observe $n$ possibly dependent random variables, the distribution of which is presumed to be stationary even though this might not be true, and we aim at estimating the stationary distribution. We establish a non-asymptotic deviation bound for the Hellinger distance between the target distribution and our estimator. If the dependence within the observations is small, the estimator performs as good as if the data were independent and identically distributed. In addition our estimator is robust to misspecification and contamination. If the dependence is too high but the observed process is mixing, we can select a subset of observations that is almost independent and retrieve results similar to what we have in the i.i.d. case. We apply our procedure to the estimation of the invariant distribution of a diffusion process and to finite state space hidden Markov models.
%We also propose a procedure to automatically choose a subset of observations when we dispose of two independent sets of observations.
\end{abstract}
\section{Introduction}
We observe $n$ random variables $X_1,\dots,X_n$ with common distribution $P$ which is assumed to belong, or at least to be close enough, to a given model $\mscr{M}$. Our aim is to estimate $P$ with an estimator $\hat{P}$ taking values in $\mscr{M}$. These random variables are not necessarily independent however we assume that for indices $i\neq j$ with $|i-j|$ large enough, the distribution of the couple $(X_i,X_j)$ is close to $P\otimes P$. We also want our estimator to be robust to contamination and outliers.\par
When we actually dispose of an independent sample, this problem has already been investigated in Baraud \etal \cite{baraudinventiones} and Baraud \& Birgé \cite{baraudrevisited}. They provide a non-asymptotic deviation bound for the Hellinger distance $h$ between $P$ and their $\rho$-estimator. For two probability distributions $P$ and $Q$ on the same measurable space, the Hellinger distance $h(P,Q)$ between $P$ and $Q$ is given by
\[ h^2\left( P,Q \right) = \frac{1}{2} \int\left( \sqrt{dP/d\mu}-\sqrt{dQ/d\mu}\right)^2 d\mu, \]
where $\mu$ is any measure that dominates both $P$ and $Q$, the result being independent of $\mu$. It is shown in those articles that the $\rho$-estimator is robust in the following sense. Even if the variables $X_i$ do not have a common distribution $P$ but marginals $P_i$ such that most of them are relatively close to a distribution $P\in\mscr{M}$, then the $\rho$-estimator is almost as efficient as when the data is i.i.d.with common distribution $P$. The obtained risk bounds are minimax, up to a logarithmic factor, when the model is well-specified and are not significantly deteriorated as long as the approximation term $n^{-1} \suml_{i=1}^n h^2(P_i,P)$ is relatively small in the misspecified case.\par
We want to obtain similar results when we do not satisfy the independence assumption but the observations are almost independent. This can happen for processes with mixing properties. We only focus on the theoretical aspects and performances of our estimation method. We prove a general result, Theorem \ref{th:main_rho}, which gives a bound in expectation for the risk of our estimator $\hat{P}$ with respect to an Hellinger-type loss. This result is free of any assumption on the data and the risk bound is the sum of three terms: the approximation term mentioned above, a dimension which measures the complexity of the model $\mscr{M}$, and a dependence term which measures how far the observations are from being independent. We quantify the dependence within the sample using Kullback-Leibler divergence of the joint distribution from the product of the marginal distributions. Our risk bound is as good as when the data is independent as long as the dependence term is not bigger than the other terms. We have the following approach for when the dependence term is too big. We split our data in order to get a subset of the original observations for which the dependence term is small enough.\par
We apply this method for the estimation of an invariant distribution of a discretely observed diffusion process. Under some condition the stationary solution of a Langevin equation is mixing and its invariant distribution has a $\log$-concave density with respect to the Lebesgue measure. We can refer to the literature on the estimation of a $\log$-concave density in the i.i.d. context and adapt our procedure to this situation. We obtain convergence rates for our estimator in any dimensions. Those rates are similar to the minimax rates for i.i.d.\ estimation, with a worse logarithmic power.\par
Our main application is hidden Markov models (HMMs). These models are widely applied to model state dependent processes where the state process is Markovian but is not observed. We refer the interested reader to Mor, Garhwal and Kumar \cite{review} for a review of applications of HMMs. Let $Y_1,\dots,Y_N,H_1,\dots,H_N$ be random variables. We say that $(Y_i,H_i)_{1\leq i\leq N}$ is a hidden Markov model (HMM) if $\left(H_i\right)_i$ is a Markov chain and each variable $Y_i$ only depends on the associated $H_i$. In particular the variables $Y_1,\dots,Y_N$ are independent conditionally on $\left(H_i\right)_i$. It is called a hidden Markov model as the Markov chain $\left(H_i\right)_i$ is typically not observed and $\left(Y_i\right)_i$ is the only accessible data.\par
We focus on homogeneous finite state space HMMs. Such processes can be completely described by the number $K$ of hidden states $h_1,\ldots,h_K$, the initial distribution $w$ and the transition matrix $Q$ of the hidden Markov chain, and the set of emission distributions $F=(F_1,\dots,F_K)$, where $F_k$ is the conditional distribution of $Y_i$ given $H_i=h_k$. In that case we say that $(Y_i,H_i)_i$ is a HMM with parameters $(K,w,Q,F)$. Because the hidden state space does not have a particular importance, we will always assume it is of the form $\{1,2,\ldots,K\}$. For a particular class of distributions $\mscr{F}$ there is a minimal value of $K$ such that $(Y_i,H_i)_i$ is a HMM with parameters $(K,w,Q,F)$ with $F_1,\dots,F_K\in\mscr{F}$. This value of $K$ is called the order of the HMM (with respect to $\mscr{F}$). Typically one aims at estimating these parameters from stationary observations $(Y_i)_{1\leq i\leq N}$.\par
Numerous estimation methods have been developed to estimate some or all of the parameters. Cappé \etal \cite{Cappe2005} provide an overall survey of the different results in the literature. Most theoretical guarantees are either asymptotic or restricted to specific parametric models. Leh\'ericy \cite{lehericyorder} provided non-parametric and non-asymptotic results for a penalized least squares estimator with the following approach. They first estimate the distribution $P_L=P_{\pi^*,Q^*,F^*}$ of $L$ consecutive observations $Y_i,Y_{i+1},\dots,Y_{i+L-1}$ of a stationary ergodic HMM with parameters $(K^*,\pi^*,Q^*,F^*)$, where $P_{w,Q,F}$ is defined by
\begin{equation}
\label{eq:p_w_q_f}
P_{w,Q,F} = \suml_{1\leq k_1,\ldots,k_L\leq K} w_{k_1} Q_{k_1,k_2} \dots Q_{k_{L-1},k_L} \bigotimes_{l=1}^L F_{k_l}.
\end{equation}
They use model selection to consistently estimate the order $K^*$. When the estimation of the order is correct, it is possible to deduce the different parameters from $P_L$ for $L$ large enough. They show that $L\geq 3$ is enough for linearly independent emission densities. They lower bound the $L^2$-distance between densities by a distance on the parameters. Therefore a risk bound for the estimation of $P_L$ is enough to obtain risk bounds for the parameter estimators.\par
However their estimator is not robust to misspecification nor to contamination and there is no estimator that tackles this problem for general finite state space HMMs. The estimation method we propose aim at solving this problem. For the sake of simplicity we do not aim at estimating the order $K^*$. We do not look into this particular aspect in this paper however model selection can be considered to choose automatically an order from the data. This is to be treated in a subsequent paper.\par
%
%We refer the interested reader to Anandkumar \etal \cite{anandkumar}, de Castro \etal \cite{decastro_lacour} and Lehéricy \cite{lehericyorder} for the problem of order estimation. \par
%
We use the tools we develop in the first part with $\mscr{M}$ containing distributions of the form $P_{w,Q,F}$ to obtain a robust estimator $\hat{P}$ of $P_L$, hence $\hat{P}$ being of the form $\hat{P}=P_{\hat{w},\hat{Q},\hat{F}}$. We have a general risk bound for $\hat{P}$ which is free of any assumption on the data from which we obtain convergence rates when we assume that the observations come from an ergodic finite state space HMM. In particular the stationarity of the observations is not necessary. We show that the performance of our estimator is not significantly worsened when the model is misspecified as long as the distance to the true distribution is small compared to the rate we have in the well-specified case. Similarly the performance of our estimator is not deteriorated by contamination as long as the contamination rate is not too big.
\par
We can deduce risk bounds for the parameter estimators $\hat{w},\hat{Q},\hat{F}$ under some conditions on the model $\mscr{M}$. We need an inequality of the form 
\begin{equation}
\label{eq:parameters_hellinger}
d\left( (w,Q,F), (\overline{w},\overline{Q},\overline{F}) \right) \leq C\left(\overline{w},\overline{Q},\overline{F}\right) h^2\left( P_{w,Q,F}, P_{\overline{w},\overline{Q},\overline{F}} \right), \forall P_{w,Q,F}\in\mscr{M}.
\end{equation}
We obtain convergence rates for the estimation of the parameters when the model is well specified. If the model is misspecified but $\overline{P}=P_{\overline{w},\overline{Q},\overline{F}}$ is the best approximation of $P_L$ within our model our estimators $\hat{w},\hat{Q},\hat{F}$ should be close to $\overline{w},\overline{Q},\overline{F}$ when this approximation is relatively good.\par
It is possible to use the results that already exist for the $L_2$-norm to obtain an inequality like (\ref{eq:parameters_hellinger}) when the densities are bounded. For two probability distributions $P$, $Q$ dominated by a positive measure $\mu$, we have
\begin{equation}
\label{eq:l_2_hellinger}
||p-q||_2^2 \leq 4 (||p||_{\infty}+||q||_{\infty}) h^2(P,Q),
\end{equation}
where $p=dP/d\mu$ and $q=dQ/d\mu$. It is also possible to prove inequalities directly for the Hellinger distance in some cases. We do so for models with emission densities that belong to exponential families with some regularity. We also consider an example with classes of emission densities that are unbounded and not even square integrable in some cases. For this example we obtain rates that are faster than the parametric rate for one of the parameters. Classical estimators such as the maximum likelihood or least-squares estimators do not apply as the considered densities are unbounded.\par
Our estimation method requires that the statistician selects themself a subset of the observations that should be almost independent. This is not possible without any knowledge on the distribution of the data. We propose to overcome this restriction and provide a way to automatically select an almost independent subset of observations when we dispose of a second set of observations independent from the first one. We obtain a general risk bound and show that for ergodic HMMs we retrieve the same rate of convergence as when the optimal way of selecting observations is known. This method is still robust to misspecification and contamination.\par
The paper is organized as follows. In Section \ref{sec:first_section}, we present our estimation procedure and our main result in a general framework. We consider the application to the estimation of the invariant distribution of a diffusion process in Section \ref{sec:exponentially_ergodic_estimation}. We dedicate Section \ref{sec:hmm} to finite state space hidden Markov models. Finally, we propose a complete procedure for situations in which we do not know the mixing regime in Section \ref{sec:selection_s}. The proofs of all the different results can be found in the appendix.\par
\textbf{Notation.} For a set $A$, we denote by $|A|$ its cardinal which can be infinite. For an integer $k$, we denote by $[k]$ the set $\{1,2,\dots,k\}$. We denote by $\mbb{R}_+$ the set of non-negative real numbers. For a real number $x$, we denote by $\lceil x\rceil$ (resp. $\lfloor x\rfloor$) the only integer $k$ satisfying $k-1 <x\leq k$ (resp. $k\leq x<k+1$). For a random variable $X$ we denote by $\mcal{L}(X)$ its probability distribution. The notation $C(\theta,\alpha,\beta)$ means that $C(\theta,\alpha,\beta)$ is a constant that depends on the parameters $\theta$, $\alpha$ and $\beta$. It can change from one inequality to the other. On the other hand a constant written $C$ will be universal. For a real number $x$ we denote by $x_+$ its positive part given by $x_+=x\vee 0$.
\section{Construction of the estimator and main result}
\label{sec:first_section}
Let $X_1,\dots,X_n$ be $n$ possibly dependent random variables on the measurable space $\left(\mscr{X},\mcal{X}\right)$. Our aim is to estimate their marginal distribution $P^*$ doing as if they were identically distributed, even though this might not be exactly the case. We denote by $\PX$ the class of all probability distribution on $\left(\mscr{X},\mcal{X}\right)$ and for $i\in[n]$ by $P_i=\mcal{L}(X_i)\in\PX$ the true marginal distribution of $X_i$. We also want our estimator of $P^*$ to be robust to misspecification, contamination and outliers. The $\rho$-estimators developed by Baraud, Birgé and Sart in \cite{baraudinventiones} and \cite{baraudrevisited} are perfectly adapted to this task when the observations are independent. We prove that their performances remain almost as good when the observations are close to being independent.
\subsection{Reminders of \texorpdfstring{$\rho$}{rho}-estimation}
\label{sec:rho_estimation}
We denote by $\psi$ the function given by
\begin{equation}
\label{eq:psi_rho}
\psi:\left|\begin{array}{l} [0,+\infty]\rightarrow[-1,1]\\x\mapsto \frac{x-1}{x+1} \end{array}\right. .
\end{equation}
Let $\mscr{M}$ be a countable subset of $\mscr{P}_X$ such that there is an associated set of density functions $\mcal{M}$ with respect to a $\sigma$-finite measure $\mu$. For $n\geq 1$, we denote by $\mbf{T}_n$ and $\mbf{\Upsilon}_n$ the functions given by 
\begin{equation}
\label{eq:t_n}
\mbf{T}_n : \left| \begin{array}{l} \mscr{X}^n\times \mcal{M}\times\mcal{M} \rightarrow [-1,1]\\
(\mbf{x},q,q') \mapsto \suml_{k=1}^n \psi\left(\sqrt{\frac{q'(x_i)}{q(x_i)}}\right)
                   \end{array}\right.
\end{equation}
with the convention $0/0=1$, $a/0=+\infty$ for all $a>0$, and
\begin{equation}
\label{eq:upsilon}
\mbf{\Upsilon}_n : \left| \begin{array}{l}
\mscr{X}^n \times \mcal{M}\\
\left(\mbf{x},q\right) \mapsto \sup_{q'\in\mcal{M}} \mbf{T}_n\left(\mbf{x},q,q'\right)
\end{array}\right. .
\end{equation}
For $\mbf{x}$ in $\mscr{X}^n$, we define the (nonvoid) set $\mscr{E}_n(\mbf{x})$ by 
\begin{equation}
\label{eq:rho_estimateur}
\mscr{E}_n(\mbf{x}) = \left\{ Q=q\cdot\mu \bigg| q\in\mcal{M}, \mbf{\Upsilon}_n\left(\mbf{x},q\right) < \inf_{q'\in\mcal{M}} \mbf{\Upsilon}_n\left(\mbf{x},q'\right)  + 11.36 \right\}.
\end{equation}
We denote by $\hat{P}\left(n,\mbf{X},\mscr{M}\right)$ any measurable element of the closure of $\mscr{E}_n(\mbf{X})$ with respect to the Hellinger distance and we call it a $\rho$-estimator on $\mscr{M}$. The constant 11.36 is given by (7) and (19) in \cite{baraudrevisited} but can be replaced by any smaller positive number.\par
One of the main results of $\rho$-estimation is Theorem 1 in \cite{baraudrevisited}. For independent random variables $X_1,\dots,X_n$, any $\rho$-estimator  $\hat{P}=\hat{P}(n,\mbf{X},\mscr{M})$ satisfies an inequality of the form
\begin{equation}
\label{eq:intro_1}
\mbb{P}\left( \frac{C}{n} \suml_{i=1}^n h^2(P_i,\hat{P}) \leq \inf_{Q\in\mscr{M}} n^{-1} \suml_{i=1}^n h^2(P_i,Q) + \frac{D_n(\mscr{M}) + \xi}{n} \right) \geq 1-e^{-\xi},
\end{equation}
where $C$ is a positive numeric constant and $D_n(\mscr{M})\geq 1$ is a dimension term that measures the complexity of the model $\mscr{M}$. This dimension term corresponds to a bound on the $\rho$-dimension. It is an important feature of $\rho$-estimation as it determines the bound on the convergence rate of the estimator. If we actually  dispose of i.i.d. observations with common distribution $\overline{P}$ in $\mscr{M}$, we get
\[ \mbb{P}\left( C h^2(\overline{P},\hat{P}) \leq \frac{D_n(\mscr{M}) + \xi}{n} \right) \geq 1-e^{-\xi}, \]
which leads to the bound $D_n(\mscr{M})/n$ on the convergence rate, up to a multiplicative constant. The notion of $\rho$-dimension is formally introduced in the appendix (Section \ref{sec:rho_application}).
\subsection{From independent to dependent data}
\label{sec:from_independent_to_dependent}
To extend the previous result to non-independent samples, we use the following idea which is not specific to our framework. We state this basic principle in a general context. 
%
%\subsubsection{Basic principle}
%
Let $\hat{\theta}:\mscr{X}^n\rightarrow \Theta$ be an estimator of some quantity $\theta\in\Theta$. The next result is proven in Section \ref{sec:proof_lem_kldiv_exponential}.\par
\begin{lemme}
\label{lem:kldiv_exponential}
Let $l:\Theta\times\Theta\rightarrow\mbb{R}_+$ be a loss function, $\mbf{P},\mbf{Q}$ two distributions on a measurable space $\left(\mscr{Y},\mcal{X}\right)$ and $\beta\in(0,1]$. Assume that when $\mbf{Y}$ has distribution $\mbf{P}$
\begin{equation}
\label{eq:exponential_deviation_inequality_general}
\mbb{P}_{\mbf{X}\sim \mbf{P}}\left( l\left(\hat{\theta}(\mbf{X}),\theta\right) \geq A + \frac{B+\xi^{\beta}}{n} \right) \leq e^{-\xi}, \forall \xi>0,
\end{equation}
then, when $\mbf{X}$ has distribution $\mbf{Q}$
\[ \mbb{E}_{\mbf{X}\sim\mbf{Q}}\left[ l\left(\hat{\theta}(\mbf{X}),\theta\right) \right] \leq A + \frac{ B + \left(2+\frac{3}{2}\mbf{K}\left(\mbf{Q}||\mbf{P}\right) \right)^{\beta} }{n}, \]
where $\mbf{K}$ is the Kullback-Leibler divergence given by
\[ \mbf{K}(Q||P) = \begin{cases}
                    \int \log\left(\frac{dQ}{dP}\right) dQ \text{  if  } Q\ll P,\\
                    +\infty \text{  otherwise}.
                   \end{cases}
\]
\end{lemme}
Deviation inequalities for $\rho$-estimators $\widehat \theta$ have been established under the assumption that one observes independent random variables $X_{1},\ldots,X_{n}$, hence when the distribution of $\mbf{X}=(X_{1},\ldots,X_{N})$ is $\mbf{P}=\mcal{L}(X_1)\otimes\dots \otimes \mcal{L}(X_n)$. Our idea is to apply Lemma~\ref{lem:kldiv_exponential} with a distribution $\mbf{Q}\ll \mbf{P}$, which is not a product probability, in order to establish a risk bound for the estimator $\widehat \theta$ when the observations $X_{1},\ldots,X_{n}$ are possibly dependent. The quantity $ \mbf{K}(\mbf{Q}||\mbf{P})$ measures thus a departure from independence. We consider subsets of the original data $X_{1},\ldots,X_{n}$ when this quantity is too big.\par
Let $n$ be larger than 2. We build subsets of observations by taking them separated by blocks of length $s\in\mbb{N}$, as described in the diagram below.\\
\setlength{\unitlength}{3mm}
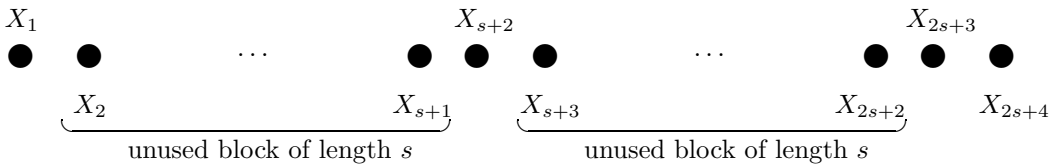
\begin{picture}(70,12)
%points
\multiput(2,6)(20,0){3}{\circle*{1}}
\multiput(5,6)(20,0){3}{\circle*{1}}
\multiput(19.5,6)(20,0){2}{\circle*{1}}
\multiput(12.3,3.2)(20,0){2}{\oval(17,1)[b]}
\multiput(6.7,1.5)(20,0){2}{{\footnotesize unused block of length $s$}}
\multiput(11.5,6)(20,0){2}{{\footnotesize $\dots$}}
%sub-sample
\put(1.3,7.3){{\footnotesize $X_1$}}
\put(21,7.3){{\footnotesize $X_{s+2}$}}
\put(40.8,7.3){{\footnotesize $X_{2s+3}$}}
%unused
\put(4.3,3.5){{\footnotesize $X_2$}}
\put(18.3,3.5){{\footnotesize $X_{s+1}$}}
\put(23.9,3.5){{\footnotesize $X_{s+3}$}}
\put(37.7,3.5){{\footnotesize $X_{2s+2}$}}
\put(44,3.5){{\footnotesize $X_{2s+4}$}}
\end{picture}
Formally, for $s\in\{0,1,\dots,s_{\max}\},s_{\max}:= \lfloor (n-2)/2 \rfloor$ and $b\in[s+1]$, we define 
\[ n(s,b) := \left\lfloor \frac{n+s+1-b}{1+s} \right\rfloor\geq 2, \]
for $i\in[n(s,b)]$
\begin{equation}
\label{eq:X_s}
X^{(s,b)}_i := X_{b+(i-1)(s+1)}\in\mscr{X},\forall i\in[n(s,b)],
\end{equation}
and
\[ \mbf{X}^{(s,b)} := \left( X^{(s,b)}_i, i\in[n(s,b)] \right). \]
We obtain $s+1$ subsets $\mbf{X}^{(s,1)},\ldots,\mbf{X}^{(s,s+1)}$ with sizes $n(s,1),\ldots,n(s,s+1)$ respectively. For each block $b\in[s+1]$, we consider the probabilities $\mbf{P}^*_{s,b}$ and $\mbf{P}^{ind}_{s,b}$ which are defined by 
\begin{equation}
\label{eq:p_s_ind}
\mbf{P}^*_{s,b} := \mcal{L}\left(\mbf{X}^{(s,b)} \right) \text{  and  }\mbf{P}^{ind}_{s,b}:=\bigotimes_{i=1}^{n(s,b)} \mcal{L}\left(X^{(s,b)}_i\right).
\end{equation}
We denote for short $\mbf{P}^*:=\mbf{P}^*_{0,1}$ the distribution of $\mbf{X}=\left(X_1,\ldots,X_n\right)$ and $\mbf{P}^{ind}:=\mbf{P}^{ind}_{0,1}=\mcal{L}(X_1) \otimes \dots \otimes \mcal{L}(X_n)$. Our estimator is obtained with the following the statistical procedure.  
\begin{enumerate}
\item Let $s$ be in $\{0,1,\ldots,s_{\max}\}$. For $b$ in $[s+1]$, we denote by $\hat{P}_{s,b}$ the estimators given by
\[ \hat{P}_{s,b} := \hat{P}\left(n(s,b),\mbf{X}^{(s,b)},\mscr{M}\right),\]
where the $\rho$-estimator $\hat{P}\left(n(s,b),\mbf{X}^{(s,b)},\mscr{M}\right)$ is defined in Section \ref{sec:rho_estimation}.
\item We denote by $\hat{P}_s=\hat{P}_s\left(\mbf{X},\mscr{M}\right)$ any element of $\mscr{M}$ that satisfies
\begin{equation}
\label{eq:p_hat_s}
\suml_{b=1}^{s+1} n(s,b) h^2\left(\hat{P}_{s,b},\hat{P}_s\right) \leq\inf_{Q\in\mscr{M}} \suml_{b=1}^{s+1} n(s,b) h^2\left(\hat{P}_{s,b},Q\right) + \iota,
\end{equation}
where $\iota$ is any fixed constant in $(0,1273]$.
\end{enumerate}
\subsection{Main result}
We assume that the $\rho$-dimension function (see Section \ref{sec:rho_application}) is uniformly bounded by a function $m\mapsto D_m(\mscr{M}) \geq 1$ which is non-decreasing. 
\begin{theoreme}
\label{th:main_rho}
For any random variables $X_1,\dots,X_n$ on $(\mscr{X},\mcal{X})$, the estimator $\hat{P}_s=\hat{P}_s\left(\mbf{X},\mscr{M}\right)$ given by (\ref{eq:p_hat_s}) satisfies
\begin{align}
\mbb{E}_{\mbf{P}^*} \left[ n^{-1} \suml_{i=1}^n h^2\left( P_i, \hat{P}_s \right) \right] &\leq \frac{c_0}{n} \inf_{Q\in\mscr{M}} \suml_{i=1}^n h^2\left( P_i, Q \right)\label{eq:th_1}\\
+ c_1 \frac{(s+1)}{n} & \left[17+D_{n(s,1)}(\mscr{M}) \right] + \frac{c_2}{n} \suml_{b=1}^{s+1} \mbf{K}\left(\mbf{P}^*_{s,b}||\mbf{P}^{ind}_{s,b}\right),\nonumber
\end{align}
where $c_0=602$, $c_1=20 056/4.7$ and $c_2=30 084$. 
\end{theoreme}
The proof of this result is postponed to Section \ref{sec:proof_th_main_rho}. One can check that we do not need any assumption on the data to obtain this result. We only need a condition on the model $\mscr{M}$ which is chosen by the statistician. However a posteriori assumptions are necessary to make this bound meaningful. It follows from the triangle inequality and $(a + b)^2 \leq 2 a^2 + 2 b^2$ for all non-negative numbers $a$ and $b$ that for any $\overline{P} \in \mscr{M}$,
\[ n h^2\left(\overline{P},\hat{P}_s\right) \leq 2 \suml_{i=1}^n h^2\left(P_i,\hat{P}_s\right) + 2 \suml_{i=1}^n h^2\left(P_i,\overline{P}\right).  \]
We derive from (\ref{eq:th_1}) the following
\begin{align}
\label{eq:th_1_comment}
C \mbb{E}_{\mbf{P}^*} \left[ h^2\left( \overline{P}, \hat{P}_s \right) \right] &\leq \frac{ (s+1)D_{n(s,1)}(\mscr{M})}{n} + n^{-1} \suml_{i=1}^n h^2\left( P_i, \overline{P} \right)\\
&+ n^{-1} \suml_{b=1}^{s+1} \mbf{K}\left(\mbf{P}^*_{s,b}||\mbf{P}^{ind}_{s,b}\right),\nonumber
\end{align}
where $C$ is a universal positive constant. Up to the factor $(s+1)$, the first term in the right-hand side of this inequality corresponds to the bound we would get if the data were truly i.i.d.\ with distribution $\overline P\in \mscr{M}$. In this ideal situation, both the second and third term vanish. When the data are not identically distributed, the second term is not zero but its size remains small when most of the true marginal distributions $P_{1},\ldots,P_{n}$ lie close enough to an element $\overline P\in \mscr{M}$. The third term accounts for the fact that the data are possibly dependent. We expect that for a choice of $s$ that is sufficiently large the observations 
\[
\mbf{X}^{(s,b)} :=\left(X_{b},X_{b+(s+1)},\ldots, X_{b+n(s,b)(s+1)}\right) \quad\text{with $b\in [s+1]$}
\]
be nearly independent and consequently that the quantity $n^{-1} \suml_{b=1}^{s+1} \mbf{K}\left(\mbf{P}^*_{s,b}||\mbf{P}^{ind}_{s,b}\right)$ be small compared to the first term.
%
%In the subsequent sections, we shall provide a control of this quantity under \emph{a posteriori} assumptions on the process that generates the data. We shall also see how we can obtain a data-driven choice of the parameter $s$ when we can observe of an independent copy of $X_{1},\ldots,X_{n}$.\par
%
%For now, we can consider a very simple case, e.g. when $\mbf{X}$ is a $m$-dependent process. In that case, for $s\geq m$ we have $\mbf{K}\left(\mbf{P}^*_{s,b}||\mbf{P}^{ind}_{s,b}\right)=0$ for all $b\in[s+1]$, and therefore we have
%\[ C \mbb{E}_{\mbf{P}^*} \left[ h^2\left( \overline{P}, \hat{P}_s \right) \right] \leq \frac{(s+1)D_{n(s,1)}(\mscr{M})}{n}, \]
%which is similar to what we had for i.i.d. observations for $s=0$, up to a multiplicative constant.
%
\subsection{Robust properties of our estimator}
\label{sec:robust}
The robustness properties of $\rho$-estimators in the independent context are illustrated in Section 5 \cite{baraudrevisited}. Let $\mbf{X}=\left(X_1,\dots,X_n\right)$ be the true process of interest such that $\mcal{L}(X_i)=\overline{P}$ for all $i$ in $[n]$. We actually observe a contaminated version of it. Let $Z_1,\dots,Z_n$ be random variables with any distributions. Let  $E_1,\dots,E_n$ be Bernoulli random variables such that
\begin{equation}
\label{eq:x_contamination}
Y_i = E_i X_i + (1-E_i) Z_i,\forall i\in[n].
\end{equation}
The next result shows that the mixing regime is not altered by independent contamination/outliers. It is proven in Section \ref{sec:proof_lem_robust_kl}.
\begin{lemme}
\label{lem:robust_kl}
If $E_1,\dots,E_n,Z_1,\dots,Z_n$ and $\mbf{X}$ are mutually independent, we have 
\[ \mbf{K}\left( \mcal{L}\left(\mbf{Y}\right) || \mcal{L}(Y_1) \otimes \dots \otimes \mcal{L}(Y_n) \right) \leq \mbf{K}\left( \mcal{L}\left(\mbf{X}\right) || \mcal{L}(X_1) \otimes \dots \otimes \mcal{L}(X_n) \right). \]
\end{lemme}
We can deduce a corollary of Theorem \ref{th:main_rho} from this. We define $p_i$ by $\mbb{P}(E_i=1)=p_i$ for $i\in[n]$.
\begin{coro}
\label{coro:robust}
Let $\hat{P}_s=\hat{P}_s\left(\mbf{Y},\mscr{M}\right)$ be the estimator given by (\ref{eq:p_hat_s}). There is a positive universal constant $C$ such that in the situation of Lemma \ref{lem:robust_kl}, we have
\begin{align*}
C \mbb{E} \left[ h^2\left( \overline{P}, \hat{P}_s \right) \right] &\leq h^2\left(\overline{P},\mscr{M}\right) + n^{-1} \suml_{i=1}^n (1-p_i)\\
&+ \frac{(s+1)D_{n(s,1)}(\mscr{M})}{n} + n^{-1} \suml_{b=1}^{s+1} \mbf{K}\left(\mbf{P}^*_{s,b}||\mbf{P}^{ind}_{s,b}\right),
\end{align*}
where $\mbf{P}^*_{s,b}$ and $\mbf{P}^{ind}_{s,b}$ are given by (\ref{eq:p_s_ind}).
\end{coro}
This result is proven in Section \ref{sec:proof_coro_robust}. Inspired by H\"{u}ber’s contamination model, we consider the situation $\overline{P}\in\mscr{M}$ and $p_i=1-\epsilon_{cont}$ for all $i\in[n]$. We get
\[ C \mbb{E}\left[ h^2\left(\overline{P},\hat{P}_s\right)\right] \leq \epsilon_{cont} + \frac{(s+1) D_{n(s,1)}(\mscr{M})}{n} + n^{-1} \suml_{b=1}^{s+1} \mbf{K}\left(\mbf{P}^*_{s,b}||\mbf{P}^{ind}_{s,b}\right). \]
Our bound on the convergence rate is not deteriorated as long as the contamination rate $\epsilon_{cont}$ is small compared to the other terms. Equally, we can consider the case where the $E_i$ are deterministic, i.e. there is a subset $I\subset [n]$ such that $\mbb{P}(E_i = 0) = \mathbbm{1}_{i\in I}$. We get
\[ C\mbb{E}\left[ h^2\left(\overline{P},\hat{P}_s\right) \right] \leq \frac{|I|}{n} + \frac{(s+1) D_{n(s,1)}(\mscr{M})}{n} + n^{-1}\suml_{b=1}^{s+1} \mbf{K}\left(\mbf{P}^*_{s,b}||\mbf{P}^{ind}_{s,b}\right). \]
As before, our bound on the convergence rate is not deteriorated as long as the proportion of outliers $|I|/n$ is small compared to the other terms on the right hand side.
\subsection{The particular case of Markov chains}
\label{sec:mixing_markov}
Under the assumption that $X_1,\dots,X_n$ is a Markov chain, the quantity $\mbf{K}\left(\mbf{P}_{s,b}^*||\mbf{P}_{ind}^{(s,b)}\right)$ can be written in a form given in the lemma below.
\begin{lemme}
\label{lem:divergence_markov}
If $\mbf{X}$ is a Markov chain,
\[ \mbf{K}\left( \mcal{L}\left(\mbf{X}\right) || \mcal{L}\left( X_1 \right) \otimes \dots \otimes \mcal{L}\left( X_n \right) \right) = \suml_{i=2}^{n} I(\sigma(X_i),\sigma(X_{i+1})), \]
where
\begin{equation}
\label{eq:coefficient_of_information}
I(\sigma(X_i),\sigma(X_{i+1})) := \mbf{K}\left( \mcal{L}(X_i,X_{i+1}) || \mcal{L}(X_i) \otimes \mcal{L}(X_{i+1}) \right).
\end{equation}
In particular for all $s$ in $\{0,1,\ldots,s_{\max}\}$ and all $b$ in $[s+1]$,
\[ \mbf{K}\left(\mbf{P}_{s,b}^*||\mbf{P}_{ind}^{(s,b)}\right) = \suml_{i=2}^{n(s,b)} I\left(\sigma(X^{(s,b)}_i),\sigma(X^{(s,b)}_{i+1}) \right), \]
where the $X^{(s,b)}_i$ are given by (\ref{eq:X_s}).
\end{lemme}
This result is proven in Section \ref{sec:proof_lem_divergence_markov}. It tells us that for Markov chains we only need to consider the simpler quantities $I(\sigma(X_i),\sigma(X_{i+s+1}))$ referred to as \emph{coefficient of information} by Bradley \cite{bradley}. This result also extends to hidden Markov models.
\begin{lemme}
\label{lem:divergence_hmm}
If $\left(X_i,H_i\right)_{1\leq i\leq n}$ is a HMM, we have
\[ \mbf{K}\left( \mcal{L}\left(\mbf{X}\right) || \mcal{L}\left( X_1 \right) \otimes \dots \otimes \mcal{L}\left( X_n \right) \right) \leq \suml_{i=2}^n I(\sigma(H_{i-1}),\sigma(H_i)). \]
In particular for all $s$ in $\{0,1,\ldots,s_{\max}\}$ and all $b$ in $[s+1]$,
\[ \mbf{K}\left(\mbf{P}_{s,b}^*||\mbf{P}_{ind}^{(s,b)}\right) \leq \suml_{i=1}^{n(s,b)-1} I\left(\sigma(H_{b+(i-1)(s+1)}),\sigma(H_{b+i(s+1)}) \right). \]
\end{lemme}
The proof of this result is postponed to Section \ref{sec:proof_lem_divergence_hmm}. This means that for HMMs we only need to consider the coefficients of information of the hidden chain. In what follows we consider different processes for which the coefficient of information has an exponential decay. In that case there exist positive constants $C$ and $r$ such that  
\[ n^{-1} \suml_{b=1}^{s+1} \mbf{K}\left( \mbf{P}^*_{s,b} || \mbf{P}^{ind}_{s,b} \right) \leq C e^{-r s}, \]
for all $s$ in $\{0,1,\ldots,s_{\max}\}$. For $s\geq r^{-1} \log n$ the quantity $n^{-1} \suml_{b=1}^{s+1} \mbf{K}\left( \mbf{P}^*_{s,b} || \mbf{P}^{ind}_{s,b} \right)$ is small compared to the first term on the right hand side in (\ref{eq:th_1_comment}), as it cannot be of order smaller than $1/n$. Such a constant $r$ is usually not known in practice but taking $s$ of order $\log^2 n$ ensures that for $n$ large enough the quantity we consider remains small compared to the term $(s+1) D_{n(s,1)}(\mscr{M})/n$. We pay the price of not knowing the constant $r$ with a worse logarithmic term in the latter quantity.
\section{Estimation of the invariant distribution of a diffusion process}
\label{sec:exponentially_ergodic_estimation}
We consider some diffusion processes that have been investigated by Royer \cite{Royer} and use the same vocabulary that they introduced.
\subsection{Langevin equation}
Let $d$ be a positive integer and $U:\mbb{R}^d\rightarrow \mbb{R}$ be a function of class $\mcal{C}^2$. The Langevin equation is the following stochastic differential equation
\begin{equation}
\label{eq:kolmogorov_equation}
dY_t = dB_t - \nabla U(Y_t) dt,
\end{equation}
where $B=(B_t)_{t\geq 0}$ is a $d$-dimensional Brownian motion. Its solution are called Kolmogorov processes in Royer \cite{Royer}. We assume that $U$ satisfies the following.
\begin{hypothese}
\label{hyp:kolmogorov}
The function $U$ is convex on $\mbb{R}^d$ and there exists a positive constant $\underline{\lambda}(U)$ such that the smallest eigenvalue of the Hessian matrix $U''(x)$ at $x\in \mbb{R}^d$ is not smaller than $\underline{\lambda}(U)$ for all $x$ in $\mbb{R}^d$. Besides we have
\begin{equation}
\label{eq:hyp_kolmogorov}
\inf_{x\in\mbb{R}^d} \left\{ \left| \left| \nabla U(x) \right| \right|_2^2 - \textup{Tr}\left( U''(x) \right) \right\} > -\infty,
\end{equation}
where \textup{Tr}$(A)$ is the trace of the matrix $A$.
\end{hypothese}
Under our assumption on the eigenvalues of $U''$, $\int_{\mbb{R}^d} e^{-\alpha U(x)} dx$ is finite for all $\alpha>0$ and we may define the probability measure $\overline{P}$ with density $\overline{p}$ with respect to the Lebesgue measure on  $\mbb{R}^d$ given by
\begin{equation}
\label{eq:invariant_kolmogorov}
\overline{p}(x) = Z^{-1} \exp(-2U(x)) \text{  with  } Z = \int_{\mbb{R}^d} e^{-2U(x)} dx.
\end{equation}
The probability $\overline{P}$ is the invariant probability distribution with respect to the semi-group associated to the Langevin equation (see Lemma 2.2.23 \cite{Royer}).
\begin{lemme}
\label{lem:kolmogorov_mixing}
Let $(Y_t)_{t\geq 0}$ be a stationary solution of the Langevin equation associated to a convex function $U$ that satisfies Assumption \ref{hyp:kolmogorov}. For all $s_0>0$, there exists a positive constant $C(U,s_0)$ such that for all $t>0$ and $s\geq s_0$, we have
\[ I( \sigma(Y_t), \sigma(Y_{t+s}) ) \leq C(U,s_0) \exp(-2 \underline{\lambda}(U) s). \]
\end{lemme}
This result is proven in Section \ref{sec:proof_lem_kolmogorov_mixing}. We aim to estimate $\overline{P}$ from discrete observations of a stationary Kolmogorov process.
\subsection{The framework}
\label{sec:kolmogorov_framework}
We consider the following statistical model for the observations $X_1,X_2,\dots,X_n$. For all $i\in [n]$, $X_i=Y_{t_i}$ where $\mbf{Y}=(Y_t)_{t\geq 0}$ is a stationary solution of the Langevin equation (\ref{eq:kolmogorov_equation}) for some unknown convex function $U$ that satisfies Assumption \ref{hyp:kolmogorov} and $t_{i+1}=t_i+\Delta_t$ for all $i\in [n-1]$. As a consequence of (\ref{eq:invariant_kolmogorov}), the $X_{i}$ are distributed according to the invariant measure $\overline{P}$ which has a $\log$-concave density $\overline{p}:x\mapsto Z^{-1}\exp(-2U(x))$ with respect to the Lebesgue measure. We therefore consider the set of distributions that admit a $\log$-concave density on $\mbb{R}^d$ with respect to the Lebesgue measure. As usual, this describes our statistical model but we do not want to assume that it perfectly describes reality. In the following section we recall some results about the problem of estimating a $\log$-concave density from i.i.d.\ observations.
%
%Following Lemma \ref{lem:kolmogorov_mixing}, in that situation we have a constant $C_{\Delta} \geq 0$  such that
%\begin{equation}
%\label{eq:hyp_kolmogorov}
%\mbf{K}\left(\mbf{P}^*_{s,b}||\mbf{P}^{ind}_{s,b}\right) \leq C_{\Delta} n(s,b) e^{-2 m s}, 
%\end{equation} 
%for all $s$ in $\{0,\dots,\lfloor (n-2)/2\rfloor \}$ and all $b$ in $[s+1]$. Additionally, the invariant measure $\overline{P}$ given by (\ref{eq:invariant_kolmogorov}) has a density function $\overline{p}:x\mapsto Z^{-1}\exp(-2U(x))$ with respect to the Lebesgue measure which is $\log$-concave. Therefore we consider models associated to classes of $\log$-concave functions.
%
\subsection{\texorpdfstring{$\log$}{log}-concave densities}
We refer to Kim \& Samworth \cite{samworth_kim} for the problem of estimating of $\log$-concave densities from i.i.d. observations in low dimensions ($d\in[3]$). Kur \etal \cite{kur_concave} investigated the same problem in higher dimensions ($d\geq 4$). We denote by $\mcal{F}_d$ the set of upper semi-continuous, $\log$-concave probability densities with respect to the Lebesgue measure, equipped with the $\sigma$-algebra it inherits as a subset of $L_1(\mbb{R}^d)$. We denote by $\mscr{F}_d$ the associated set of probability distributions on $\mbb{R}^d$. For $f\in\mcal{F}_d$, we define
\[ \overline{x}_f:=\int_{\mbb{R}^d} x f(x) dx \in \mbb{R}^d \text{  and  } \Sigma_f := \int_{\mbb{R}^d} (x-\mu_f)(x-\mu_f)^T f(x) dx\in\mbb{R}^{d\times d}.\]
For a symmetric, positive-definite $d\times d$ matrix $\Sigma$, we denote by $\lambda_{\min}(\Sigma)$ and $\lambda_{\max}(\Sigma)$ the smallest and largest eigenvalues respectively of $\Sigma$. For $0 < \lambda_- < \lambda_+ < \infty$ and $M>0$, we define
\[ \mcal{F}_{\lambda_-,\lambda_+,M} := \left\{ f\in\mcal{F}_d; ||\overline{x}_f||\leq M, \Sigma\in\text{Sym}(\lambda_-,\lambda_+) \right\}, \]
where
\[ \text{Sym}(\lambda_-,\lambda_+)=\left\{ \Sigma \text{ covariance matrix}, \lambda_-\leq \lambda_{\min}(\Sigma)\leq \lambda_{\max}(\Sigma) \leq \lambda_+ \right\}. \]
We denote by $\mscr{F}_{\lambda_-,\lambda_+,M}$ the class of probability distributions associated to $\mcal{F}_{\lambda_-,\lambda_+,M}$.\par
Given a subset $\mscr{A}$ of a class $\mscr{P}$ of probability distributions and $\epsilon\geq0$, we say that $\mscr{A}[\epsilon]$ is an $\epsilon$-net of $\mscr{A}$ if $\mscr{A}[\epsilon] \subset \mscr{P}$ and for all $Q$ in $\mscr{A}$ there exists $R$ in $\mscr{A}[\epsilon]$ such that $h\left(Q,R\right)\leq \epsilon$. The case $\epsilon=0$ corresponds to $\mscr{A}[\epsilon]$ being dense in $\mscr{A}$. The following result is proven in Section \ref{sec:proof_lem_entropy_log_concave} and based on the work of Kim \& Samworth \cite{samworth_kim} for $d\in[3]$ and Kur \etal \cite{kur_concave} for $d\geq 4$. 
\begin{lemme}
\label{lem:entropy_log_concave}
For all positive $\epsilon$ there exists an $\epsilon$-net $\mscr{F}_{\lambda_-,\lambda_+,M}[\epsilon]$ such that
\[ | \mscr{F}_{\lambda_-,\lambda_+,M}[\epsilon]| \leq \begin{cases}
\ds{\frac{9}{\eta_1} \frac{M (\lambda_+-\lambda_-)}{\lambda_-^{3/2}} e^{\overline{K_1}\epsilon^{-1/2}}} \text{ for }d=1,\\
\ds{\frac{3^8 \pi}{\eta_2^3} \frac{M^2 (\lambda_+-\lambda_-)^2 \lambda_+}{\lambda_-^{4}} e^{\overline{K_2}\epsilon^{-1}\log_{++}^{3/2}(1/\epsilon)}} \text{ for } d=2,\\
\ds{\frac{2^7 3^{27/2} \pi^3}{\eta_3^6} \frac{M^3 (\lambda_+-\lambda_-)^3 \lambda_+^3}{\lambda_-^{15/2}} e^{\overline{K_3}\epsilon^{-2}}} \text{ for } d=3,
\end{cases} \]
where $\eta_d$ and $\overline{K}_d$ are constants given in Theorem 4 \cite{samworth_kim} that only depend on $d$, and with $\log_{++}(x)=\max(1,\log x)$. For $d\geq 4$ and all positive $\epsilon$ there exists an $\epsilon$-net $\mscr{F}_{\lambda_-,\lambda_+,M}[\epsilon]$ such that
\[ |\mscr{F}_{\lambda_-,\lambda_+,M}[\epsilon]| \leq C_d \frac{\lambda_+^{d(d-1)/2} M^d (\lambda_+-\lambda_-)^d}{ \lambda_-^{d(d+1)/2} } \exp\left( \overline{K}_d \epsilon^{-(d-1)} \log^{(d+1)(d+2)/2}(\epsilon^{-1}) \right), \]
where $\eta_d$ and $\overline{K}_d$ are constants that only depend on $d$.
\end{lemme}
%

%
%\begin{lemme}
%\label{lem:entropy_log_concave_high_dim}
%Let $d$ be larger than 3. There are positive constants $C_d$ and $K_d$ such that
%\[ |\mcal{F}_{\lambda_-,\lambda_+,M}[\epsilon]| \leq C_d \frac{\lambda_+^{d(d-1)/2} M^d (\lambda_+-\lambda_-)^d}{ \lambda_-^{d(d+1)/2} } \exp\left( K_d \epsilon^{-(d-1)} \log(\epsilon^{-1})^{(d+1)(d+2)/2} \right), \]
%for all $\epsilon$ in $(0,1)$.
%\end{lemme}
%
%The proof of this result can be found in Section \ref{sec:proof_lem_entropy_log_concave_high_dim}. 
%
\subsubsection{The case \texorpdfstring{$d\in\{1,2,3\}$}{dleq3}}
Let $\mscr{F}_{\lambda_-,\lambda_+,M}[\epsilon]$ be a $\epsilon$-net of $\mscr{F}_{\lambda_-,\lambda_+,M}$ that satisfies the bound given in Lemma \ref{lem:entropy_log_concave} for
\begin{equation}
\label{eq:m_log_concave}
\lambda_+ = \lambda_-^{-1} = M :=  \begin{cases}
\ds{ \exp\left( \overline{K}_1 \left( n / \log n \right)^{1/5} \right) } \text{  for  } d=1,\\
\ds{ \exp\left( \overline{K}_2 n^{1/3} \log^{2/3} n  \right) } \text{  for  } d=2,\\
\ds{ \exp\left( \overline{K}_3 \left( n / \log n \right)^{1/2} \right) } \text{  for  } d=3,
                  \end{cases}
\end{equation}
and
\begin{equation}
\label{eq:delta_log_concave}
\epsilon := \begin{cases}
\ds{ n^{-2/5} \log^{2/5}n }\text{  for  } d=1,\\
\ds{ n^{-1/3} \log^{5/6}n }\text{  for  } d=2,\\
\ds{ n^{-1/4} \log^{1/4}n }\text{  for  } d=3.
                  \end{cases}
\end{equation}
The following result holds and its proof can be found in Section \ref{sec:proof_th_kolmogorov}.
\begin{theoreme}
\label{th:kolmogorov}
Let $n\geq 3$ and $X_1,X_2,\dots,X_n$ be arbitrary random variables with marginal distributions $P_1,\dots,P_n$. The $\rho$-estimator $\hat{P}_s$ given by (\ref{eq:p_hat_s}) with $\mscr{M}=\mscr{F}_{\lambda_-,\lambda_+,M}[\epsilon]$ satisfies for all $\overline{P}\in\mscr{P}_X$
\begin{align}
C_d \mbb{E} \left[ h^2\left( \overline{P}, \hat{P}_s \right) \right] &\leq h^2\left( \overline{P}, \mscr{F}_{\lambda_-,\lambda_+,M} \right) + n^{-1} \suml_{i=1}^n h^2\left(P_i,\overline{P}\right)\label{eq:th_kolmogorov}\\
&+ n^{-1} \suml_{b=1}^{s+1} \mbf{K}\left(\mbf{P}^*_{s,b} || \mbf{P}^{ind}_{s,b}\right)\nonumber\\
&+ \begin{cases}
\ds{n^{-4/5}  \left( \log^{4/5} n + s \log^{-1/5} n \right) } \text{  for  } d=1,\\
\ds{n^{-2/3} \left( \log^{5/3} n + s \log^{2/3} n \right) } \text{  for  } d=2,\\
\ds{n^{-1/2} \left( \log^{1/2} n + s \log^{-1/2} n \right) } \text{  for  } d=3,
\end{cases}\nonumber
\end{align}
for positive constants $C_1,C_2,C_3$. In particular if the model described in Section \ref{sec:kolmogorov_framework} is exact and $s\geq (2\underline{\lambda}(U))^{-1} \log n $, there exists a positive constant $C(U,d,\Delta_t)$ such that for $n$ large enough
\[ C(U,d,\Delta_t) \mbb{E} \left[ h^2\left( \overline{P}, \hat{P}_s \right) \right] \leq \begin{cases}
\ds{n^{-4/5}  \left( \log^{4/5} n + s \log^{-1/5} n \right) } \text{  for  } d=1,\\
\ds{n^{-2/3} \left( \log^{5/3} n + s \log^{2/3} n \right) } \text{  for  } d=2,\\
\ds{n^{-1/2} \left( \log^{1/2} n + s \log^{-1/2} n \right) } \text{  for  } d=3,
\end{cases}
\]
where $\overline{P}$ is the invariant distribution given by (\ref{eq:invariant_kolmogorov}).
\end{theoreme}
Inequality (\ref{eq:th_kolmogorov}) is a consequence of Theorem \ref{th:main_rho} and does not require any assumption on the data. The last term comes from the control of the dimension of the net $\mscr{F}_{\lambda_-,\lambda_+,M}[\epsilon]$ and the choice of $\epsilon$ given by (\ref{eq:delta_log_concave}).  Ideally, most of the distributions $P_i$ lie in a small neighborhood of a distribution $\overline{P}$ in $\mscr{F}_{\lambda_-,\lambda_+,M}$ so that the first two terms in the bound remain small compared to the last term. Those two terms vanish when the model is exact and a good choice of $s$ guarantees the term $n^{-1} \suml_{b=1}^{s+1} \mbf{K}\left( \mbf{P}^*_{s,b} || \mbf{P}^{ind}_{s,b} \right)$ is negligible with respect to the last one.\par
We can derive convergence rates for the optimal choice of $s$ given $\underline{\lambda}(U)$. One can check that up to a logarithmic factor, we obtain the same rates as Theorem 5 \cite{samworth_kim} in the i.i.d.\ case. Our power of $\log n$ is even better for $d=3$. As mentioned in Section \ref{sec:mixing_markov}, the knowledge of $\underline{\lambda}(U)$ is not necessary to obtain convergence rates. We obtain slightly worse powers of $\log n$ in the convergence rates for $s$ of order $\log^2 n$. We can also derive results for i.i.d. observations from (\ref{eq:th_kolmogorov}) by taking the term $n^{-1} \suml_{b=1}^{s+1} \mbf{K}\left( \mbf{P}^*_{s,b} || \mbf{P}^{ind}_{s,b} \right)$ down to 0 which provides a result for the robust estimation of a $\log$-concave density from i.i.d. observations.\par
In order to illustrate the robustness of our estimators we consider the situation of Section \ref{sec:robust}. Let $Z_1,\dots,Z_n$ be random variables with any distributions and $E_1,\dots,E_n$ be Bernoulli random variables such that for all $i\in[n]$,
\[ X_i = E_i Y_{t_1+(i-1)\Delta_t} + (1-E_i) Z_i, \]
where $(Y_t)_t$ is a stationary solution of the Langevin equation (\ref{eq:kolmogorov_equation}) for some unknown convex function $U$ that satisfies Assumption \ref{hyp:kolmogorov}.
\begin{coro}
\label{coro:kolmorov_robust}
Let $\hat{P}_s$ be the estimator given by (\ref{eq:p_hat_s}) with $\mscr{M}=\mscr{F}_{\lambda_-,\lambda_+,M}[\epsilon]$. If $E_1,\dots,E_n,Z_1,\dots,Z_n$ and $\mbf{X}$ are mutually independent, there exists a positive constant $C(U,d,\Delta_t)$ such that for $s\geq  (2\underline{\lambda}(U))^{-1}\log n$ we have
\begin{align}
C(U,d,\Delta_t) \mbb{E} \left[ h^2\left( \overline{P}, \hat{P}_s \right) \right] &\leq n^{-1} \suml_{i=1}^n (1-p_i)\label{eq:coro_kolmogorov_robust}\\
&+ \begin{cases}
\ds{n^{-4/5}  \left( \log^{4/5} n + s \log^{-1/5} n \right) } \text{  for  } d=1,\\
\ds{n^{-2/3} \left( \log^{5/3} n + s \log^{2/3} n \right) } \text{  for  } d=2,\\
\ds{n^{-1/2} \left( \log^{1/2} n + s \log^{-1/2} n \right) } \text{  for  } d=3,
\end{cases},\nonumber
\end{align}
where $p_i=\mbb{P}(E_i=1)$ for all $i\in[n]$.
\end{coro}
One can see that our deviation bound is not significantly worse as long as the average proportion of contamination $n^{-1} \suml_{i=1}^n (1-p_i)$ remains small compared to the last term on the right hand side of (\ref{eq:coro_kolmogorov_robust}).
\subsubsection{The case \texorpdfstring{$d\geq 4$}{dgeq4}}
Let $\mscr{F}_{\lambda_-,\lambda_+,M}[\epsilon]$ be an $\epsilon$-net of $\mscr{F}_{\lambda_-,\lambda_+,M}$ that satisfies the bound given in Lemma \ref{lem:kolmogorov_mixing} with
\begin{align}
\label{eq:lambda_d_high}
\lambda_+ = \lambda_-^{-1} &= \exp\left( \frac{\epsilon^{-(d-1)} \log^{(d+1)(d+2)/2}(\epsilon^{-1})}{d^2} \right)\\
\label{eq:m_d_high}
M &= \exp\left( \frac{\epsilon^{-(d-1)} \log^{(d+1)(d+2)/2}(\epsilon^{-1})}{d} \right),
\end{align}
with
\begin{equation}
\label{eq:epsilon_d_high}
\epsilon = n^{-\frac{1}{d+1}} \log^{\frac{1}{d+1}+\frac{d+2}{2}}n.
\end{equation}
The following result holds and its proof can be found in Section \ref{sec:proof_th_kolmogorov}.
\begin{theoreme}
\label{th:kolmogorov_high}
Let $n\geq 3$ and $X_1,X_2,\ldots,X_n$ be arbitrary random variables with marginal distributions $P_1,\ldots,P_n$. The $\rho$-estimator $\hat{P}_s$ given by (\ref{eq:p_hat_s}) with $\mscr{M}=\mscr{F}_{\lambda_-,\lambda_+,M}[\epsilon]$ satisfies for all $\overline{P}\in\PX$
\begin{align*}
C_d \mbb{E}_{\mbf{P}^*} \left[ h^2\left( \overline{P}, \hat{P}_s \right) \right] &\leq h^2\left( \overline{P}, \mscr{F}_{\lambda_-,\lambda_+,M} \right) + n^{-1} \suml_{i=1}^n h^2\left(P_i,\overline{P}\right)\\
&+ n^{-1} \suml_{b=1}^{s+1} \mbf{K}\left( \mbf{P}^*_{s,b} || \mbf{P}^{ind}_{s,b} \right)\\
&+ n^{-\frac{2}{d+1}} \left( \log^{d+2+\frac{2}{d+1}} n + s \log^{d+1+\frac{2}{d+1}} n \right).
\end{align*}
In particular if the model described in Section \ref{sec:kolmogorov_framework} is exact and $s\geq (2 \underline{\lambda}(U) )^{-1} \log n$, there exists a positive constant $C(U,d,\Delta_t)$ such that for $n$ large enough
\[ C(U,d,\Delta_t) \mbb{E} \left[ h^2\left( \overline{P}, \hat{P}_s \right) \right] \leq n^{-\frac{2}{d+1}} \left( \log^{d+2+\frac{1}{d+1}}n + s \log^{d+1+\frac{2}{d+1}} n \right), \]
where $\overline{P}$ is the invariant distribution given by (\ref{eq:invariant_kolmogorov}).
\end{theoreme}
This result is equivalent to Theorem \ref{th:kolmogorov} and the comments that applied to it also apply now. Our estimator is also robust and tolerates a higher contamination rate as the convergence rate is slower. One can check that up to a logarithmic factor, we have the same rate that Kur \etal \cite{kur_concave} obtain for the estimation of $\log$-concave estimation from i.i.d. observations. We can derive a result equivalent to Corollary \ref{coro:kolmorov_robust} for $d\geq 4$. Our estimator can tolerate an average proportion of contamination of order not larger than $n^{-\frac{2}{d+1}} \log^{d+2+\frac{2}{d+1}} n$ without its performance being significantly deteriorated.
\section{Hidden Markov models}
\label{sec:hmm}
\subsection{Stationary hidden Markov models}
Let $\left(Y_i,H_i\right)_i$ be a finite state space HMM with parameters $(K^*,w^*,Q^*,F^*)$. If $w^*$ is invariant with respect to $Q^*$, then the process $\left(Y_i,H_i\right)_i$ is stationary. As explained in the introduction, we aim at estimating the different parameters through the distribution of consecutive observations. For $L\geq 2$ we define $
P_L = P_{w^*,Q^*,F^*}$ with $P_{W^*,Q^*,F^*}$ defined by (\ref{eq:p_w_q_f}), and we have $\mcal{L}(Y_i,Y_{i+1},\dots,Y_{i+L-1})=P_L$ for all $i$. We have identically distributed but dependent random variables from which we can estimate $P_L$. It is possible to relax the stationary assumption.
\begin{hypothese}
\label{hyp:hmm_ergodic}
Let $\left(Y_i,H_i\right)_i$ be a finite state space HMM with parameters $(K^*,w^*,Q^*,F^*)$ such that $Q^*$ is irreducible and aperiodic.
\end{hypothese}
In this case we do not have identically distributed observations anymore. However the distribution $\mcal{L}\left(Y_i,\ldots,Y_{i+L-1}\right)$ converges exponentially fast to the distribution
\begin{equation}
\label{eq:p_star}
P^* = P_{\pi^*,Q^*,F^*},
\end{equation}
where $\pi^*$ is the only invariant distribution with respect to $Q^*$.
%vWe also have an exponential decay for the coefficient of information (see Lemma \ref{lem:k_s_hmm}).
%
\subsection{The framework}
Let $Y_1,Y_2,\ldots,Y_N$ be random variables taking values in a measurable space $\left(\mscr{Y},\mcal{Y}\right)$. Let $L$ be in $\{2,3,\ldots,\lfloor N/2\rfloor \}$ and $n$ be the integer given by $n=N+1-L$. We define the new random variables
\begin{equation}
\label{eq:x_y_hmm}
X_i = \left(Y_i,Y_{i+1},\dots,Y_{i+L-1}\right),i=1,\dots,n,
\end{equation}
taking values in the measurable space $\left(\mscr{X},\mcal{X}\right)=\left(\mscr{Y}^L,\mcal{Y}^{\otimes L}\right)$. We follow the notation established in Section \ref{sec:first_section}.\par
We denote $\mscr{P}_Y$ the class of all probability distributions on $\left(\mscr{Y},\mcal{Y}\right)$. For $K\geq 2$ and subsets $\overline{\mscr{F}}_1,\dots,\overline{\mscr{F}}_K$ of $\mscr{P}_Y$, we denote by $\mscr{H}\left(K,\overline{\mscr{F}}_1,\dots,\overline{\mscr{F}}_K\right)$ the set of distributions defined by
\begin{equation}
\label{eq:model_h_hmm}
\mscr{H}\left(K,\overline{\mscr{F}}_1,\dots,\overline{\mscr{F}}_K\right) := \left\{ P_{w,Q,F}; \begin{array}{c}
\forall k\in [K], w\in\mcal{W}_K,\\
Q\in\mcal{T}_K, F_k \in \overline{\mscr{F}}_k 
\end{array}
\right\}\subset \mscr{P}_X,
\end{equation}
where $P_{w,Q,F}$ is given by (\ref{eq:p_w_q_f}),
\begin{align}
\label{eq:t_mcal_k}
&\mcal{T}_K = \left\{ Q\in [0,1]^{K\times K}; \suml_{j=1}^K Q_{ij} =1, \forall i \in\{1,\dots,K\} \right\},\\
&\text{  and  } \mcal{W}_K = \left\{ w\in[0,1]^K; w_1+\dots+w_K=1\right\}.
\end{align}
We call \emph{emission models} the sets $\overline{\mscr{F}}_1,\ldots,\overline{\mscr{F}}_K$. Let $\overline{\mscr{M}}$ be a  non-empty subset of $\mscr{H}\left(K,\overline{\mscr{F}}_1,\dots,\overline{\mscr{F}}_K\right)$.
\subsection{Estimation}
\label{sec:hmm_estimation_L}
%
%We need a countable set of density functions to perform the procedure described in Section \ref{sec:from_independent_to_dependent}.
Let $\nu$ be a $\sigma$-finite measure on $\left(\mscr{Y},\mcal{Y}\right)$ and we denote by $\mu$ the associated $\sigma$-finite measure on $\left(\mscr{X},\mcal{X}\right)$ given by $\mu:=\nu^{\otimes L}$. We consider emission models that satisfy the following.
\begin{hypothese}
\label{hyp:hmm_vc}
We dispose of countable sets $\mcal{F}_i,i=1,\dots,K$ of probability density functions (with respect to $\nu$) such that
\begin{enumerate}
\item for all $k$ in $[K]$, the set of distributions $\mscr{F}_i:=\left\{f\cdot\nu; f\in\mcal{F}_i \right\}$ is an $\epsilon$-net of $\overline{\mscr{F}}_i$ with respect to the Hellinger distance;
\item for any $k_1,\dots,k_L\in[K]$, the class of functions
\begin{equation*}
\mcal{F}_{k_1,\dots,k_L} = \left\{ \mbf{x}\in\mscr{Y}^L \mapsto f_1(x_1) \dots f_L(x_L) ; f_l \in \mcal{F}_{k_l}, \forall l\in[L] \right\}
\end{equation*}
is VC-subgraph with VC-index not larger than $V_{k_1,\dots,k_L}$. Then we write
\begin{equation}
\label{eq:vc_v_hmm}
\overline{V} := \suml_{1\leq k_1,\dots,k_L\leq K} V_{k_1,\dots,k_L}.
\end{equation}
\end{enumerate}
\end{hypothese}
We refer to van der Vaart \& Wellner \cite{VanDerVaart} (Section 2.6.5) and Baraud \etal \cite{baraudinventiones} (Section 8) as an introduction to VC-subgraph classes of functions. We just mention the following example. Any finite set $\mcal{F}$ of real-valued functions is VC-subgraph with VC-index $V(\mcal{F})$ that satisfies
\begin{equation}
\label{eq:vc_finite}
V(\mcal{F}) \leq 1+\log_2(|\mcal{F}|). 
\end{equation}
Therefore we can consider finite $\epsilon$-nets as we did in Section \ref{sec:exponentially_ergodic_estimation}. We also show in Section \ref{sec:exponential_family} that exponential families satisfy our assumption.\par
We consider countable approximations of $\mcal{W}_K$ and $\mcal{T}_K$ given by
\begin{equation}
\label{eq:countable_simplex}
\mcal{W}_{\delta,K}:=\mcal{W}_K\cap\left([\delta,1]\cap\mbb{Q}\right)^K \text{  and  } \mcal{T}_{\delta,K}:=\mcal{T}_K\cap\left([\delta,1]\cap\mbb{Q}\right)^{K\times K},
\end{equation}
for $0<\delta\leq 1/K$. We define $\mscr{H}_{\delta}$ by 
\begin{equation}
\label{eq:h_delta}
\mscr{H}_{\delta} := \left\{ P_{w,Q,f}; w\in\mcal{W}_{\delta,K}, Q\in\mcal{T}_{\delta,K}, f_k\in\mscr{F}_k,\forall i\in[K] \right\},
\end{equation}
where the sets $\left(\mscr{F}_k\right)_{1\leq k\leq K}$ are given in Assumption \ref{hyp:hmm_vc}. This lower bound $\delta$ is a technicality for bounding the dimension of our model. We define the countable set of distributions
\begin{equation}
\label{eq:m_delta}
\mscr{M}_{\delta} := \left\{ P_{w,Q,F} \in \mscr{H}_{\delta}; \exists P_{w',Q',F'} \in \overline{\mscr{M}}, \begin{array}{c}
h^2\left(Q_{k\cdot},Q'_{k\cdot}\right) \leq (K-1)\delta\\
h\left( F_k,F'_k \right) \leq \epsilon,\forall k\in[K],\\
h^2\left(w,w'\right) \leq (K-1)\delta,
\end{array} \right\},
\end{equation}
which is a good approximation of $\overline{\mscr{M}}$ for small values of $\delta$ and $\epsilon$. We denote by $\hat{P}_{s,\delta}$ the estimator
\begin{equation}
\hat{P}_{s,\delta} := \hat{P}_s\left( \mscr{M}_{\delta}, \mbf{X} \right), \label{eq:p_hat_s_delta}
\end{equation}
as defined by (\ref{eq:p_hat_s}). The following theorem is proven in Section \ref{sec:proof_th_hmm}.
\begin{theoreme}
\label{th:first_hmm}
Let $N\geq K+L$ and $Y_1,\ldots,Y_N$ be arbitrary random variables. Under Assumption \ref{hyp:hmm_vc}, let $\hat{P}_s=\hat{P}_{s,\delta}$ be the estimator given by (\ref{eq:p_hat_s_delta}) with
\begin{equation}
\label{eq:delta_s}
\delta = \frac{\overline{V}}{n(s,1) (K-1)}\wedge\frac{1}{K}.
\end{equation}
There exists a positive constant $C$ such that for all $\overline{P}\in\PX,$
\begin{align}
C \mbb{E}\left[ h^2\left( \overline{P}, \hat{P}_s \right) \right] &\leq h^2\left( \overline{P}, \overline{\mscr{M}} \right) + n^{-1} \suml_{i=1}^n h^2\left( \overline{P}, P_i \right) + n^{-1} \suml_{b=1}^{s+1} \mbf{K}\left(\mbf{P}^*_{s,b}||\mbf{P}^{ind}_{s,b}\right)\nonumber\\
&+ L \epsilon^2 + (s+1) L \overline{V} \frac{\log n}{n}.\label{eq:th_first_hmm}
\end{align}
In particular under Assumption \ref{hyp:hmm_ergodic}, there exist positive constants $C(Q^*)$ and $c(Q^*)$ such that for $s \geq c(Q^*) \log n  \vee (L-1)$ we have
\begin{equation}
\label{eq:bound_hmm}
C(Q^*) \mbb{E}\left[ h^2\left( P^*, \hat{P}_s \right) \right] \leq h^2\left( P^*, \overline{\mscr{M}} \right) + L \epsilon^2 + L \overline{V} \frac{s \log n}{n},
\end{equation}
where $P^*$ is given by (\ref{eq:p_star}).
\end{theoreme}
Inequality (\ref{eq:th_first_hmm}) is a consequence of Theorem \ref{th:main_rho} and does not require any assumption on the data. The last two terms come respectively from the approximation of $\overline{\mscr{M}}$ by $\mscr{M}$ and the control of the dimension of $\mscr{M}$. Ideally, we can take $\overline{P}$ in $\overline{\mscr{M}}$ such that most of the distributions $P_i$ lie in a small neighborhood of $\overline{P}$ so that the first two terms in the bound remain small compared to the last term. Under Assumption \ref{hyp:hmm_ergodic} the quantity $\sum_{i=1}^n h^2(P^*,P_i)$ is bounded and a good choice of $s$ guarantees the term $n^{-1} \sum_{b=1}^{s+1} \mbf{K}( \mbf{P}^*_{s,b} || \mbf{P}^{ind}_{s,b} )$ to be negligible with respect to the last one. The optimal choice of $s$ depends on a constant $c(Q^*)$ which relates to the spectral gap of $Q^*$. We distinguish two cases in order to obtain convergence rates over the class
\begin{align}
\label{eq:h_star}
\mscr{H}^* &\left(K,\overline{\mscr{F}}_1,\ldots,\overline{\mscr{F}}_K\right)\\
&:= \left\{ P_{w,Q,F}\in\mscr{H}\left(K,\overline{\mscr{F}}_1,\ldots,\overline{\mscr{F}}_K\right); \begin{array}{c}
Q \text{ irreducible },\\
Q \text{ aperiodic},\\                                                                                                                                                                                
\text{ and } w=Qw                                                                                                                                                                                \end{array}
\right\}.\nonumber
\end{align}
The first case is when we satisfy Assumption \ref{hyp:hmm_vc} with $\epsilon=0$. In that situation and for $P^*$ in $\overline{\mscr{M}}=\mscr{H}\left(K,\overline{\mscr{F}}_1,\dots,\overline{\mscr{F}}_K\right)$ the first two terms in (\ref{eq:bound_hmm}) vanish. For the optimal choice of $s$ our estimator achieves the convergence rate $n^{-1} \log^2 n$ with respect to the squared Hellinger distance over $\mscr{H}^*\left(K,\overline{\mscr{F}}_1,\dots,\overline{\mscr{F}}_K\right)$. This means that up to a logarithmic term we achieve the optimal rate $1/n$ in the independent context (see Birgé \cite{birge_minimax_hellinger}). As mentioned in Section \ref{sec:mixing_markov}, the knowledge of $c(Q^*)$ is not necessary to obtain convergence rates. We only obtain slightly worse powers of $\log n$ in the convergence rates for $s=\log^2 n$.\\
The second case is when we cannot take $\epsilon= 0$. In that situation the term $\overline{V}$ depends on $\epsilon$ and we proceed as in Section \ref{sec:exponentially_ergodic_estimation}. We obtain a convergence rate taking $\epsilon$ that goes to 0 with $n$ at a rate that balances the last two terms in (\ref{eq:bound_hmm}). This happens when $\epsilon^2/\overline{V}$ is of order $n^{-1}$ up to a logarithmic term. We put it in application in Section \ref{sec:hmm_log_concave}.\par
In order to illustrate the robustness of our estimators we consider the situation of Section \ref{sec:robust}. Let $Z_1,\dots,Z_N$ be random variables with any distributions and $E_1,\dots,E_N$ be Bernoulli random variables such that for all $i\in[N]$,
\[ Y_i = E_i Y'_i + (1-E_i) Z_i, \]
where $\mbf{Y}'$ satisfy Assumption \ref{hyp:hmm_ergodic}. The following result is proven in Section \ref{sec:proof_coro_robust_hmm}.
\begin{coro}
\label{coro:hmm_robust}
Let $N\geq K+L$ and $\hat{P}_s=\hat{P}_{s,\delta}$ be the estimator given by (\ref{eq:p_hat_s_delta}) with $\delta$ given by (\ref{eq:delta_s}). If $E_1,\dots,E_N,Z_1,\dots,Z_N$ and $\mbf{Y}'$ are mutually independent, there exist positive constants $C(Q^*)$ and $c(Q^*)$ such that for $s \geq  c(Q^*) \log n$ we have
\begin{align}
C(Q^*) \mbb{E} \left[ h^2\left( P^*, \hat{P}_s \right) \right] &\leq h^2 \left( P^*,\overline{\mscr{M}} \right) + \frac{L}{N} \suml_{i=1}^N (1-p_i)\label{eq:coro_hmm_robust}\\
&+ L\epsilon^2 + L \overline{V} \frac{s \log n}{n},\nonumber
\end{align}
where $p_i=\mbb{P}(E_i=1)$ for all $i\in[N]$ and $\delta$ is given by (\ref{eq:delta_s}).
\end{coro}
One can see that our deviation bound is not significantly worse as long as the average proportion of contamination $\frac{L}{N} \suml_{i=1}^N (1-p_i)$ remains small compared to the last two terms. One would typically look at the following situation. We assume that the model is well specified, i.e. $P^*\in\overline{\mscr{M}}$. For H\"{u}ber's contamination model, i.e.\ $p_i=1-\alpha_{cont}$ for all $i\in[N]$, we get
\begin{equation}
\label{eq:hmm_contamination_bound}
C(Q^*) \mbb{E} \left[ h^2\left( P^*, \hat{P}_s \right) \right] \leq L \left[ \alpha_{cont} + \epsilon^2 + \overline{V} \frac{s\log n}{n} \right], 
\end{equation}
for $s\geq c(Q^*) \log n$.
The bound on the convergence rate is not deteriorated as long as the contamination rate $\alpha_{cont}$ is small compared to $\epsilon^2 + \overline{V} \frac{s\log n}{n}$. We can also consider the situation where $\mbb{P}\left(E_i=0\right)=\mathbbm{1}_{i\in I}$ for some subset $I\subset[N]$. We get
\begin{equation}
\label{eq:hmm_outliers_bound}
C(Q^*) \mbb{E} \left[ h^2\left( P^*, \hat{P}_s \right) \right] \leq L\left[ \frac{|I|}{N} + \epsilon^2 + \overline{V} \frac{s\log n}{n} \right], 
\end{equation}
for $s\geq c(Q^*) \log n$. As before, our bound on the convergence rate is not deteriorated as long as the proportion of outliers $|I|/N$ is small compared to $\epsilon^2 + \overline{V} \frac{s\log n}{n}$. 
\subsubsection{\texorpdfstring{$\log$}{log}-concave emission densities}
\label{sec:hmm_log_concave}
We use results and notation given in Section \ref{sec:exponentially_ergodic_estimation}. Let $d$ be a positive integer and $\epsilon\in(0,1)$. Let $\mscr{F}_{\lambda_-,\lambda_+,M}[\epsilon]$ be an $\epsilon$-net of $\mscr{F}_{\lambda_-,\lambda_+,M}$ that satisfies the bound given in Lemma \ref{lem:kolmogorov_mixing}. We take $\overline{\mscr{F}}_k=\mscr{F}_{\lambda_-,\lambda_+,M}$ for all $k\in[K]$ and satisfy Assumption \ref{hyp:hmm_vc} with
\begin{equation}
\label{eq:v_hmm_log_concave}
\overline{V} = K^L \left( 1 + L \log_2\left( |\mcal{F}_{\lambda_-,\lambda_+,M}[\epsilon]|\right) \right).
\end{equation}
We take $\overline{\mscr{M}}=\mscr{H}\left(K,\mscr{F}_{\lambda_-,\lambda_+,M},\dots,\mscr{F}_{\lambda_-,\lambda_+,M}\right)$. We distinguish the two cases $d\in\{1,2,3\}$ and $d\geq 4$.\par
For $d\in\{1,2,3\}$ we take $\lambda_+,\lambda_-,M$ as in (\ref{eq:m_log_concave}) and $\epsilon$ as in (\ref{eq:delta_log_concave}). The following result holds and its proof can be found in Section \ref{sec:proof_th_hmm_log_concave}.
\begin{theoreme}
\label{th:hmm_log_concave}
Let $N\geq K+L$ and $\hat{P}_s=\hat{P}_{s,\delta}$ be the estimator given by (\ref{eq:p_hat_s_delta}) with $\delta$ given by (\ref{eq:delta_s}). There exist positive constants $C_1,C_2,C_3$ such that for all $\overline{P}\in\PX$,
\begin{align}
C_d \mbb{E}\left[ h^2\left( \overline{P}, \hat{P}_s \right) \right] &\leq h^2\left( \overline{P}, \overline{\mscr{M}} \right) + n^{-1} \suml_{i=1}^n h^2\left(P_i,\overline{P}\right)\label{eq:th_hmm_log_concave_1}\\
&+ n^{-1} \suml_{b=1}^{s+1} \mbf{K}\left( \mbf{P}^*_{s,b} || \mbf{P}^{ind}_{s,b} \right)\nonumber\\
&+ (s+1) L^2 K^L  \times \begin{cases}
\ds{n^{-4/5} \log^{4/5} n} \text{  for } d=1,\\
\ds{n^{-2/3} \log^{5/3} n} \text{ for } d=2,\\
\ds{n^{-1/2} \log^{1/2} n} \text{ for }d=3.
\end{cases}\nonumber
\end{align}
In particular under Assumption \ref{hyp:hmm_ergodic}, there exist positive constants $C(Q^*)$ and $c(Q^*)$ such that for $s\geq c(Q^*) \log n$ we have
\begin{equation*}
C(Q^*) \mbb{E} \left[ h^2\left( P^*, \hat{P}_s \right) \right] \leq h^2\left( P^*, \overline{\mscr{M}} \right) + s L^2 K^L  \times \begin{cases}
\ds{n^{-4/5} \log^{4/5} n} \text{  for } d=1,\\
\ds{n^{-2/3} \log^{5/3} n} \text{ for } d=2,\\
\ds{n^{-1/2} \log^{1/2} n} \text{ for }d=3,
\end{cases}
\end{equation*}
where $P^*$ is given by (\ref{eq:p_star}).
\end{theoreme}
Inequality (\ref{eq:th_hmm_log_concave_1}) is a consequence of Theorem \ref{th:first_hmm} and does not require any assumption on the data. We can deduce convergence rates over the class $\mscr{H}^*\left(K,\mscr{F}_d,\ldots,\mscr{F}_d\right)$, where $\mscr{F}_d$ is the set of distributions with $\log$-concave densities defined in Section \ref{sec:exponentially_ergodic_estimation}. For the optimal choice of $s$, we have
\begin{equation}
\label{eq:rates_hmm_log_concave}
C(Q^*) \mbb{E}\left[ h^2\left( P^*, \hat{P}_s \right) \right]  \leq L^2 K^L \times \begin{cases}
\ds{n^{-4/5} \log^{9/5} n} \text{  for } d=1,\\
\ds{n^{-2/3} \log^{8/3} n} \text{ for } d=2,\\
\ds{n^{-1/2} \log^{3/2} n} \text{ for }d=3,
\end{cases}
\end{equation} 
for all $P^*$ in $\mscr{H}^*\left(K,\mscr{F}_d,\ldots,\mscr{F}_d\right)$. We see that we have a worse power of $\log n$ compared to Theorem \ref{th:kolmogorov}. It comes from an additional logarithmic factor in the dimension term for HMMs. Corollary \ref{coro:hmm_robust} tells us our estimator is also robust to contamination and outliers. Let us illustrate it for $d=1$. We can see from (\ref{eq:hmm_contamination_bound}) that our bound is not significantly worse as long as the contamination rate $\alpha_{cont}$ is of order not larger than $n^{-4/5} \log^{9/5} n$. Similarly (\ref{eq:hmm_outliers_bound}) tells us that a number $|I|$ of outliers of order not larger than $n^{1/5} \log^{9/5} n$ does not significantly deteriorate our bound on the convergence rate of our estimator. We can follow the same train of thought for $d=2$ and $d=3$ and deduce the level of contamination or outliers our estimator can tolerate before its performance significantly worsens.\par
For $d\geq 4$ we take $\lambda_+,\lambda_-^{-1}$ as in 
(\ref{eq:lambda_d_high}), $M$ as in (\ref{eq:m_d_high}) and $\epsilon$ as in (\ref{eq:epsilon_d_high}). The following result holds and its proof can be found in Section \ref{sec:proof_th_hmm_log_concave}.
\begin{theoreme}
\label{th:hmm_log_concave_high}
Let $N\geq K+L$ and $\hat{P}_s=\hat{P}_{s,\delta}$ be the estimator given by (\ref{eq:p_hat_s_delta}) with $\delta$ given by (\ref{eq:delta_s}). There exist a positive constant $C_d$ such that for all $\overline{P}\in\PX$,
\begin{align*}
C_d \mbb{E}\left[ h^2\left( \overline{P}, \hat{P}_s \right) \right] &\leq h^2\left( \overline{P}, \overline{\mscr{M}} \right) + n^{-1} \suml_{i=1}^n h^2\left(P_i,\overline{P}\right)\\
&+ n^{-1} \suml_{b=1}^{s+1} \mbf{K}\left( \mbf{P}^*_{s,b} || \mbf{P}^{ind}_{s,b} \right)\\
&+ (s+1) L^2 K^L n^{-\frac{2}{d+1}} \log^{d+2+\frac{2}{d+1}} n.
\end{align*}
In particular under Assumption \ref{hyp:hmm_ergodic}, there exist positive constants $C(Q^*)$ and $c(Q^*)$ such that for $s\geq c(Q^*) \log n$ we have
\begin{equation}
C(Q^*) \mbb{E} \left[ h^2\left( P^*, \hat{P}_s \right) \right] \leq h^2\left( P^*, \overline{\mscr{M}} \right) + s L^2 K^L n^{-\frac{2}{d+1}} \log^{d+2+\frac{2}{d+1}} n,
\end{equation}
where $P^*$ is given by (\ref{eq:p_star}).
\end{theoreme}
Inequality (\ref{eq:th_hmm_log_concave_1}) does not require any assumption on the data. We can deduce convergence rates over the class $\mscr{H}^*\left(K,\mscr{F}_d,\ldots,\mscr{F}_d\right)$. For the optimal choice of $s$, we have
\[ C(Q^*) \mbb{E}\left[ h^2\left(P^*,\hat{P}_s \right) \right] \leq L^2 K^L n^{-\frac{2}{d+1}} \log^{d+3+\frac{2}{d+1}} n \]
for all $P^*\in\mscr{H}^*\left(K,\mscr{F}_d,\ldots,\mscr{F}_d\right)$. As for $d\leq 3$, we have the same rate as in Section \ref{sec:exponentially_ergodic_estimation} with a worse power of $\log n$ due to the higher complexity of HMMs. Our estimator is also robust to contamination and outliers. We can see from (\ref{eq:hmm_contamination_bound}) that our bound is not significantly worse as long as the contamination rate $\alpha_{cont}$ is of order not larger than $n^{-\frac{2}{d+1}} \log^{d+3+\frac{2}{d+1}} n$. Similarly (\ref{eq:hmm_outliers_bound}) tells us that a number of outliers of order not larger than $n^{\frac{d-1}{d+1}} \log^{d+3+\frac{2}{d+1}} n$ does not significantly deteriorate our bound on the convergence rate of our estimator.
\subsubsection{Exponential families as emission models}
\label{sec:exponential_family}
We introduce exponential families as follow. Let $d$ be a positive integer and $\eta:\overline{\Theta}\rightarrow \mbb{R}^d$ be a function over a non-empty set $\overline{\Theta}$. Let $T:\mscr{Y} \rightarrow \mbb{R}^d$ and $B:\mscr{Y}\rightarrow \mbb{R}$ be measurable functions such that
\[ \int_{\mscr{Y}} e^{\langle \eta(\theta), T(x)\rangle + B(x)} \nu(dx) <\infty, \forall  \theta\in\overline{\Theta},\]
we denote by $\mcal{E}\left(\overline{\Theta},\eta,T,d,B\right)$ the exponential family defined by
\begin{equation}
\label{eq:f_theta}
\mcal{E}\left(\overline{\Theta},\eta,T,d,B\right) := \left\{ f_{\theta} : x \mapsto e^{\langle \eta(\theta), T(x)\rangle + A(\theta) + B(x)}; \theta\in\overline{\Theta}\right\}, 
\end{equation}
where
\[ A(\theta) := - \log\left( \int_{\mscr{Y}} e^{\langle \eta(\theta), T(x)\rangle + B(x)} \nu(dx) \right).\]
It is a set of probability density functions with respect to $\nu$.
\begin{hypothese}
\label{hyp:exponential_family}
For all $k\in \{1,\dots,K\}$,
\begin{enumerate}
\item $\overline{\mscr{F}}_k$ is of the form
\begin{equation}
\label{eq:exponenatial_fam_i}
\overline{\mscr{F}}_k = \left\{ q \cdot \nu; q\in\mcal{E}\left(\overline{\Theta}_k,\eta_k,T_k,d_k,B_k\right) \right\},
\end{equation}
\item $\Theta_k$ is a countable subset of $\overline{\Theta}_k$ such that
\[ \mscr{F}_k=\left\{ q\cdot\nu; q\in\mcal{E}\left(\Theta_k,\eta_{k|\Theta_k},T_k,d_k,B_k\right) \right\} \]
is a dense subset of $\overline{\mscr{F}}_k$.
\end{enumerate}
\end{hypothese}
The next result is proven in Section \ref{sec:proof_vc_exponential} and shows that the last assumption is sufficient to satisfy our main assumption.
\begin{prop}
\label{prop:exponential_family}
Under Assumption \ref{hyp:exponential_family}, we satisfy Assumption \ref{hyp:hmm_vc} with $\epsilon=0$ and $V_{k_1,\dots,k_L} = 3 + \suml_{k_l=1}^L d_{k_l}$. Therefore we have
\begin{equation}
\label{eq:vc_exponential_family}
\overline{V} = 3 K^L + L K^{L-1} \left( d_1+\dots+d_K\right). 
\end{equation}
\end{prop}
We can see that the constant $\overline{V}$ does not depend on $\mscr{X}$ but on the dimensions $d_1,\dots,d_K$ which is the actual measure of the complexity of the exponential families. To our knowledge, the existence of a countable dense subset is satisfied for all the common exponential families. We obtain the following result for $\overline{\mscr{M}}\subset \mscr{H}\left(K,\overline{\mscr{F}}_1,\dots,\overline{\mscr{F}}_K\right)$.
\begin{coro}
\label{coro:exponential_family}
Let $N\geq K+L$ and $\hat{P}_s=\hat{P}_{s,\delta}$ be the estimator given by (\ref{eq:p_hat_s_delta}) with $\delta$ given by (\ref{eq:delta_s}). There exists a positive constant $C$ such that for all $\overline{P}\in\PX$, we have
\begin{align*}
C \mbb{E}\left[ h^2\left( \overline{P}, \hat{P}_{s,\delta} \right) \right] &\leq h^2\left( \overline{P}, \overline{\mscr{M}} \right) + n^{-1} \suml_{i=1}^n h^2\left( \overline{P},P_i \right)\\
&+ n^{-1} \suml_{b=1}^{s+1} \mbf{K}\left(\mbf{P}^*_{s,b}||\mbf{P}^{ind}_{s,b}\right)\nonumber\\
&+ (s+1) L K^{L-1} \left( K + L (d_1+\dots+d_K) \right) \log n.
\end{align*}
In particular under Assumption \ref{hyp:hmm_ergodic}, there exist positive constants $C(Q^*)$ and $c(Q^*)$ such that for $s \geq c(Q^*) \log n$ we have
\begin{align}
\label{eq:bound_exponential_family}
C(Q^*) \mbb{E}\left[ h^2\left(P^*,\hat{P}_s \right) \right] &\leq h^2\left( P^*, \overline{\mscr{M}} \right)\\
&+ L K^{L-1} \left( K + L (d_1+\dots+d_K) \right) \frac{s \log n}{n},\nonumber
\end{align}
where $P^*$ is given by (\ref{eq:p_star}).
\end{coro}
This result is a direct consequence of Theorem \ref{th:first_hmm} and Proposition \ref{prop:exponential_family}. We can deduce a bound on the convergence rate over $\mscr{H}^*\left(K,\overline{\mscr{F}}_1,\dots,\overline{\mscr{F}}_K\right)$. For the optimal choice of $s$, we have
\[ C(Q^*) \mbb{E}\left[ h^2\left( P^*, \hat{P}_s \right) \right] \leq L K^{L-1} \left( K + L (d_1+\dots+d_K) \right) \frac{\log^2 n}{n}, \]
for all $P^*$ in $\mscr{H}^*\left(K,\overline{\mscr{F}}_1,\dots,\overline{\mscr{F}}_K\right)$. We obtain the optimal $1/n$ rate with respect to the squared Hellinger distance, up to a logarithmic factor. Corollary \ref{coro:hmm_robust} shows that our estimator is also robust to contamination and outliers. From (\ref{eq:hmm_contamination_bound}) we see that our bound is not significantly worse as long as the contamination rate $\alpha_{cont}$ is of order not larger than $n^{-1} \log^2 n$. Similarly, we get from (\ref{eq:hmm_outliers_bound}) that the performance of our estimator is not altered as long as the number of outliers $|I|$ is of order not larger than $\log^2 n$.\par
Let us illustrate how Corollary \ref{coro:exponential_family} applies with the following example. Let $d$ be a positive integer and $\text{Cov}_{+*}(d)$ be the set of  $d\times d$ symmetric and positive-definite matrices.  For $z\in\mbb{R}^d$ and $\Sigma\in\text{Cov}_{+*}(d)$, we denote by $g_{z,\Sigma}$ the density function of  the normal distribution $\mcal{N}(z,\Sigma)$ with respect to the Lebesgue measure given by
\begin{equation}
\label{eq:g_mu_sigma}
g_{z,\Sigma}(x) := \frac{1}{\sqrt{(2\pi)^d |\Sigma|}} \exp\left( -\frac{(z-m)^T\Sigma^{-1}(z-m)}{2} \right), 
\end{equation}
where $|\Sigma|$ denotes the determinant of $\Sigma$. Let $\mcal{G}_d$ be the location-scale family of densities given by $\mcal{G}_d := \{ g_{z,\Sigma}; z\in\mbb{R}^d, \Sigma\in\text{Cov}_{+*}(d)\}$. One can check it is an exponential family with $\mcal{G}_d=\mcal{E}\left( \mbb{R}^d \times \text{Cov}_{+*}(d), \eta, T, \frac{d(d+3)}{2}, 0 \right)$ where
\begin{align*}
T(x) &= \left( x, \left( x_i^2 \right)_{1\leq i\leq d}, \left( x_i x_j \right)_{1\leq i<j\leq d} \right) \text{  and  }\\
\eta(z,\Sigma) &= \left( \Sigma^{-1}z, -\frac{1}{2}\left(\Sigma^{-1}_{ii}\right)_{1\leq i\leq d}, -\left( \Sigma^{-1}_{ij}\right)_{1\leq i<j\leq d} \right).
\end{align*}
For a fixed $\Sigma$ we denote by $\mcal{G}_{loc}(\Sigma)$ the associated location family given by $\mcal{G}_{loc}(\Sigma) := \{ g_{z,\Sigma}; z\in\mbb{R}^d\}$. It is also an exponential family with $\mcal{G}_{loc}(\Sigma)=\mcal{E}\left( \mbb{R}^d\times\text{Cov}_{+*}(d), \eta, T, d, B \right)$, where
\[ \eta(z) = \Sigma^{-1}z, T(x) = x\text{  and  } B(x) = -\frac{x^T \Sigma^{-1} x}{2}.\]
We denote by $\mscr{G}_d$ and $\mscr{G}_{loc}(\Sigma)$ respectively, the sets of probability distributions associated to $\mcal{G}_d$ and $\mcal{G}_{loc}(\Sigma)$. The next result is a consequence of Corollary \ref{coro:exponential_family}.
\begin{theoreme}
\label{th:hmm_multivariate_gaussian}
Let $N\geq K+L$ and $Y_1,\dots,Y_N$ be arbitrary random variables.
\begin{itemize}
\item  Let $\hat{P}_s=\hat{P}_{s,\delta}$ be the estimator given by (\ref{eq:p_hat_s_delta}) with $\overline{\mscr{M}}=\mscr{H}\left(K,\mscr{G}_d,\dots,\mscr{G}_d\right)$ and $\delta$ given by (\ref{eq:delta_s}). There exists a positive constant $C$ such that for all $\overline{P}\in\PX$
\begin{align}
C \mbb{E}\left[ h^2\left( \overline{P}, \hat{P}_s \right) \right] &\leq h^2\left( \overline{P}, \overline{\mscr{M}} \right) + n^{-1} \suml_{i=1}^n h^2\left( \overline{P},P_i \right)\nonumber\\
&n^{-1} \suml_{b=1}^{s+1} \mbf{K}\left( \mbf{P}^*_{s,b} || \mbf{P}^{ind}_{s,b} \right)\nonumber\\
&+ (s+1) L^2 K^L d(d+3) \frac{\log n}{n}.\label{eq:bound_hmm_gauss}
\end{align}
In particular under Assumption \ref{hyp:hmm_ergodic} there exist positive constants $C(Q^*)$ and $c(Q^*)$ such that for $s\geq c(Q^*) \log n$ we have
\[ C(Q^*) \mbb{E}\left[ h^2\left( P^*, \hat{P}_s \right) \right] \leq h^2\left( P^*,\mscr{M} \right) + (s+1) L^2 K^L d(d+3) \frac{\log n}{n}, \]
where $P^*$ is given by (\ref{eq:p_star}). 
\item Let $\hat{P}_s=\hat{P}_{s,\delta}$ be the estimator given by (\ref{eq:p_hat_s_delta}) with $\overline{\mscr{M}}=\mscr{H}\left(K,\mscr{G}_{loc}(\Sigma),\dots,\mscr{G}_{loc}(\Sigma)\right)$ and $\delta$ given by (\ref{eq:delta_s}). There exists a positive constant $C$ such that for all $\overline{P}\in\PX$
\begin{align}
C \mbb{E}\left[ h^2\left( \overline{P}, \hat{P}_s \right) \right] &\leq h^2\left( \overline{P}, \overline{\mscr{M}} \right) + n^{-1} \suml_{i=1}^n h^2\left( \overline{P},P_i \right)\nonumber\\
&+ n^{-1} \suml_{b=1}^{s+1} \mbf{K}\left( \mbf{P}^*_{s,b} || \mbf{P}^{ind}_{s,b} \right)\nonumber\\
&+ (s+1) L^2 K^L d \frac{\log n}{n},\label{eq:bound_hmm_gauss_fixed_covariance}
\end{align}
for any $\Sigma$ in $\text{Cov}_{+*}(d)$. In particular under Assumption \ref{hyp:hmm_ergodic} there exist positive constants $C(Q^*)$ and $c(Q^*)$ such that for $s\geq c(Q^*) \log n$ we have
\[ C(Q^*) \mbb{E}\left[ h^2\left( P^*, \hat{P}_s \right) \right] \leq h^2\left( P^*, \overline{\mscr{M}} \right) + (s+1) L^2 K^L d \frac{\log n}{n}, \]
where $P^*$ is given by (\ref{eq:p_star}). 
\end{itemize}
\end{theoreme}
Inequalities (\ref{eq:bound_hmm_gauss}) and (\ref{eq:bound_hmm_gauss_fixed_covariance}) are consequences of Corollary \ref{coro:exponential_family} and do not require any assumption on the data. We deduce bounds on the convergence rate of our estimator over $\mscr{H}^*\left(K,\mscr{G}_d,\dots,\mscr{G}_d\right)$ and $\mscr{H}^*\left(K,\mscr{G}_{loc}(\Sigma),\dots,\mscr{G}_{loc}(\Sigma)\right)$. For the optimal choice of $s$ we obtain the rate $n^{-1} \log^2 n$ with respect to the squared Hellinger distance both for $P^*\in\mscr{H}^*\left(K,\mscr{G}_d,\dots,\mscr{G}_d\right)$ and $P^*\in\mscr{H}^*\left(K,\mscr{G}_{loc}(\Sigma),\dots,\mscr{G}_{loc}(\Sigma)\right)$. This rate is optimal up to a logarithmic factor. We can see that the dependence on the dimension $d$ is linear for the model $\mscr{H}\left(K,\mscr{G}_{loc}(\Sigma),\dots,\mscr{G}_{loc}(\Sigma)\right)$ while its quadratic for $\mscr{H}^*\left(K,\mscr{G}_d,\dots,\mscr{G}_d\right)$.\par
We can obtain similar results for any exponential family. It is also possible to consider hidden Markov models with different exponential families as emission models. The next section investigates the estimation of the parameters.
\subsubsection*{Estimation of the parameters with emission exponential families}
We say that $\hat{\pi}$, $\hat{Q}$ and $\hat{F}$ are $\rho$-estimators of $\pi^*$, $Q^*$ and $F^*$ if $P_{\hat{w},\hat{Q},\hat{F}}=\hat{P}_{s,\delta}$ is an estimator of $P^*$ given by (\ref{eq:p_hat_s_delta}). If we consider models of densities that are uniformly bounded, we can use (\ref{eq:l_2_hellinger}) and Theorem 9 of Lehéricy \cite{lehericyorder} to deduce risk bounds for the parameter estimators. It is also possible to use the results of Ibragimov and Has’minski\u{\i} \cite{Ibragimov} for regular parametric models.\par
We consider that Assumption \ref{hyp:exponential_family} is satisfied with $\overline{\Theta}_k\subset\mbb{R}^{e_k}$ for all $k\in[K]$. For $k\in[K]$ we denote by $F_{\theta_k}$ the probability distribution given by the parameter $\theta_k\in\overline{\Theta}_k$, i.e. $F_{\theta_k}=f_{\theta_k}\cdot\nu$ with $f_{\theta}$ given by (\ref{eq:f_theta}). Let $\overline{\Phi}$ be an open convex subset of $O_K^{K+1} \times \overline{\Theta}_1 \times \dots \times \overline{\Theta}_K$, where
\[O_K=\left\{ \mbf{a}\in(0,1)^{K-1}, a_1+\dots+a_{K-1}<1 \right\}.\]
For $\phi$ in $\overline{\Phi}$, we can define $w\in\mcal{W}_K$, $Q\in\mcal{T}_K$ and $\mbf{\theta}\in \overline{\Theta}_1 \times \dots \times \overline{\Theta}_K$ by $\phi = (\phi_w,\phi_{Q,1},\dots,\phi_{Q,K},\phi_{\theta})$ with
\begin{align*}
(w_1,\ldots,w_{K-1}) &= \phi_w \in O_K,\\
(Q_{k,1},\ldots,Q_{K-1,1}) &= \phi_{Q,k} \in O_K,\\
(\theta_1,\ldots,\theta_K) &= \phi_{\theta} \in \overline{\Theta}_1\times\dots\times\overline{\Theta}_K.
\end{align*}
We denote by $\overline{\mscr{M}}$ the model given by
\begin{equation}
\label{eq:model_parametric_regular}
\overline{\mscr{M}} := \left\{ P_{\phi} = p(\cdot;\phi)\cdot \mu ; \phi \in \overline{\Phi}  \right\}
\end{equation}
and
\[ p(\mbf{x};\phi)= \suml_{1\leq k_1,\dots,k_L\leq K} w_{k_1} Q(k_2|k_1) \dots Q(k_L|k_{L-1}) \prod\limits_{l=1}^L f_{\theta_{k_l}}(x_l). \]
We need the following assumption to make sure we can deduce $\phi$ from $P_{\phi}$.
\begin{hypothese}
\label{hyp:regular_parametric_model}
For all $k$ in $[K]$,
\begin{itemize}
\item the map $\theta_k\mapsto F_{\theta_k}$ is continuous on $\overline{\Theta}_k$ with respect to the Hellinger distance;
\item the functions $\eta_k$ and $A_k$ are of class $\mcal{C}^1$ on $\overline{\Theta}_k$;
\item for all $\theta_k$ in $\overline{\Theta}_k$, we have $\int ||T_k(x)||^2 f_{\theta_k}(x) \nu(dx) < \infty$ and
\[ \int ||T_k(x)||^2 \left| f_{\theta_k}(x) - f_{\theta'_k}(x) \right| \nu(dx) \xrightarrow[||\theta_k-\theta'_k||\rightarrow 0]{} 0. \]
\end{itemize}
\end{hypothese}
The next result is proven in Section \ref{sec:proof_prop_regular} and shows that under some conditions we can deduce the parameters from the distribution $P_{\phi}$. 
\begin{prop}
\label{prop:regular_parametric_model}
Under Assumption \ref{hyp:regular_parametric_model} the information matrix $I$ function given by 
\[ I_{ij}:\phi\mapsto I\left(\phi\right)_{ij} = \int_{\mscr{X}^L} \partial_{\phi_i} p(\mbf{x};\phi) \partial_{\phi_j} p(\mbf{x};\phi) \frac{\mu(d\mbf{x})}{p(\mbf{x};\phi)} \]
is well-defined and continuous on $\overline{\Phi}$. We define the subset $\Phi^*\subset\overline{\Phi}$ by
\begin{equation}
\label{eq:phi_star}
\Phi^* := \left\{ \overline{\phi}\in\overline{\Phi}; \begin{array}{c}
I\left(\overline{\phi}\right) \text{ is definite positive and}\\
\inf\limits_{\substack{||\overline{\phi}-\phi||\geq a\\ \phi\in \overline{\Phi}}} h^2\left(P_{\overline{\phi}},P_{\phi}\right) >0, \forall a>0,
\end{array}
\right\}
\end{equation}
For all $\phi^* \in \Phi^*$, there exists a positive constant $C( \phi^*)$ such that
\begin{equation}
\label{eq:prop_ibragimov}
C(\phi^*) \left[ \left|\left| w^* - w \right|\right|_2^2 + \left|\left| Q^* - Q \right|\right|_2^2 + \suml_{k=1}^K \left|\left| \theta^*_k - \theta_k \right|\right|_2^2\wedge 1 \right]  \leq  h^2\left( P_{\phi^*}, P_{\phi} \right),
\end{equation}
for all $\phi$ in $\overline{\Phi}$.
\end{prop}
The constant $C(\phi^*)$ depends on the inverse of the smallest eigenvalue of $I(\phi^*)$ and the geometry of $\overline{\Phi}$ around $\phi^*$ induced by the Hellinger distance on $\overline{\mscr{M}}$. The next result is a consequence of Proposition \ref{prop:regular_parametric_model} and Corollary \ref{coro:exponential_family}.
\begin{theoreme}
\label{th:ibragimov}
Let $N\geq K+L$ and $Y_1,\dots,Y_N$ be arbitrary random variables. Let $P_{\hat{\phi}}=\hat{P}_{s,\delta}$ be the estimator given by (\ref{eq:p_hat_s_delta}) with $\delta$ given by (\ref{eq:delta_s}). Under Assumption \ref{hyp:regular_parametric_model}, for all $\overline{\phi}\in\Phi^*$ there exists a positive constant $C(\overline{\phi})$ such that
\begin{align}
C\left(\overline{\phi}\right) &\mbb{E}\left[ \left|\left|\overline{w}-\hat{w}\right|\right|_2^2 + \left|\left| \overline{Q} - \hat{Q} \right|\right|_2^2 + \suml_{k=1}^K \left|\left|\overline{\theta}_k - \hat{\theta}_k \right|\right|_2^2\wedge 1 \right]\nonumber\\
&\leq n^{-1} \suml_{i=1}^n h^2\left( P_{\overline{\phi}}, P_i \right) + n^{-1} \suml_{b=1}^{s+1} \mbf{K}\left(\mbf{P}^*_{s,b}||\mbf{P}^{ind}_{s,b}\right)\nonumber\\
&+ (s+1) L K^{L-1} \left( K + L (d_1+\dots+d_K) \right) \frac{\log n}{n}.\label{eq:th_ibragimov}
\end{align}
In particular under Assumption \ref{hyp:hmm_ergodic}, there exist positive constants $C(\overline{\phi},Q^*)$ and $c(Q^*)$ such that for $s \geq c(Q^*) \log n$ we have
\begin{align}
C\left(\overline{\phi},Q^*\right) &\mbb{E}\left[ \left|\left|\overline{w}-\hat{w}\right|\right|_2^2 + \left|\left| \overline{Q} - \hat{Q} \right|\right|_2^2 + \suml_{k=1}^K \left|\left|\overline{\theta}_k - \hat{\theta}_k \right|\right|_2^2\wedge 1 \right]\label{eq:th_ibragimov_2}\\
&\leq h^2\left( P^*, P_{\overline{\phi}} \right) + L K^{L-1} \left( K + L (d_1+\dots+d_K) \right) \frac{s \log n}{n},\nonumber
\end{align}
where $P^*$ is given by (\ref{eq:p_star}).
\end{theoreme}
Inequality (\ref{eq:th_ibragimov}) is a consequence of Proposition \ref{prop:regular_parametric_model} and Corollary \ref{coro:exponential_family}. It does not require any assumption on the data and shows that the estimators of the parameters can be meaningful even if the model is misspecified. Ideally there exists $\overline{\phi}$ in $\Phi^*$ such that most of the distributions $P_i$ lie in a small neighborhood of $P_{\overline{\phi}}$ so that the first term of our bound is small compared to the last term. In that case the estimators $\hat{w}, \hat{Q}, \hat{\theta}_1,\dots,\hat{\theta}_K$ converge to a small neighborhood around $\overline{w},\overline{Q},\overline{\theta}_1,\ldots,\overline{\theta}_K$, where $P_{\overline{\phi}}$ should be seen as the best approximation of the true distribution in the model. We can deduce bounds on the convergence rate of our parameter estimators in the well-specified case from (\ref{eq:th_ibragimov_2}). For $P^*=P_{\phi^*}\in\mscr{H}^*\left(K,\overline{\mscr{F}}_1,\dots,\overline{\mscr{F}}_K\right)$ with $\phi^*\in\Phi^*$ and for the optimal choice of $s$, we retrieve the usual parametric rate for each parameter estimator, up to a logarithmic factor. Let us illustrate this with the following example.\par
We consider exponential distributions for the emission models, i.e. we have $\overline{\mscr{F}}_i=\overline{\mscr{E}}$ for all $i$ in $[K]$ with
\begin{equation}
\label{eq:exponential_distribution}
\overline{\mscr{E}} := \left\{ f_{\theta} \cdot \nu; f_{\theta}\in\mcal{E}\left(\overline{\Theta},\text{id}_{\Theta}, -\text{id}_{\mscr{X}},1,0 \right) \right\}
\end{equation}
where $\overline{\Theta}=(0,\infty)$, $\mscr{X}=[0,\infty)$, $\nu$ is the Lebesgue measure on $\mscr{X}$, and we can deduce $A:\theta\mapsto\log \theta$. This means we have $f_{\theta}:x\mapsto \theta e^{-\theta x}\mathbbm{1}_{x\geq 0}$ for any $\theta>0$. One can easily check that we satisfy Assumption \ref{hyp:regular_parametric_model}, the last condition being a direct consequence of the dominated convergence theorem. We define $\overline{\Phi}$ by
\begin{equation}
\label{eq:phi_exponential}
\overline{\Phi} = O_K^{K+1} \times \left\{ \mbf{\theta}\in \Theta^K; \theta_1>\theta_2>\dots>\theta_K \right\},
\end{equation}
and $\overline{\mscr{M}}$ as in (\ref{eq:model_parametric_regular}). The condition on the parameters $\theta$ ensures identifiability over $\overline{\Phi}$ and $\overline{\Phi}^*=\overline{\Phi}$. The choice $L=3$ is enough to obtain the result of Proposition \ref{prop:regular_parametric_model}. The next theorem is proven in Section \ref{sec:proof_th_ibragimov_exponential}.
\begin{theoreme}
\label{th:ibragimov_exponential}
Let $N\geq K+3$ and $Y_1,\dots,Y_N$ be arbitrary random variables. Let $P_{\hat{\phi}}=\hat{P}_{s,\delta}$ be the estimator given by (\ref{eq:p_hat_s_delta}) with $\delta$ given by (\ref{eq:delta_s}). For any $\overline{\phi}$ in $\overline{\Phi}$ there exists a positive constant $C(\overline{\phi})$ such that we have 
\begin{align*}
&C\left(\overline{\phi}\right) \mbb{E}\left[ \left|\left|\overline{w}-\hat{w}\right|\right|_2^2 + \left|\left| \overline{Q} - \hat{Q} \right|\right|_2^2 + \suml_{k=1}^K \left( \overline{\theta}_k - \hat{\theta}_k \right)^2\wedge 1 \right]\\
&\leq n^{-1} \suml_{i=1}^n h^2\left( P_{\overline{\phi}}, P_i \right) + n^{-1} \suml_{b=1}^{s+1} \mbf{K}\left(\mbf{P}^*_{s,b}||\mbf{P}^{ind}_{s,b}\right) + (s+1) K^3 \frac{\log n}{n}.\nonumber
\end{align*}
In particular under Assumption \ref{hyp:hmm_ergodic}, there exist positive constants $C(\overline{\phi},Q^*)$ and $c(Q^*)$ such that for $s \geq c(Q^*) \log n$ we have
\begin{align}
\label{eq:th_ibragimov_exponential}
C\left(\overline{\phi},Q^*\right) &\mbb{E}\left[ \left|\left|\overline{w}-\hat{w}\right|\right|^2 + \left|\left|\overline{Q}-\hat{Q}\right|\right|^2 + \suml_{k=1}^K \left( \overline{\theta}_k - \hat{\theta}_k \right)^2\wedge 1 \right]\\
&\leq  h^2\left(P^*,P_{\overline{\phi}}\right) + s K^3 \frac{\log n}{n},\nonumber
\end{align}
where $P^*$ is given by (\ref{eq:p_star}).
\end{theoreme}
Our different parameter estimators all reach the usual parametric rate up to a logarithmic factor. One can notice that the ordering of the $\theta_k$ in (\ref{eq:phi_exponential}) can be replaced by considering only distinct values and taking the infimum over permutation of the hidden states.\par
It is possible to follow the same scheme to obtain similar results for other exponential families, including HMMs with different exponential families as emission models. The difficulty relies in determining the set $\Phi^*$ given by (\ref{eq:phi_star}).   
\subsubsection{Another example}
\label{sec:specific_example}
In this section we consider a relatively simple example that does not fit any framework already investigated but for which we can obtain risk bounds for the estimation of the parameters. Let $\nu$ be the Lebesgue measure on $\mbb{R}$ and $\alpha$ be in $(0,1)$. We denote by $f_{\alpha}$ the probability density function with respect to $\nu$ defined by
\[ f_{\alpha} : x\in\mbb{R}\mapsto \frac{1-\alpha}{2} \frac{\mathbbm{1}_{|x|\in[0,1]}}{|x|^{\alpha}},\]
with the convention $1/0=+\infty$.  For $z$ in $\mbb{R}$, we denote by $F_{\alpha,z}$ the probability distribution associated to the density $x\mapsto f_{\alpha}(x-z)$. We fix $L=2$ and consider the model $\overline{\mscr{M}}$ defined by
\[ \overline{\mscr{M}} = \left\{ P_{w,q,z}; w,q_{12},q_{21}\in[0,1], z\in\mbb{R} \right\},\]
where
\begin{align*}
P_{w,q,z} &= w F_{\alpha,0} \otimes \left[ (1-q_{12}) F_{\alpha,0} + q_{12} F_{\alpha,z} \right]\\
&+ (1-w) F_{\alpha,z} \otimes \left[ q_{21}  F_{\alpha,0} + (1-q_{21}) F_{\alpha,z} \right]. 
\end{align*}
The distributions $P_{w,q,z}$ correspond to translation hidden Markov models with one known location parameter. The following result is proven in Section \ref{sec:proof_faster_rates_prop} and shows that we can deduce the parameters from the distribution $P_{w,q,z}$.
\begin{prop}
\label{prop:faster_rates}
For $z^*\neq 0$, $w^*<1$ and $q^*_{21}<1$, there is a constant $C(\alpha,z^*,w^*,q^*)$ such that we have
\begin{align*}
C(\alpha,z^*,w^*,q^*) h^2\left( P_{w,q,z}, P_{w^*,q^*,z^*} \right) &\geq \left(|z-z^*|\wedge 1\right)^{1-\alpha} + (w^*)^2 \left( q_{12} - q_{12}^* \right)^2\\
&+ (1-w^*)^2 \left( q_{12} - q_{12}^* \right)^2 + \left(w-w^*\right)^2,
\end{align*}
for all $w,q_{12},q_{21}\in[0,1]$ and all $z\in\mbb{R}$.
\end{prop}
We can deduce a deviation bound for the parameter estimators. The model $\overline{\mscr{M}}$ is a subset of $\mscr{H}(2,\overline{\mscr{F}}_{\alpha}, \overline{\mscr{F}}_{\alpha})$, with $\overline{\mscr{F}}_{\alpha}=\{F_{\alpha,z};z\in\mbb{R}\}$. We satisfy Assumption \ref{hyp:hmm_vc} with $\epsilon=0$, $\mcal{F}_{\alpha}=\{f_{\alpha}(\cdot-z);z\in\mbb{Q}\}$ and $\overline{V}=784$. The next result is proven in Section \ref{sec:proof_th_faster_rates}.
\begin{theoreme}
\label{th:faster_rates}
Let $N\geq K+2$ and $P_{\hat{w},\hat{q},\hat{z}}=\hat{P}_{s,\delta}$ be the estimator given by (\ref{eq:p_hat_s_delta}) with $\delta$ given by (\ref{eq:delta_s}). For all $\overline{z}\neq0$, $\overline{w}<1$, $\overline{q}_{12}\in[0,1]$ and $\overline{q}_{21}<1$, there exists a positive constant $C(\alpha,\overline{z},\overline{w},\overline{q})$ such that we have 
\begin{align}
&C\left(\alpha,\overline{z},\overline{w},\overline{q}\right) \mbb{E}\left[ (\overline{w}-\hat{w})^2 + (\overline{q}_{12} - \hat{q}_{12})^2 + (\overline{q}_{12} - \hat{q}_{12})^2 + (|\overline{z} - \hat{z}|\wedge 1)^2 \right]\nonumber\\
&\leq \frac{1}{n} \suml_{i=1}^n h^2\left( P_{\overline{w},\overline{q},\overline{z}}, P_i \right) + \frac{1}{n} \suml_{b=1}^{s+1} \mbf{K}\left(\mbf{P}^*_{s,b}||\mbf{P}^{ind}_{s,b}\right)+ (s+1) \frac{\log n}{n}.\label{eq:th_faster_rates}
\end{align}
In particular under Assumption \ref{hyp:hmm_ergodic}, there exist positive constants $C(\overline{\phi},Q^*)$ and $c(Q^*)$ such that for $s \geq c(Q^*) \log n$ we have
\begin{align}
C\left( \overline{\phi}, Q^* \right) &\mbb{E}\left[ (\overline{w}-\hat{w})^2 + (\overline{q}_{12} - \hat{q}_{12})^2 + (\overline{q}_{21} - \hat{q}_{21})^2 + \left( |\overline{z} - \hat{z}|\wedge 1 \right)^{1-\alpha} \right]\nonumber\\
&\leq  h^2\left( P^*, P_{\overline{w},\overline{q},\overline{z}}\right) + \frac{s \log n}{n},\label{eq:th_faster_rates_2}
\end{align}
where $P^*$ is given by (\ref{eq:p_star}).
\end{theoreme}
Inequality (\ref{eq:th_faster_rates}) does not require any assumption on the data. It is a consequence of Proposition \ref{prop:faster_rates} and Theorem \ref{th:first_hmm}. We can deduce convergence rates for our parameter estimators from (\ref{eq:th_faster_rates_2}) for $P^*=P_{\pi^*,q^*,z^*}$ with $z^*\neq0$, $w^*<1$ and $q^*_{21}<1$. The estimators $\hat{w}$ and $\hat{q}$ achieve the usual parametric rate up to a logarithmic factor. However the location estimator $\hat{z}$ reaches the faster rate $(n^{-1}\log^2 n)^{1/(1-\alpha)}$. This rate is optimal up the logarithmic factor. It is a consequence of Theorem 1.1 in \cite{Ibragimov} (Chapter VI), noticing that $f_{\alpha}$ has a singularity of order $-\alpha$ in 0, and with the fact that we cannot do better than $1/n$ for the Hellinger distance. One should notice that $f_{\alpha}$ is unbounded for all $\alpha\in(0,1)$. Therefore the maximum likelihood and the least squares estimators are undefined and those methods do not apply on $\overline{\mscr{M}}$. In addition, we can see that $f_{\alpha}$ is not square integrable for $\alpha\in[1/2,1)$.
\section{Selection of the spacing parameter}
\label{sec:selection_s}
Until now we gave results that required a good choice of the spacing parameter $s$, given some bound on the dependence term $\mbf{K}\left( \mbf{P}^*_{s,b} || \mbf{P}^{ind}_{s,b} \right)$. This section propose a way to automatically select a value of $s$ from the data, assuming that we dispose of two independent sets of observations. We use the first set to produce an estimator $\hat{P}_s$ for different values of $s$. We then use the second set to produce an estimator $\hat{s}$ of the optimal value of $s$.
\subsection{Framework and result}
Let $X^{(1)}_1,\ldots,X^{(1)}_{n_1},X^{(2)}_1,\dots,X^{(2)}_{n_2}$ be $n_1+n_2$ random variables on the measurable space $\left(\mscr{X},\mcal{X}\right)$. We define $P^{(j)}_i$ by $P^{(j)}_i := \mcal{L}(X^{(j)}_i)$ for all $j$ in $[2]$ and all $i$ in $[n_j]$. We also write
\[ \mbf{P}^*_{s,b} = \mcal{L}\left(X^{(1)}_b,\ldots,X_{b+n_1(s,b)(s+1)}\right) \text{  and  } \mbf{P}^{ind}_{s,b} = \bigotimes_{i=1}^{n_1(s,b)} \mcal{L}\left( X^{(1)}_{b+(i-1)(s+1)} \right), \]
with
\begin{equation}
\label{eq:n_1_s_b}
n_1(s,b) = \left\lfloor \frac{n_1+s+1-b}{1+s} \right\rfloor. 
\end{equation}
Let $S$ be a subset of $\{0,1,\ldots,s_{\max}\}$, $s_{\max}=\lfloor (n_1-2)/2\rfloor$. Let $\left(\mscr{M}_s\right)_{s\in S}$ be countable subsets of $\PX$ such that the $\rho$-dimension function (see Section \ref{sec:rho_application}) is uniformly bounded over $\mscr{M}_s$ by a non-decreasing function $m\mapsto D_m(\mscr{M}_s)\geq 1$ for all $s\in S$. We follow the procedure below.
\begin{enumerate}
\item For $s$ in $S$, let $\hat{P}_s=\hat{P}_s\left(\mscr{M}_s,\mbf{X}^{(1)}\right)$ be the estimator given by (\ref{eq:p_hat_s}). Conditionally on $\mbf{X}^{(1)}$, we define the finite model
\[ \widehat{\mscr{M}}_S = \widehat{\mscr{M}}_S\left(\mbf{X}^{(1)}\right) := \left\{ \hat{P}_s : s \in S\right\}. \]
\item Let $\hat{P}$ be the $\rho$-estimator $\hat{P} = \hat{P}\left(n_2,\mbf{X}^{(2)},\widehat{\mscr{M}}_S\right)$ given by (\ref{eq:rho_estimateur}). We denote by $\hat{s}$ the value of $s$ such that $\hat{P}=\hat{P}_{\hat{s}}$ and we write
\begin{equation}
\label{eq:p_hat_s_hat}
\hat{P}=\hat{P}_{\hat{s}}\left(\mbf{X}^{(1)},\mbf{X}^{(2)}\right). 
\end{equation}
\end{enumerate}
We make the following assumption.
\begin{hypothese}
\label{hyp:selection_s_independent}
The random variables
\[ \mbf{X}^{(1)} := \left(X^{(1)}_1,\ldots,X^{(1)}_{n_1}\right) \text{  and  } \mbf{X}^{(2)} := \left(X^{(2)}_1,\ldots,X^{(2)}_{n_2}\right) \]
are independent.
\end{hypothese}
The following result is proven in Section \ref{sec:proof_th_rho_selection_s}.
\begin{theoreme}
\label{th:rho_selection_s}
Let $n_1,n_2\geq 3$ and $\hat{P}=\hat{P}_{\hat{s}}\left(\mbf{X}^{(1)},\mbf{X}^{(2)}\right)$ be the estimator given by (\ref{eq:p_hat_s_hat}). Under Assumption \ref{hyp:selection_s_independent}, there exists a positive constant $C>0$ such that for all $\overline{P}\in\PX$
\begin{align}
C \mbb{E}\left[ h^2\left(\overline{P}, \hat{P}_{\hat{s}}\right) \right] &\leq n_1^{-1} \suml_{i=1}^{n_1} h^2\left( P_i^{(1)},\overline{P} \right) + n_2^{-1} \suml_{i=1}^{n_2} h^2\left( P_i^{(2)},\overline{P} \right)\label{eq:th_selection_s_main}\\
+ \inf_{t\in[n_2]} &\bigg\{ \frac{t}{n_2} \left( 1 + \log(|S|) \right) + \lceil n_2/t\rceil \beta_t\left( \mbf{X}^{(2)} \right) \bigg\}\nonumber\\
+ \inf_{s\in S} \bigg\{ & h^2\left( \overline{P} , \mscr{M}_s \right) + \frac{(s+1) D_{n_1(s,1)}(\mscr{M}_s)}{n_1} + n_1^{-1} \suml_{b=1}^{s+1} \mbf{K}\left(\mbf{P}^*_{s,b}||\mbf{P}^{ind}_{s,b}\right) \bigg\},\nonumber
\end{align}
where the mixing coefficient $\beta_{t}\left(\mbf{X}^{(2)}\right)$ is given by (1.2.5) in Dedecker \etal \cite{weakdependence}.
\end{theoreme}
One can check that we do not need any assumption other than Assumption \ref{hyp:selection_s_independent} to obtain this result. We need to make additional assumptions a posteriori to make this bound meaningful. Let us interpret this inequality in simpler cases. We consider there is $\mscr{M}$ such that $\mscr{M}_s=\mscr{M}$ for all $s\in S$. If the data were truly i.i.d.\ with distribution $\overline{P}\in\mscr{M}$, we would get
\[ C \mbb{E}\left[ h^2\left(\overline{P}, \hat{P}\right) \right] \leq \frac{(s+1) D_{n_1(s,1)}(\mscr{M})}{n_1}
+ \frac{\left( 1 + \log(|S|) \right)}{n_2}. \]
The second term is the bound we get for i.i.d. estimation from a $n_2$-sample over a finite model of cardinal $|S|$. When the data are not identically distributed, the quantity 
\[ n_2^{-1} \suml_{i=1}^{n_2} h^2\left( P_i^{(2)},\overline{P} \right) + n_1^{-1} \suml_{i=1}^{n_1} h^2\left( P_i^{(1)},\overline{P} \right) \]
is not zero but it remains small when most of the true marginal distributions $P^{(j)}_i$ lie close enough to some distribution $\overline{P}$ in $\mscr{M}$. The terms $n_1^{-1} \sum_{b=1}^{s+1} \mbf{K}\left(\mbf{P}^*_{s,b}||\mbf{P}^{ind}_{s,b}\right)$ and $\lceil n_2/t\rceil \beta_t( \mbf{X}^{(2)})$ account for the possible dependence within $\mbf{X}^{(1)}$ and $\mbf{X}^{(2)}$ respectively. They vanish if the observations $X^{(1)}_1,\dots,X^{(1)}_{n_1},X^{(2)}_1,\dots,X^{(2)}_{n_2}$ are all independent. Contrary to Theorem \ref{th:first_hmm} we do not have to choose a good value of $s$ as the method automatically select a reasonable $s$ in $S$ as long as the $P^{(j)}_i$ can be well approximated by a distribution $\overline{P}\in\mscr{M}$.
\subsection{Robustness}
\label{sec:robust_selection_s}
Let $\mbf{\overline{X}}^{(1)}=\left(\overline{X}^{(1)}_1,\ldots,\overline{X}^{(1)}_{n_1}\right)$ and $\mbf{\overline{X}}^{(2)}=\left(\overline{X}^{(2)}_1,\ldots,\overline{X}^{(2)}_{n_2}\right)$ be the true processes of interest such that $P^{(j)}_i=\overline{P}$ for all $j\in[2]$ and $i\in[N_j]$. We actually observe a contaminated version of it. Let $Z^{(1)}_1,\ldots,Z^{(1)}_{N_1},Z^{(2)}_1,\ldots,Z^{(2)}_{N_2}$ be random variables with any distributions and $E^{(1)}_1,\ldots,E^{(1)}_{N_1},E^{(2)}_1,\ldots,E^{(2)}_{N_2}$ be Bernoulli random variables such that for all $j\in[2]$ and all $i\in[N_j]$,
\begin{equation}
\label{eq:x_contamination_selection_s}
X^{(j)}_i = E_i \overline{X}^{(j)}_i + (1-E^{(j)}_i) Z^{(j)}_i.
\end{equation}
For $s\in\{0,1,\ldots,s_{\max}\}$ and $b\in[s+1]$, we define the distributions
\[ \overline{\mbf{P}}^*_{s,b} = \mcal{L}\left( \overline{X}^{(1)}_b,\ldots,\overline{X}^{(1)}_{b+n_1(s,b)(s+1)}\right) \text{  and  }\overline{\mbf{P}}^{ind}_{s,b} = \bigotimes_{i=1}^{n_1(s,b)} \mcal{L}\left( \overline{X}^{(1)}_{b+(i-1)(s+1)} \right). \]
The next result is a complement of Lemma \ref{lem:robust_kl} and is proven in Section \ref{sec:proof_lem_robust_beta}.
\begin{lemme}
\label{lem:robust_beta}
If $E^{(1)}_1,Z^{(1)}_1,\ldots,E^{(1)}_{n_1},Z^{(1)}_{n_1},E^{(2)}_1,Z^{(2)}_1,\ldots,E^{(2)}_{n_2},Z^{(2)}_{n_2},\mbf{\overline{X}}^{(1)}$ and $\mbf{\overline{X}}^{(2)}$ are mutually independent, we have
\begin{equation}
\label{eq:robust_beta_1}
\mbf{K}\left( \mbf{P}^*_{s,b} || \mbf{P}^{ind}_{s,b} \right) \leq
\mbf{K}\left(\overline{\mbf{P}}^*_{s,b} || \overline{\mbf{P}}^{ind}_{s,b} \right),\forall s\in\{0,1,\ldots,s_{\max}\},\forall b\in[s+1],
\end{equation}
and
\[ \beta_t\left( \mbf{X}^{(2)} \right) \leq \beta_t\left( \overline{\mbf{X}}^{(2)} \right),\forall t\geq 1. \]
\end{lemme}
We define $p^{(j)}_i$ by $\mbb{P}\left( E^{(j)}_i = 1 \right)=p^{(j)}_i$ for $j\in[2]$ and $i\in[N_j]$.
\begin{coro}
\label{coro:robust_selection_s}
Let $n_1,n_2\geq 3$ and $\hat{P}=\hat{P}_{\hat{s}}\left(\mbf{X}^{(1)},\mbf{X}^{(2)}\right)$ be the estimator given by (\ref{eq:p_hat_s_hat}). There exists a positive constant $C$ such that in the situation of Lemma \ref{lem:robust_beta} and for all $\overline{P}\in\PX$,
\begin{align*}
C \mbb{E}\left[ h^2\left( \overline{P}, \hat{P}_{\hat{s}} \right) \right] &\leq n_1^{-1} \suml_{i=1}^{n_1} (1-p^{(1)}_i) + n_2^{-1} \suml_{i=1}^{n_2} (1-p^{(2)}_i)\\
+ \inf_{t\in[n_2]} &\bigg\{ \frac{t}{n_2} \left( 1 + \log(|S|) \right) + \lceil n_2/t\rceil \beta_t\left( \overline{\mbf{X}}^{(2)} \right) \bigg\}\nonumber\\
+ \inf_{s\in S} &\bigg\{ h^2\left( \overline{P} , \mscr{M}_s \right) + \frac{(s+1) D_{n_1(s,1)}(\mscr{M}_s)}{n_1} + n_1^{-1} \suml_{b=1}^{s+1} \mbf{K}\left( \overline{\mbf{P}}^*_{s,b} || \overline{\mbf{P}}^{ind}_{s,b}\right) \bigg\}.
\end{align*}
\end{coro}
This result is a direct consequence of Theorem \ref{th:rho_selection_s} and Lemma \ref{lem:robust_beta}. We illustrate the performance of our estimator with hidden Markov models.
\subsection{Application to hidden Markov models}
Let $Y_1^{(1)},\dots,Y_{N_1}^{(1)},Y^{(2)}_1,\dots,Y^{(2)}_{N_2}$ be random variables taking values in the measurable space $\left(\mscr{Y},\mcal{Y}\right)$. Let $L$ be in $\{2,3,\dots,\lfloor (N_1\wedge N_2)/2\rfloor \}$ and $n_j=N_j+1-L$ for $j\in[2]$. We define the new random variables
\[ X^{(j)}_i = \left( Y^{(j)}_i,Y^{(j)}_{i+1},\dots,Y^{(j)}_{i+L-1}\right),i\in[n_j],j\in[2], \]
taking values in the measurable space $\left(\mscr{X},\mcal{X}\right) = \left(\mscr{Y}^L,\mcal{Y}^{\otimes L}\right)$. We adapt Assumption \ref{hyp:hmm_ergodic} to this context.
\begin{hypothese}
\label{hyp:hmm_ergodic_s}
Let $\left(Y^{(1)}_i,H^{(1)}_i\right)_i$ and $\left(Y^{(2)}_i,H^{(2)}_i\right)_i$ be finite state space HMM with parameters $(K^*,w_1^*,Q^*,F^*)$ and $(K^*,w_2^*,Q^*,F^*)$ such that $Q^*$ is irreducible and aperiodic.
\end{hypothese}
Under this assumption $Q^*$ has only one invariant distribution $\pi^*$ and we define the distribution $P^*$ by (\ref{eq:p_star}). Let $\tau\geq e$ and $J=\left\lfloor \log_{\tau}\left( \lfloor (n_1-2)/2\rfloor \right) \right\rfloor$. Let $S$ be the set given by
\begin{equation}
\label{eq:S}
S = \{ 0 \} \cup \left\{ \left\lceil \tau^j \right\rceil ; j\in \{ 0,1,\dots,J\} \right\}.
\end{equation}
Let $\overline{\mscr{F}}_1,\ldots,\overline{\mscr{F}}_K$ be subsets of $\mscr{P}_Y$ such that Assumption \ref{hyp:hmm_vc} is satisfied. Let $\overline{\mscr{M}}$ be a non-empty subset of the model $\mscr{H}\left(K,\overline{\mscr{F}}_1,\dots,\overline{\mscr{F}}_K\right)$ defined by (\ref{eq:model_h_hmm}). For $s$ in $S$, we take $\mscr{M}_s=\mscr{M}_{\delta(s)}$ with 
\[ \delta(s)=\frac{\overline{V}}{n_1(s,1)(K-1)}\bigwedge \frac{1}{K}, \]
where $\mscr{M}_{\delta}$ is given by (\ref{eq:m_delta}) and $n_1(s,1)$ given by (\ref{eq:n_1_s_b}). The following result is proven in Section \ref{sec:proof_th_selection_s_hmm}.
\begin{theoreme}
\label{th:selection_s_hmm}
Let $N_1,N_2 \geq K+L$ and $\hat{P}=\hat{P}_{\hat{s}}\left(\mbf{X}^{(1)},\mbf{X}^{(2)}\right)$ be the estimator given by (\ref{eq:p_hat_s_hat}). Under Assumption \ref{hyp:selection_s_independent}, there is a numeric constant $C>0$ such that for all $\overline{P}\in\PX$
\begin{align}
C \mbb{E}\left[ h^2\left(\overline{P}, \hat{P}_{\hat{s}}\right) \right] &\leq h^2\left(\overline{P},\overline{\mscr{M}}\right) + n_1^{-1} \suml_{i=1}^{n_1} h^2\left( P_i^{(1)},\overline{P} \right) + n_2^{-1} \suml_{i=1}^{n_2} h^2\left( P_i^{(2)},\overline{P} \right)\nonumber\\
&+ L\epsilon^2 + \inf_{t\in[n_2]} \bigg\{ \frac{t \log \log n_1}{n_2} + \lceil n_2/t\rceil \beta_t\left( \mbf{X}^{(2)} \right) \bigg\}\label{eq:th_selection_s_hmm}\\
+ \inf_{s\in S} &\bigg\{ \frac{(s+1) L \overline{V} \log n_1}{n_1} + n_1^{-1} \suml_{b=1}^{s+1} \mbf{K}\left(\mbf{P}^*_{s,b}||\mbf{P}^{ind}_{s,b}\right) \bigg\}\nonumber.
\end{align}
In particular under Assumption \ref{hyp:hmm_ergodic_s}, there exists a positive constant $C(Q^*)$ such that
\begin{equation}
C(Q^*) \mbb{E}\left[ h^2\left( P^*, \hat{P}_{\hat{s}} \right) \right] \leq h^2\left( P^*,\overline{\mscr{M}} \right) + L \epsilon^2 + \tau L \overline{V} \frac{\log^2 n_1}{n_1} + \frac{\log n_2 \log \log n_1}{n_2},\label{eq:th_selection_s_hmm_2}
\end{equation}
where $P^*$ is given by (\ref{eq:p_star}).
\end{theoreme}
Inequality (\ref{eq:th_selection_s_hmm}) is a consequence of Theorem \ref{th:rho_selection_s} and only requires Assumption \ref{hyp:selection_s_independent}. Under Assumption \ref{hyp:hmm_ergodic_s} we can control the different terms and obtain (\ref{eq:th_selection_s_hmm_2}). If $\epsilon=0$, the ideal situation is to have the same number of observations in each set, i.e. $n_1=n_2=n$. In this case we have
\[ C\left(Q^*\right) \mbb{E}\left[ h^2\left( P^* , \hat{P} \right) \right] \leq h^2\left( P^*,\overline{\mscr{M}} \right) + L \tau \overline{V} \frac{\log^2 n}{n}, \]
and the first vanishes when the model is well specified which gives the rate $n^{-1} \log^2 n$ with respect to the squared Hellinger distance over $\mscr{H}^*\left(K,\overline{\mscr{F}}_1,\dots,\overline{\mscr{F}}_K\right)$. When $\epsilon>0$ the quantity $\overline{V}$ depends on $\epsilon$ and we need to balance the second and third term in (\ref{eq:th_selection_s_hmm_2}), i.e. $\epsilon^2/\overline{V}$ is of order $n_1^{-1}$ up to a logarithmic term. Then the ideal situation only requires $n_2$ to be of order $\epsilon^{-2}$ up to logarithmic term and the bound on the convergence rate is of order $\epsilon^2$. For example, we would have $\epsilon^{-2}=n_1^{\frac{2}{d+1}} \log^{-\frac{2}{d+1}-(d+2)}n_1$ in the situation of Theorem \ref{th:hmm_log_concave_high}. In both cases, it shows that we recover a value of $s$ that allows to obtain the same rate as when the optimal value is known. This is especially interesting for the robustness aspect of our estimator.\par
Let us consider a situation similar to Section \ref{sec:robust_selection_s}. Let $Z^{(1)}_1,\ldots,Z^{(1)}_{N_1},Z^{(2)}_1,\ldots,Z^{(2)}_{N_2}$ be random variables with any distributions and $E^{(1)}_1,\ldots,E^{(1)}_{N_1},E^{(2)}_1,\ldots,E^{(2)}_{N_2}$ be Bernoulli random variables such that for all $j\in[2]$ and all $i\in[N_j]$,
\[ Y^{(j)}_i = E_i \overline{Y}^{(j)}_i + (1-E^{(j)}_i) Z^{(j)}_i. \]
The following result is proven in Section \ref{sec:proof_coro_robust_hmm_selection_s}.
\begin{coro}
\label{coro:robust_hmm_selection_s}
Let $\hat{P}_{\hat{s}}=\hat{P}_{\hat{s}}\left(\mbf{X}^{(1)},\mbf{X}^{(2)}\right)$ be the estimator given by (\ref{eq:p_hat_s_hat}). If $E^{(1)}_1,Z^{(1)}_1,\ldots,E^{(1)}_{n_1},Z^{(1)}_{n_1}$, $E^{(2)}_1,Z^{(2)}_1,\ldots,E^{(2)}_{n_2},Z^{(2)}_{n_2},\mbf{\overline{X}}^{(1)}$ and $\mbf{\overline{X}}^{(2)}$ are mutually independent, and if $\mbf{\overline{Y}}^{(1)}$ and $\mbf{\overline{Y}}^{(2)}$ satisfy Assumption \ref{hyp:hmm_ergodic_s}, there exists a positive constant $C(Q^*)$ such that
\begin{align*}
C(Q^*) \mbb{E}\left[ h^2\left( P^*, \hat{P}_{\hat{s}} \right) \right] &\leq \frac{L}{N_1} \suml_{i=1}^{N_1} \left( 1-p^{(1)}_i \right) + \frac{L}{N_2} \suml_{i=1}^{N_2} \left( 1-p^{(2)}_i \right)\\
&+ L \epsilon^2 + \tau L \overline{V} \frac{\log^2 n_1}{n_1} + \frac{\log n_2 \log \log n_1}{n_2},
\end{align*}
where $P^*$ is given by (\ref{eq:p_star}) and $p^{(j)}_i=\mbb{P}\left(E^{(j)}_i=1\right)$ for all $j\in[2]$ and $i\in[N_j]$.
\end{coro}
One can see that our deviation bound is not significantly worse as long as the average proportions of contamination $N_1^{-1} \sum_{i=1}^{N_1} (1-p^{(1)}_i)$ and $N_2^{-1} \sum_{i=1}^{N_2} (1-p^{(2)}_i)$ are small compared to $\epsilon^2 + \tau \overline{V} \frac{\log^2 n_1}{n_1}$ and $\frac{\log n_2 \log \log n_1}{n_1}$ respectively. We interpret this result further for $\epsilon^2$ and $n_1=n_2=n$. Let us consider H\"{u}ber's contamination model with $p^{(j)}_i=1-\alpha_{cont}$ for all $j\in[2]$ and $i\in[N]$. In this situation we get
\begin{equation*}
C(Q^*) \mbb{E} \left[ h^2\left( P^*, \hat{P}_s \right) \right] \leq L \left[ \alpha_{cont} + \frac{\tau \overline{V} \log^2 n}{n} \right].
\end{equation*}
Our bound on the convergence rate is not deteriorated as long as the contamination rate $\alpha_{cont}$ is small compared to $\epsilon^2 + \frac{\tau \overline{V} \log^2 n}{n}$. We can also consider the situation $\mbb{P}(E^{(j)}_i=0)=\mathbbm{1}_{i\in I_j}$ for some subsets $I_1\subset[N]$ and $I_2\subset[N]$. We get
\begin{equation*}
C(Q^*) \mbb{E} \left[ h^2\left( P^*, \hat{P}_s \right) \right] \leq L \left[ \frac{|I_1|+|I_2|}{N} + \frac{\tau \overline{V} \log^2 n}{n} \right].
\end{equation*}
Our bound on the convergence rate is not deteriorated as long as the proportions of outliers $|I_1|/N,|I_2|/N$ are small compared to the other terms. 
\printbibliography
\newpage
\appendix
\section{Auxiliary results}
\label{sec:first_step_appendix}
We denote by $\text{C}(\mscr{X})$ the set given by
\[ \text{C}(\mscr{X})=\bigcup_{n\geq1} \{n\}\times \mscr{X}^n. \]
Let $d:\mscr{A}\times\mscr{A}\rightarrow\mbb{R}$ be a loss function where $\mscr{A}\subset\mscr{P}_X$ denotes a set of admissible probability distributions. Let $\mscr{M}$ be a subset of $\mscr{A}$. Let $\hat{P}: C\left(\mscr{X}\right) \rightarrow \mscr{M}$ be an estimation method.
\begin{hypothese}
\label{hyp:estimator_oracle}
There exist constants $C_0>0,\beta\in(0,1]$ and non decreasing functions $f,g$ such that for all independent random variables $X_1,\dots,X_n$ with distributions $P_1,\dots,P_n\in\mscr{A}$ and for all $\xi>0$
\[ \mbb{P} \left( \suml_{i=1}^n d\left( P_i,\hat{P}\left(n,\mbf{X}\right) \right) \leq  C_0 \inf_{Q\in\mscr{M}} \suml_{i=1}^n d\left( P_i,Q \right) + f(n) + g(n) \xi^{\beta} \right) \geq 1-e^{-\xi}. \]
\end{hypothese}
Many estimators satisfy such an assumption, see for instance mean discrepancy estimators \cite{alquier1}, $T$-estimators \cite{testimator} or $l$-estimators \cite{lestimator}. We can get rid of the independence assumption with the following result.
\begin{prop}
\label{prop:motivation}
Under Assumption \ref{hyp:estimator_oracle}, for all random variables $X_1,\dots,X_n$ with distributions $P_1,\dots,P_n\in\mscr{A}$ we have
\begin{align*}
\mbb{E}\left[ \suml_{i=1}^n d\left( P_i,\hat{P}\left(n,\mbf{X}\right) \right) \right] &\leq C_0 \inf_{Q\in\mscr{Q}} \suml_{i=1}^n d(P_i,Q) + f(n)\\
&+ g(n) \left[ 2 + \frac{3}{2} \mbf{K}\left(\mbf{P}^*||\mbf{P}^{ind}\right) \right]^{\beta},
\end{align*}
where
\[ \mbf{P}^* = \mcal{L}\left( X_1,\ldots,X_n \right) \text{  and  } \mbf{P}^{ind} = \mcal{L}\left(X_1\right) \otimes \ldots \otimes \mcal{L}\left(X_n\right). \]
\end{prop}
This result is obtained by applying Lemma \ref{lem:kldiv_exponential} that we prove hereafter, with $\mbf{P}=\mbf{P}^{ind}$ and $\mbf{Q}=\mbf{P}^*$.
\subsection{Proof of Lemma \ref{lem:kldiv_exponential}}
\label{sec:proof_lem_kldiv_exponential}
We use Lemma 48 in \cite{baraudinventiones}. For $\lambda\in(0,a^{-1/\beta})$, we have
\begin{align*}
&\mbb{E}_{\mbf{Q}} \left[ \lambda \left(  n l\left(\hat{\theta}(\mbf{X}),\theta\right) - n A - B \right)_+^{1/\beta} \right]\\
&\leq \log\left( 1 + \int_0^{+\infty} e^{\xi} \mbf{P}\left( l\left(\hat{\theta}(\mbf{X}),\theta\right) > A + \frac{B+(\xi/\lambda)^{\beta}}{n} \right) d\xi \right) + \mbf{K}\left(\mbf{Q}||\mbf{P}\right)\\
&\leq \log\left( 1 + \int_0^{+\infty} e^{\xi} e^{- \xi/\lambda} d\xi \right) + \mbf{K}\left(\mbf{Q}||\mbf{P}\right) = \log\left(\frac{1}{1-\lambda} \right) + \mbf{K}\left(\mbf{Q}||\mbf{P}\right).
\end{align*}
We have
\begin{align*}
\mbb{E}_{\mbf{Q}} \left[ \left(  n l\left(\hat{\theta}(\mbf{X}),\theta\right) - n A - B \right)_+^{1/\beta} \right] \leq \lambda^{-1} \left[ \log\left(\frac{1}{1-\lambda} \right) + \mbf{K}\left(\mbf{Q}||\mbf{P}\right) \right].
\end{align*}
Assuming $\mbf{K}(\mbf{Q}||\mbf{P})<\infty$, minimization over $\lambda$ demands
\[ \log\left(1-\lambda\right) - \mbf{K}\left(\mbf{Q}||\mbf{P}\right) + \frac{\lambda}{1-\lambda} = 0. \]
Let $\lambda^*$ be such a number. In that case
\[ (\lambda^*)^{-1} \left[ \log\left( \frac{1}{1-\lambda^*} \right) + \mbf{K}\left(\mbf{Q}||\mbf{P}\right) \right] = \frac{1}{1-\lambda^*}. \]
We set $a(x)=x-\log(1 +x)$ for $x$ in $(0,+\infty)$. Following the proof of Proposition 5 \cite{baraudinventiones}, $a$ is increasing and
\[\forall x>0, a^{-1}(x)\leq x+\sqrt{2x}.\]
Since $\frac{\lambda^*}{1-\lambda^*} = a^{-1}\left( \mbf{K}\left(\mbf{Q}||\mbf{P}\right) \right)$, we get
\begin{align*}
\frac{1}{1-\lambda^*} = 1 + \frac{\lambda^*}{1-\lambda^*} &\leq 1 + \mbf{K}\left(\mbf{Q}||\mbf{P}\right) + \sqrt{2 \mbf{K}\left(\mbf{Q}||\mbf{P}\right) }\\
&\leq 2 + \frac{3}{2} \mbf{K}\left(\mbf{Q}||\mbf{P}\right). 
\end{align*}
Finally, with Jensen's inequality we get
\begin{align*}
\mbb{E}_{\mbf{Q}} \left[ l\left( \hat{\theta}(\mbf{X}),\theta \right) \right] &\leq A + \frac{B + \left( 2 + \frac{3}{2} \mbf{K}\left(P||Q\right) \right)^{\beta}}{n}.
\end{align*}
\section{Main results}
\label{sec:rho_application}
This section gathers the proofs of Theorem \ref{th:main_rho}, Corollary \ref{coro:robust} and Lemmas \ref{lem:robust_kl}, \ref{lem:divergence_markov}, \ref{lem:divergence_hmm}. We first give a formal definition of the $\rho$-dimension function that is originally introduced in Baraud \& Birgé \cite{baraudrevisited}. We slightly modify some notation to adapt it to our context. The function $\psi$ defined by (\ref{eq:psi_rho}) satisfies Assumption 2 \cite{baraudrevisited} with $a_0 = 4, a_1 = 3/8$ and $a_2^2 = 3\sqrt{2}$ (see Proposition 3 \cite{baraudrevisited}). Let $n$ be a positive integer and $\mscr{M}$ be a countable subset of $\mscr{P}_X$. For $y>0$, $\mbf{P}^{ind}= \bigotimes_{i=1}^n P^{ind}_1 \in\mscr{P}_X^{\otimes n}$ and $P\in\mscr{M}$ we write
\[ \mscr{B}^{\mscr{M}}\left(\mbf{P}^{ind},\overline{P},y\right) := \left\{ Q\in\mscr{M}; \suml_{i=1}^n h^2\left(P^{ind}_i,P\right) + h^2\left(P^{ind}_i,Q\right) < y^2 \right\}. \]
If $\mcal{M}$ is a countable set of probability density functions with respect to a $\sigma$-finite measure $\nu$ such that $\mscr{M}=\{ Q = q\cdot \nu; q\in\mcal{M}\}$, we write
\[ w\left(\nu,\mcal{M},\mscr{M},\mbf{P}^{ind},P,y\right) = \mbb{E}_{\mbf{X}\sim\mbf{P}^{ind}} \left[ \sup\limits_{Q\in\mscr{B}^{\mscr{M}}(\mbf{P}^{ind},P,y)} \left| \mbf{Z}_n\left(\mbf{X},p,q\right) \right| \right], \]
where
\[ \mbf{Z}_n(\mbf{X},q,q') := \mbf{T}_n(\mbf{X},q,q') - \mbb{E}_{\mbf{P}^{ind}} \mbf{T}_n(\mbf{X},q,q'),\]
and $\mbf{T}_n$ is given by (\ref{eq:t_n}). We define $\mbf{w}^{\mscr{M}}\left(\mbf{P}^{ind},P,y\right) = \inf_{(\nu,\mcal{M})} w\left(\nu,\mcal{M},\mscr{M},\mbf{P}^{ind},P,y\right)$, where the infimum is taken over all couples $(\nu,\mcal{M})$ such that $\mcal{M}$ is the class of density functions associated to $\mscr{M}$ with respect to a $\sigma$-finite measure $\nu$. We define the $\rho$-dimension function by 
\[ D^{\mscr{M}}\left(\mbf{P}^{ind},P^{\otimes n} \right) = \left[ \frac{3}{2^{21/2}} \sup\left\{ y^2 ;\mbf{w}^{\mscr{M}}\left(\mbf{P}^{ind},P,y\right)>\frac{3 y^2}{64} \right\} \right] \bigvee 1. \]
As mentioned at the beginning of Section \ref{sec:first_section} we consider cases for which we have a uniform bound over the $\rho$-dimension function. More precisely we assume there is a non-increasing function $m\mapsto D_m(\mscr{M})$ such that
\[ D^{\mscr{M}}\left(\mbf{P}^{ind},P^{\otimes m} \right) \leq D_m(\mscr{M}), \forall \mbf{P}^{ind}\in \PX^{\otimes m}, \forall P\in\mscr{M}. \]
\subsection{Proof of Theorem \ref{th:main_rho}}
\label{sec:proof_th_main_rho}
From Theorem 1 of Baraud \& Birgé \cite{baraudrevisited}, we have that for all independent random variables $X_1,\dots,X_n$ with respective distributions $P_1,\dots,P_n$, for all $Q\in\mscr{M}$ and for all $\xi>0$, we have
\begin{align*}
\suml_{i=1}^n  h^2\left( P_i, \hat{P}(n,\mbf{X},\mscr{M}) \right) \leq \gamma \suml_{i=1}^n h^2\left( P_i, Q \right) + \frac{4\kappa}{a_1} \left( \frac{D_n(\mscr{M})}{4.7} + 1.49 + \xi \right),
\end{align*}
with probability at least $1-e^{-\xi}$, where $\gamma$ and $\kappa$ are given in \cite{baraudrevisited} and satisfy $\gamma\leq 150$ and $\frac{4\kappa}{a_1} \leq 5014$ (see proof of Theorem 1 \cite{baraud_glm}, page 32). We can take the infimum for $Q$ over $\mscr{M}$ and it shows we satisfy Assumption \ref{hyp:estimator_oracle} with $C_0=150$, $f(n) = 5014 \left( \frac{D_n(\mscr{M}))}{4.7} + 1.49 \right)$, $g(n)=5014$ and $\beta=1$. From Proposition \ref{prop:motivation}, we have
\begin{align*}
\mbb{E}\left[ \suml_{i=1}^{n(s,b)} h^2\left( P_{b+(i-1)(s+1)}, \hat{P}_s \right) \right] &\leq 150 \inf_{Q\in\mscr{Q}} \suml_{i=1}^{n(s,b)} h^2\left( P_{b+(i-1)(s+1)},Q\right)\\
&+ 5014 \left( \frac{D_{n(s,b)}(\mscr{M})}{4.7} + 3.49 + \frac{3}{2}\mbf{K}\left(\mbf{P}_{s,b}^*||\mbf{P}_{s,b}^{ind}\right) \right),
\end{align*}
for all $b\in[s+1]$. From (\ref{eq:p_hat_s}), we have
\begin{align*}
\suml_{i=1}^n h^2\left( P_i, \hat{P}_s \right) &= \suml_{b=1}^{s+1} \suml_{i=1}^{n(s,b)} h^2\left( P_{b+(i-1)(s+1)}, \hat{P}_s \right)\\
&\leq 2 \suml_{b=1}^{s+1} \suml_{i=1}^{n(s,b)} h^2\left( P_{b+(i-1)(s+1)}, \hat{P}_{s,b} \right) + 2 \suml_{b=1}^{s+1} n(s,b) h^2\left( \hat{P}_{s,b}, \hat{P}_s \right)\\
&\leq 2 \suml_{b=1}^{s+1} \suml_{i=1}^{n(s,b)} h^2\left( P_{b+(i-1)(s+1)}, \hat{P}_{s,b} \right) + 2 \inf\limits_{Q\in\mscr{M}} \suml_{b=1}^{s+1} n(s,b) h^2\left( \hat{P}_{s,b}, Q \right) + 2 \iota\\
&\leq 4 \suml_{b=1}^{s+1} \suml_{i=1}^{n(s,b)} h^2\left( P_{b+(i-1)(s+1)}, \hat{P}_{s,b} \right) + 2 \inf\limits_{Q\in\mscr{M}} \suml_{i=1}^N h^2\left( P_i, Q \right) + 2 \iota.
\end{align*}
Combining the inequalities above, we obtain
\begin{align*}
\mbb{E}\left[ \suml_{i=1}^n h^2\left( P_i, \hat{P}_s \right)\right] &\leq 600 \suml_{b=1}^{s+1} \inf_{Q\in\mscr{M}} \suml_{i=1}^{n(s,b)} h^2(P_{b+(i-1)(s+1)},Q) + 2 \inf\limits_{Q\in\mscr{M}} \suml_{i=1}^n l\left( P_i, Q \right)\\
&+ 20 056 \suml_{b=1}^{s+1} \left( \frac{D_{n(s,b)}(\mscr{M})}{4.7} + 3.49 + \frac{3}{2}\mbf{K}\left(\mbf{P}_{s,b}^*||\mbf{P}_{s,b}^{ind}\right) \right) + 2\iota\\
&\leq 602 \inf_{Q\in\mscr{M}} \suml_{i=1}^n h^2(P_i,Q) + 20 056 (s+1) \left( \frac{D_{n(s,1)}(\mscr{M})}{4.7} + 3.49 \right)\\
&+ 30 084 \suml_{b=1}^{s+1} \mbf{K}\left(\mbf{P}_{s,b}^*||\mbf{P}_{s,b}^{ind}\right) + 2\iota.
\end{align*}
Since $\iota \leq 2546 < 20 056 \times \frac{0.597}{4.7}$, we get
\begin{align*}
\mbb{E}_{\mbf{P}^*} \left[ \suml_{i=1}^n h^2\left( P_i, \hat{P}_s \right) \right] &\leq 602 \inf\limits_{Q\in\mscr{M}} \suml_{i=1}^n h^2\left( P_i, Q \right) + \frac{20 056}{4.7} (s+1) \left[ D_{n(s,1)}(\mscr{M}) + 17 \right]\\
&+ 30 084 \suml_{b=1}^{s+1} \mbf{K}\left(\mbf{P}^*_{s,b}||\mbf{P}^{ind}_{s,b}\right).
\end{align*}
\subsection{Proof of Lemma \ref{lem:robust_kl}}
\label{sec:proof_lem_robust_kl}
For $\mbf{e}\in\{0,1\}^n$, we denote by $I(\mbf{e})$ the set given by $I(\mbf{e})=\{i\in[n]; e_i=1\}$. From the convexity property of the Kullback-Leibler divergence, we have
\begin{align*}
\mbf{K}&\left( \mcal{L}\left(\mbf{Y}\right) || \mcal{L}(Y_1) \otimes \dots \otimes \mcal{L}(Y_n) \right)\\
&\leq \suml_{\mbf{e}\in\{0,1\}^n} \mbb{P}(\mbf{E}=\mbf{e}) \mbf{K}\left( \mcal{L}\left(\mbf{Y}|\mbf{E}=\mbf{e}\right) || \mcal{L}(Y_1|E_1=e_1) \otimes \dots \otimes \mcal{L}(Y_n|E_n=e_N) \right)\\
&= \suml_{\mbf{e}\in\{0,1\}^n} \mbb{P}(\mbf{E}=\mbf{e}) \mbf{K}\left( \mcal{L}\left( \left(X_i\right)_{i\in I(\mbf{e})}\right) \otimes \bigotimes_{i\notin I(\mbf{e})}  \mcal{L}\left( Z_i\right) || \bigotimes_{i\in I(\mbf{e})} \mcal{L}(X_i) \otimes \bigotimes_{i\notin I(\mbf{e})} \mcal{L}(Z_i) \right)\\
&= \suml_{\mbf{e}\in\{0,1\}^n} \mbb{P}(\mbf{E}=\mbf{e}) \mbf{K}\left( \mcal{L}\left( \left(X_i\right)_{i\in I(\mbf{e})}\right) || \bigotimes_{i\in I(\mbf{e})} \mcal{L}(X_i) \right).
\end{align*}
We need an auxiliary result before ending the proof.
\begin{lemme}
\label{lem:kullback_formula}
For random variables $A,B,C$ such that $\mcal{L}(A)\ll\mcal{L}(B)$, we have
\begin{equation}
\label{eq:kl_indep}
\mbf{K}\left( \mcal{L}(A) || \mcal{L}(B) \right) \leq \mbf{K}\left( \mcal{L}\left(A,C\right) || \mcal{L}(B) \otimes \mcal{L}(C) \right).
\end{equation}
\end{lemme}
With this result we have
\begin{align*}
\mbf{K}\left( \mcal{L}\left( \left(X_i\right)_{i\in I(\mbf{e})}\right) || \bigotimes_{i\in I(\mbf{e})} \mcal{L}(X_i) \right) &\leq \mbf{K}\left( \mcal{L}\left(\mbf{X}\right) || \mcal{L}(X_1) \otimes \dots \otimes \mcal{L}(X_n) \right),
\end{align*}
which allows to conclude.
\subsubsection{Proof of Lemma \ref{lem:kullback_formula}}
\label{sec:lem_kullback_formula}
Let $\mu_1$ and $\mu_2$ be measures dominating $\mcal{L}(B)$ and $\mcal{L}(C)$ respectively. We write
\[ p_{B,C}=\frac{d\mcal{L}(B,C)}{d\mu_1\otimes\mu_2}, p_{A,C}=\frac{d\mcal{L}(A,C)}{d\mu_1\otimes\mu_2}, p_A=\frac{d\mcal{L}(A)}{d\mu_1},p_B=\frac{d\mcal{L}(B)}{d\mu_1}, p_C=\frac{d\mcal{L}(C)}{d\mu_2}. \]
We have
\begin{align*}
\mbf{K}\left( \mcal{L}(A,C) || \mcal{L}(B) \otimes \mcal{L}(C) \right) &= \int p_{A,C}(x,z) \log\left(\frac{ p_{A,C}(x,z) }{ 
p_B(x) p_C(z) } \right) \mu_1(dx) \mu_2(dz)\\
&= \int p_{A,C}(x,z) \log\left(\frac{ p_{A,C}(x,z) }{ 
p_A(x) p_C(z) } \right) \mu_1(dx) \mu_2(dz)\\
&+ \int p_{A,C}(x,z) \log\left(\frac{ p_A(x) }{ 
p_B(x) } \right) \mu_1(dx) \mu_2(dz)\\
&=  \mbf{K}\left( \mcal{L}(A,C) || \mcal{L}(A) \otimes \mcal{L}(C) \right) + \mbf{K}\left( \mcal{L}(A) || \mcal{L}(B) \right).
\end{align*}
The non-negativity of the Kullback-Leibler divergence concludes the proof.
\subsection{Proof of Corollary \ref{coro:robust}}
\label{sec:proof_coro_robust}
One can check that we have
\begin{align*}
h^2\left(\overline{P},\hat{P}_s\right) &\leq 2 n^{-1} \suml_{i=1}^n h^2\left(\mcal{L}(Y_i),\overline{P}\right) + 2 n^{-1} \suml_{i=1}^n h^2\left(\mcal{L}(Y_i),\hat{P}_s\right)\\
&\leq 2 n^{-1} \suml_{i=1}^n (1-p_i) + 2 n^{-1} \suml_{i=1}^n h^2\left(\mcal{L}(Y_i),\hat{P}_s\right),
\end{align*}
and for $Q$ in $\mscr{M}$
\begin{align*}
\suml_{i=1}^n h^2\left( \mcal{L}(Y_i),Q \right) &\leq 2 \suml_{i=1}^n h^2\left( \mcal{L}(Y_i),\overline{P} \right) + 2 \suml_{i=1}^n h^2\left(\overline{P},Q\right)\\
&\leq 2 \suml_{i=1}^n (1-p_i) + 2 n h^2\left(\overline{P},Q\right).
\end{align*}
We can conclude with Theorem \ref{th:main_rho} and Lemma \ref{lem:robust_kl}.
\subsection{Proof of Lemma \ref{lem:divergence_markov}}
\label{sec:proof_lem_divergence_markov}
We have
\begin{align*}
\mbf{K}\left( \mcal{L}\left(\mbf{X}\right) || \mcal{L}\left(X_1\right)\otimes\dots\otimes \mcal{L}\left(X_n\right) \right) &= \mbb{E}\left[ \mbf{K}\left( \mcal{L}\left(X_n|X_1,\dots,X_{n-1}\right) || \mcal{L}\left(X_n\right) \right) \right]\\
&+ \mbf{K}\left( \mcal{L}\left(X_1,\dots,X_{n-1}\right) || \mcal{L}\left(X_1\right)\otimes\dots\otimes \mcal{L}\left(X_{n-1}\right) \right),
\end{align*}
and with the Markov property
\begin{align*}
\mbb{E}\left[ \mbf{K}\left( \mcal{L}\left(X_n|X_1,\dots,X_{n-1}\right) || \mcal{L}\left(X_n\right) \right) \right] &= \mbb{E}\left[ \mbf{K}\left( \mcal{L}\left(X_n|X_{n-1}\right) || \mcal{L}\left(X_n\right) \right) \right]\\
&= \mbf{K}\left( \mcal{L}\left(X_{n-1},X_n\right) || \mcal{L}\left(X_{n-1}\right)\otimes \mcal{L}\left(X_n\right) \right).
\end{align*}
Therefore
\begin{align*}
\mbf{K}\left( \mcal{L}\left(\mbf{X}\right) || \mcal{L}\left(X_1\right)\otimes\dots\otimes \mcal{L}\left(X_n\right) \right) &= \mbf{K}\left( \mcal{L}\left(X_1,\ldots,X_{n-1}\right) || \mcal{L}\left(X_1\right)\otimes\dots\otimes \mcal{L}\left(X_{n-1}\right) \right)\\
&+ \mbf{K}\left( \mcal{L}\left(X_{n-1},X_n\right) || \mcal{L}\left(X_{n-1}\right)\otimes \mcal{L}\left(X_n\right) \right),
\end{align*}
and we can conclude by induction.
\subsection{Proof of Lemma \ref{lem:divergence_hmm}}
\label{sec:proof_lem_divergence_hmm}
If $\left(\mbf{X},\mbf{H}\right)$ a hidden Markov chain, with Lemma \ref{lem:divergence_markov} we have
\begin{align*}
&\mbf{K}\left( \mcal{L}\left(\mbf{X}\right) || \mcal{L}\left( X_1 \right) \otimes \dots \otimes \mcal{L}\left( X_n \right) \right)\\
&\leq \suml_{i=2}^n \mbf{K}\left( \mcal{L}\left(X_{i-1},H_{i-1},X_i,H_i \right) || \mcal{L}\left(X_{i-1},H_{i-1}\right) \otimes \mcal{L}\left(X_i,H _i\right) \right). 
\end{align*}
We need the following result. For random variables $A_1,A_2,B_1,B_2$, we have
\begin{align*}
&\mbf{K}\left( \mcal{L}\left(A_1,B_1,A_2,B_2 \right) || \mcal{L}\left(A_1,B_1 \right) \otimes \mcal{L}\left(A_2,B_2 \right) \right)\\
&= \mbf{K}\left( \mcal{L}\left(A_1,A_2 \right) || \mcal{L}\left(A_1 \right) \otimes \mcal{L}\left(A_2 \right) \right)\\
&+ \mbb{E}\left[ \mbf{K}\left( \mcal{L}\left(B_1,B_2|A_1,A_2 \right) || \mcal{L}\left(B_1|A_1 \right) \otimes \mcal{L}\left(B_2|A_2 \right) \right) \right].
\end{align*}
With the non-negativity of the Kullback-Leibler divergence we get
\[ \mbf{K}\left( \mcal{L}\left(\mbf{X}\right) || \mcal{L}\left( X_1 \right) \otimes \dots \otimes \mcal{L}\left( X_n \right) \right) \leq \suml_{i=2}^n \mbf{K}\left( \mcal{L}\left(H_{i-1},H_i \right) || \mcal{L}\left(H_{i-1}\right) \otimes \mcal{L}\left(H_i\right) \right). \]
\section{Kolmogorov processes}
\label{sec:proof_kolmogorov}
This section gathers the proofs of Theorems \ref{th:kolmogorov}, \ref{th:kolmogorov_high} and Lemmas \ref{lem:kolmogorov_mixing}, \ref{lem:entropy_log_concave}.
\subsection{Proof of Theorems \ref{th:kolmogorov} and \ref{th:kolmogorov_high}}
\label{sec:proof_th_kolmogorov}
From Proposition 6 \cite{baraudrevisited}, we can take $D_n(\mscr{F}_{\lambda_-,\lambda_+,M}[\epsilon]) = 9 \log( 2 |\mscr{F}_{\lambda_-,\lambda_+,M}[\epsilon]| )$. From Theorem \ref{th:main_rho} there exists a positive constant $C$ such that
\begin{align*}
C \mbb{E}_{\mbf{P}^*} \left[ h^2\left( \overline{P}, \hat{P}_s \right) \right] &\leq h^2\left( \overline{P}, \mscr{F}_{\lambda_-,\lambda_+,M} \right) + \epsilon^2 + n^{-1}\suml_{i=1} h^2\left(P_i,\overline{P} \right)\\
&+ n^{-1} \mbf{K}\left(\mbf{P}^*_{s,b} || \mbf{P}^{ind}_{s,b} \right)\\
&+ \frac{s+1}{n} \left[ 1 + \log(2|\mcal{F}_{\lambda_-,\lambda_+,M}[\epsilon]|) \right].
\end{align*}
Given the bounds on $\log(2|\mcal{F}_{\lambda_-,\lambda_+,M}[\epsilon]|)$ given by Lemma \ref{lem:entropy_log_concave}, we obtain the following inequalities.
\begin{itemize}
\item For $d=1$ we have $\epsilon^2 = n^{-4/5} \log^{4/5} n$ and
\begin{align*}
\log(2|\mcal{F}_{\lambda_-,\lambda_+,M}[\epsilon]|) &\leq \log(9/\eta_1) + \frac{7}{2}\log M + \overline{K}_1 \epsilon^{-1/2}\\
&= \log(9/\eta_1) + \frac{9}{2} \overline{K}_1 n^{1/5} \log^{-1/5} n.
\end{align*}
\item For $d=2$ we have $\epsilon^2 = n^{-2/3} \log ^{5/3} n$ and
\begin{align*}
\log( 2 |\mcal{F}_{\lambda_-,\lambda_+,M}[\epsilon]| ) &\leq \log\left( \frac{3^8\pi}{\eta_2^3} \right) + 9 \log M + \overline{K}_2 \epsilon^{-1}\log_{++}^{3/2}( 1 / \epsilon )\\
&\leq \log\left( \frac{3^8\pi}{\eta_2^3} \right) + \frac{28}{3} \overline{K}_2 n^{1/3}\log^{2/3} n.
\end{align*}
\item For $d=3$ we have $\epsilon^2 = n^{-1/4} \log^{1/4} n$ and
\begin{align*}
\log( 2 |\mcal{F}_{\lambda_-,\lambda_+,M}[\delta]| ) &\leq \log\left( \frac{2^73^{27/2}\pi^3}{\eta_3^6} \right) + \frac{33}{2} \log M + \overline{K}_3 \epsilon^{-2}\\
&= \log\left( \frac{2^73^{27/2}\pi^3}{\eta_3^6} \right) + \frac{33}{2} \overline{K}_3 n^{1/2} \log^{-1/2} n.
\end{align*}
\end{itemize}
This proves the bound (\ref{eq:th_kolmogorov}). Lemma \ref{lem:kolmogorov_mixing} allows to conclude the proof of Theorem \ref{th:kolmogorov}.\par
For $d\geq 4$ we have $\epsilon^2 = n^{-\frac{2}{d+1}} \log^{d+2+\frac{2}{d+1}} n$ and
\begin{align*}
\log(|\mcal{F}_{\lambda_-,\lambda_+,M}[\epsilon]|) &\leq \log C_d + \left( \overline{K}_d + 2 + \frac{1}{d}+\frac{1}{d^2}\right)\epsilon^{-(d-1)} \log^{(d+1)(d+2)/2}( \epsilon^{-1} )\\
&\leq \log C_d + \frac{1}{d+1} \left( \overline{K}_d + 2 + \frac{1}{d}+\frac{1}{d^2}\right) n^{\frac{d-1}{d+1}} \log^{\frac{2}{d+1}+d+1} n.
\end{align*}
Lemma \ref{lem:kolmogorov_mixing} allows to conclude the proof of Theorem \ref{th:kolmogorov_high}.
\subsection{Proof of Lemma \ref{lem:kolmogorov_mixing}}
\label{sec:proof_lem_kolmogorov_mixing}
We have
\[ I\left(\sigma(Y_t),\sigma(Y_{t+s})\right) = \mbf{K}\left( \mcal{L}\left(Y_t,Y_{t+s}\right) || \mcal{L}\left(Y_t\right) \otimes \mcal{L}\left(Y_{t+s}\right) \right) = \mbb{E}\left[ \mbf{K}\left( \mcal{L}\left(Y_{t+s}|Y_t\right) || \mcal{L}\left(Y_{t+s}\right) \right) \right]. \]
Since $(Y_t)_{t\geq 0}$ is stationary we have $\mcal{L}(Y_{t+s})=\overline{P}$. For $x\in\mbb{R}^d$ fixed, we write
\[ A_x(s) = \mbf{K}\left( \mcal{L}(Y_s^x) ||\overline{P} \right), \]
where $Y_t^x$ is the solution of (\ref{eq:kolmogorov_equation}) satisfying $Y^x_0=x$. We follow the proof of Theorem 3.2.7 \cite{Royer} with their notation. From (44) therein we have
\begin{equation}
\label{eq:royer_0}
A_x(s) \leq \mbb{E}\left[ \left( \log(Z) + U(x) + U(W_s) - 2v(W_s) -
\frac{1}{2}\int_0^s [|\nabla U|^2 - \Delta U](W_t) dt \right) F \right],
\end{equation}
where
\begin{itemize}
\item $W$ is the Brownian motion starting from $x$,
\item $F$ is the density of the distribution of $X^x$ over $\mcal{C}([0, s])$ with respect to the distribution $P$ of $W$ given by
\[ F = \exp\left( U (x) - U (W_s ) - \frac{1}{2} \int_0^s [|\nabla U |^2 - \Delta U ] (W_t ) dt \right), \]
\item $v$ is such that $\exp(-2v)$ is the Gaussian density of $\mcal{L}(W_s)$ with respect to the Lebesgue measure, i.e.
\begin{equation}
\label{eq:exp_v}
\exp(-2v(y)) = (2\pi t)^{-d/2} \exp\left( -\frac{(x-y)^2}{2s} \right), \forall y\in\mbb{R}^d. 
\end{equation}
\end{itemize}
Let us check that the right-hand side of (\ref{eq:royer_0}) is finite. From (\ref{eq:exp_v}), we have $-2v(y)\leq -\frac{d}{2}\log (2\pi s)$. Also
\[ - \frac{1}{2}\int_0^s \left[ |\nabla U|^2 - \Delta U \right](W_t) dt \leq - \frac{C s}{2}, \]
where $C$ is given by (\ref{eq:hyp_kolmogorov}). Since $\mbb{E}F=1$, we get
\[ A_x(s) \leq  \log(Z) + U(x) - \frac{d}{2}\log(2\pi s) - \frac{C s}{2} + \mbb{E}\left[ U(W_s) F \right]. \]
We only need to consider the last term $\mbb{E}\left[ U(W_s) F \right]$. We have
\begin{align*}
\mbb{E}\left[ U(W_s) F \right] &= \mbb{E}\left[ U(W_s)  \exp\left( U(x) - U (W_s) - \frac{1}{2} \int_0^s [|\nabla U |^2 - \Delta U ] (W_t) dt \right)  \right]\\
&= e^{U(x)} \mbb{E}\left[ U(W_s) \exp\left( - U (W_s) - \frac{1}{2} \int_0^s [|\nabla U |^2 - \Delta U ] (W_t) dt \right) \right]\\
&\leq e^{U(x)-\frac{C s}{2}} \mbb{E}\left[ U(W_s) \exp\left( - U (W_s)\right)  \right]\\
&\leq e^{U(x)-\frac{C s}{2}} \mbb{E}\left[ U^+(W_s) \exp\left( - U (W_s)\right)  \right]\\
&\leq e^{U(x)-\frac{C s}{2}} ||g||_{\infty},
\end{align*}
where $g$ is defined on $\mbb{R}^+$ by $g(x)=x\exp(-x)$. We end up with
\begin{align}
A_x(s) &\leq  \log(Z) + U(x) - \frac{d}{2}\log(2\pi s) - \frac{C s}{2} + e^{U(x)-\frac{C s}{2}} ||g||_{\infty}\nonumber\\
&\leq \log(Z) - \frac{d}{2}\log(2\pi s) - \frac{C s}{2} + e^{U(x)} ||g||_{\infty} \left( 1 + e^{-\frac{C s}{2}} \right).\label{eq:royer_1}
\end{align}
Therefore, $A_x(s)$ is finite for all $s>0$ and all $x\in\mbb{R}^d$. From Theorem 3.1.29 and Theorem 3.2.5 of Royer \cite{Royer}, for all $s_0>0$, we have
\begin{equation}
\label{eq:royer_2}
A_x(s) \leq A_x(s_0)\exp\left(-2 m (s-s_0)\right), \forall s>s_0.
\end{equation}
Therefore with (\ref{eq:royer_1}) and (\ref{eq:royer_2}), we have
\begin{align*}
I\left( \sigma(Y_t), \sigma(Y_{t+s}) \right) &= \mbb{E}\left[ A_{Y_t}(s) \right]\\
&\leq \exp\left( - 2 m (s-s_0) \right) \mbb{E}\left[ A_{Y_t}(s_0) \right]\\
&\leq e^{- 2 m (s-s_0)} \left[ \log(Z) - \frac{d}{2}\log(2\pi s_0) - \frac{C s_0}{2} + \mbb{E}\left[ e^{U(Y_t)} \right] ||g||_{\infty} (1+e^{-\frac{C s_0}{2}}) \right]\\
&= e^{- 2 m (s-s_0)} \left[ \log(Z) - \frac{d}{2}\log(2\pi s_0) - \frac{C s_0}{2} + ||g||_{\infty} (1+e^{-\frac{C s_0}{2}}) Z^{-1} \int_{\mbb{R}^d} e^{-U(x)} dx \right]\\
&=: C(s_0) e^{-2 m s},
\end{align*}
for $s\geq s_0>0$ with $C(s_0)<\infty$ since $\int_{\mbb{R}^d} e^{-\alpha U(x)} dx<\infty$ for all $\alpha$.
\subsection{Proof of Lemma \ref{lem:entropy_log_concave}}
\label{sec:proof_lem_entropy_log_concave}
We divide the proof in two parts, first the case $d\leq 3$ and the case $d\geq 4$ in a second time.\par
\textbf{Case $d\in\{1,2,3\}$}. For $\xi> 0$ and $\nu\in(0,1)$, let
\[ \tilde{\mcal{F}}_d^{\xi,\nu}= \left\{ \tilde{f}\in\mcal{F}_d: || \overline{x}_{\tilde{f}} ||_2 \leq \xi \text{  and  } 1-\nu < \lambda_{\min}(\Sigma_{\tilde{f}}) \leq \lambda_{\max}(\Sigma_{\tilde{f}}) \leq 1 + \nu \right\}.\]
We first state the classic bound
\begin{equation}
\label{eq:entropy_moyenne}
N( B_2(M), ||\cdot||_2, \epsilon) \leq \left( \frac{3M}{\epsilon} \right)^d,
\end{equation}
where $B_2(M)$ is the ball of radius $M$ in $\mbb{R}^d$ with respect to the Euclidean distance $||\cdot||_2$. Let $B_2(M)\left[\sqrt{\lambda_-}\right]$ be a $\sqrt{\lambda_-}$-net of $B_2(M)$ with respect to the Euclidean distance $||\cdot||_2$, with $\left|B_2(M)\left[\sqrt{\lambda_-}\right]\right| \leq (3M/\lambda_-)^d$. Let $\text{Sym}(\lambda_-,\lambda_+)[\eta_d\lambda_-]$ be a $\eta_d\lambda_-$-net of $\text{Sym}(\lambda_-,\lambda_+)$ with respect to the operator norm $||\cdot||_{op}$, with $|\text{Sym}(\lambda_-,\lambda_+)[\eta_d\lambda_-]|\leq N_{\Sigma}(\lambda_+,\lambda_-,d,\eta_d\lambda_-)$. Let $\tilde{F}_d^{1,\eta_d}[\epsilon]$ be an $\epsilon$-net of $\tilde{F}_d^{1,\eta_d}$ with respect to the Hellinger distance. We define
\[ \mcal{F}_{\lambda_-,\lambda_+,M}[\epsilon] := \left\{ (\det \Sigma)^{-1/2} g\left( \Sigma^{-1/2}\left(\cdot - \overline{x} \right) \right); \begin{array}{l}
\overline{x} \in B_2(M)\left[\sqrt{\lambda_-}\right], \\
\Sigma \in \text{Sym}(\lambda_-,\lambda_+)[\eta_d\lambda_-],\\
g\in \tilde{F}_d^{1,\eta_d}[\epsilon]
\end{array}    \right\}  \]
and we show it is an $\epsilon$-net of $\mcal{F}_{\lambda_-,\lambda_+,M}$ with respect to the Hellinger distance. For $f\in\mcal{F}_{\lambda_-,\lambda_+,M}$, there is $\Sigma$ in $\text{Sym}(\lambda_-,\lambda_+)[\eta_d \lambda_-]$ and $\overline{x}$ in $B_2(M)[\sqrt{\lambda_-}]$ such that
\[ || \overline{x}_f - \overline{x} ||_2 \leq \sqrt{\lambda_-} \text{  and  } || \Sigma_f - \Sigma ||_{op}\leq \lambda_- \eta_d. \]
We write $\tilde{f}=( \det \Sigma)^{1/2} f\left( \Sigma^{1/2} \cdot + \overline{x} \right)$. Let us check that $\tilde{f}$ belongs to $\tilde{F}_d^{1,\eta_d}$. We have
\[ || \overline{x}_{\tilde{f}} ||_2 = || \Sigma^{-1/2} (\overline{x}_f - \overline{x}) ||_2 \leq \frac{ || \overline{x}_f - \overline{x} ||_2 }{ \sqrt{\lambda_-} } \leq  1,\]
and
\[ || \Sigma_{\tilde{f}} - I ||_{op} = || \Sigma^{-1/2} \Sigma_f  \Sigma^{-1/2} - I ||_{op} = || \Sigma^{-1/2} ( \Sigma_f - \Sigma ) \Sigma^{-1/2} ||_{op} \leq \frac{ || \Sigma_f - \Sigma ||_{op} }{\lambda_-}\leq \eta_d. \]
Therefore $\tilde{f} \in \tilde{F}_d^{1,\eta_d}$ and there is $g \in \tilde{F}_d^{1,\eta_d}[\epsilon]$ such that $h\left( \tilde{f},g\right) \leq \epsilon$. Since the Hellinger distance is invariant by translation and scaling, we have
\[ h\left( f, (\det \Sigma)^{-1/2} g\left( \Sigma^{-1/2}(\cdot - \mu) \right) \right) = h\left( \tilde{f}, g \right) \leq \epsilon, \]
which proves that $\mcal{F}_{\lambda_-,\lambda_+,M}[\epsilon]$ is an $\epsilon$-net of $\mcal{F}_{\lambda_-,\lambda_+,M}$. Therefore
\[ |\mcal{F}_{\lambda_-,\lambda_+,M}[\epsilon]| \leq \left(\frac{3M}{\sqrt{\lambda_-}}\right)^d \times N_{\Sigma}(\lambda_+,\lambda_-,d,\eta_d\lambda_-) \times |\tilde{F}_d^{1,\eta_d}[\epsilon]|. \]
We need to bound the different entropy numbers now. For a metric space $(\mscr{A},d)$ and $\epsilon>0$, we denote by $N(\epsilon,\mscr{A},d)$ the minimal number of balls of radius $\epsilon$, with respect to $d$, to cover $\mscr{A}$.\par 
The next result provides a bound on the entropy for the class of covariance matrices we are considering. Let $||\cdot||_{op}$ denote the operator norm on square matrices induced by the Euclidean distance. For matrices with real-valued eigenvalues, it is equivalent to the largest absolute value of its eigenvalues.
\begin{lemme}
\label{lem:entropy_symmetric}
We have
\begin{equation}
\label{eq:entropy_symmetric}
N\left( \epsilon, \text{Sym}(\lambda_-,\lambda_+), ||\cdot||_{op}\right) \leq
\begin{cases}
\frac{3(\lambda_+-\lambda_-)}{\epsilon} \text{ for }d=1,\\
\left(\frac{9}{\epsilon}\right)^3 (\lambda_+-\lambda_-)^2 \lambda_+\pi \text{ for } d=2,\\
2 \left( \frac{2 \cdot 3^{5/4} \sqrt{\lambda_+(\lambda_+-\lambda_-)\pi}}{\epsilon} \right)^6 \text{ for } d=3.                                                           \end{cases}
\end{equation}
In higher dimensions, we have
\begin{align*}
N\left( \epsilon, \text{Sym}(\lambda_-,\lambda_+), ||\cdot||_{op} \right) &\leq C \left( \frac{3}{4} \right)^d \frac{\pi^{d(d-1)/2}}{e^{(d-1)(d-2)/4}} (2\lambda_+)^{d(d-1)/2} (\lambda_+-\lambda_-)^d\\
&\times (d+1)^{d(d+1)/2} d^{(d-1)(d+2)/2} (d-1)^{(d-1)/2} \epsilon^{-d(d+1)/2},
\end{align*}
with $C=\frac{e^{1/2}}{3^{1/2}2^3}$.
\end{lemme}
Theorem 4 \cite{samworth_kim} gives a bound on $|\tilde{F}_d^{1,\eta_d}[\epsilon]|$ which allows to conclude the proof of Theorem \ref{th:kolmogorov}.\par 
\textbf{Case $d\geq 4$}. We use Theorem 3 of Kur \etal \cite{kur_concave}. We follow some of their notation. Let $d\geq 4$. There exist positive constants $\xi_d$ and $\overline{K}_d$ such that
\[ \log N(\epsilon,\mscr{F}_{d,\tilde{I}},h) \leq \overline{K}_d \epsilon^{-(d-1)} \log_{++}(\epsilon^{-1})^{(d+1)(d+2)/2},\]
where $\mscr{F}_{d,\tilde{I}}$ is the set of distributions associated to
\[ \mcal{F}_{d,\tilde{I}} = \left\{ \tilde{f}\in\mcal{F}_d: || \overline{x}_{\tilde{f}} ||_2 \leq \xi_d \text{  and  } 1/2 < \lambda_{\min}(\Sigma_{\tilde{f}}) \leq \lambda_{\max}(\Sigma_{\tilde{f}}) \leq 2 \right\}.\]
Let $\mcal{F}_{d,\tilde{I}}[\epsilon]$ be a set of probability densities with respect to the Lebesgue measure such that $\mscr{F}_{d,\tilde{I}}[\epsilon]=\{ f(x)dx;f\in\mcal{F}_{d,\tilde{I}}\}$ is an $\epsilon$-net of $\mscr{F}_{d,\tilde{I}}$ with respect to the Hellinger distance and
\[ \log |\mscr{F}_{d,\tilde{I}}[\epsilon]| \leq \overline{K}_d \epsilon^{-(d-1)} \log_{++}(\epsilon^{-1})^{(d+1)(d+2)/2}. \]
Let $B_2(M)\left[\xi_d\sqrt{\lambda_-}\right]$ be a $\xi_d\sqrt{\lambda_-}$-net of $B_2(M)$ with respect to the Euclidean distance $||\cdot||_2$, with $\left|B_2(M)\left[\xi_d\sqrt{\lambda_-}\right]\right| \leq (3M/\xi_d\sqrt{\lambda_-})^d$. Let $\text{Sym}(\lambda_-,\lambda_+)[\lambda_-/3]$ be a $\lambda_-/3$-net of $\text{Sym}(\lambda_-,\lambda_+)$ with respect to the operator norm $||\cdot||_{op}$, with $|\text{Sym}(\lambda_-,\lambda_+)[\lambda_-/3]| \leq N_{\Sigma}(\lambda_+,\lambda_-)$. We define
\[ \mcal{F}_{\lambda_-,\lambda_+,M}[\epsilon] := \left\{ (\det \Sigma)^{-1/2} g\left( \Sigma^{-1/2}\left(\cdot - \overline{x} \right) \right); \begin{array}{l}
\overline{x} \in B_2(M)\left[\xi_d \sqrt{\lambda_-}\right], \\
\Sigma \in \text{Sym}(\lambda_-,\lambda_+)[\lambda_-/3],\\
g\in \mcal{F}_{d,\tilde{I}}[\epsilon]
\end{array}    \right\}  \]
and we show that $\mscr{F}_{\lambda_-,\lambda_+,M}[\epsilon]=\{ f(x) dx; f\in \mcal{F}_{\lambda_-,\lambda_+,M}[\epsilon]\}$ is an $\epsilon$-net of $\mscr{F}_{\lambda_-,\lambda_+,M}$ with respect to the Hellinger distance. For $f\in\mcal{F}_{\lambda_-,\lambda_+,M}$, there is $\Sigma$ in $\text{Sym}(\lambda_-,\lambda_+)[\lambda_-/3]$ and $\overline{x}$ in $B_2(M)[\xi_d \sqrt{\lambda_-}]$ such that
\[ || \overline{x}_f - \overline{x} ||_2 \leq \xi_d \sqrt{\lambda_-} \text{  and  } || \Sigma_f - \Sigma ||_{op}\leq \lambda_-/3. \]
We write $\tilde{f}=(\det \Sigma)^{1/2} f\left( \Sigma^{1/2} \cdot + \overline{x} \right)$. Let us check that $\tilde{f}$ belongs to $\mcal{F}_{d,\tilde{I}}$. We have
\[ || \overline{x}_{\tilde{f}} ||_2 = || \Sigma^{-1/2} (\overline{x}_f - \overline{x}) ||_2 \leq \frac{ || \overline{x}_f - \overline{x} ||_2 }{ \sqrt{\lambda_-} } \leq  \xi_d,\]
and
\[ || \Sigma_{\tilde{f}} - I || = || \Sigma^{-1/2} \Sigma_f  \Sigma^{-1/2} - I || = || \Sigma^{-1/2} (\Sigma_f-\Sigma) \Sigma^{-1/2} || \leq \frac{|| \Sigma_f-\Sigma ||}{\lambda_-} \leq 1/3. \]
Hence
\[ \lambda_{\min}(\Sigma_{\tilde{f}})\geq 2/3 > 1/2 \text{  and  } \lambda_{\max}(\Sigma_{\tilde{f}}) \leq 4/3 < 2. \]
Therefore we have $\tilde{f} \in \mcal{F}_{d,\tilde{I}}$ and there is $g \in \mcal{F}_{d,\tilde{I}}[\epsilon]$ such that $h\left( \tilde{f}(x) dx, gx) dx\right) \leq \epsilon$. Since the Hellinger distance is invariant by translation and scaling, we have
\[ h\left( f(x) dx, (\det \Sigma)^{-1/2} g\left( \Sigma^{-1/2}( x - \overline{x}) \right) dx \right) = h\left( \tilde{f}(x) dx, g(x) dx \right) \leq \epsilon, \]
which proves that $\mscr{F}_{\lambda_-,\lambda_+,M}[\epsilon]$ is an $\epsilon$-net of $\mscr{F}_{\lambda_-,\lambda_+,M}$. Therefore
\[ |\mscr{F}_{\lambda_-,\lambda_+,M}[\epsilon]| \leq \left( \frac{3M}{\xi_d \sqrt{\lambda_-}} \right)^d \times N_{\Sigma}(\lambda_+,\lambda_-,d) \times |\mscr{F}_{d,\tilde{I}}[\epsilon]|. \]
With Lemma \ref{lem:entropy_symmetric} we get
\begin{align*}
|\mscr{F}_{\lambda_-,\lambda_+,M}[\epsilon]| &\leq C \left( \frac{3M}{\xi_d \sqrt{\lambda_-}} \right)^d \left(\frac{3}{4}\right)^d \frac{\pi^{d(d-1)/2}}{e^{(d-1)(d-2)/4}} (2\lambda_+)^{d(d-1)/2} \left(\lambda_+-\lambda_-\right)^d\\
&\times (d+1)^{d(d+1)/2} d^{(d-1)(d+2)/2} (d-1)^{(d-1)/2} \left(\frac{\lambda_-}{3}\right)^{-d(d+1)/2}\\
&\times \exp\left( \overline{K}_d \epsilon^{-(d-1)} \log(\epsilon^{-1})^{(d+1)(d+2)/2} \right)\\
&\leq C_d \frac{\lambda_+^{d(d-1)/2} M^d (\lambda_+-\lambda_-)^d}{ \lambda_-^{d(d+1)/2} } \exp\left( \overline{K}_d \epsilon^{-(d-1)} \log(\epsilon^{-1})^{(d+1)(d+2)/2} \right).
\end{align*}
\subsubsection{Proof of Lemma \ref{lem:entropy_symmetric}}
For $d=1$, we have $\text{Sym}(\lambda_-,\lambda_+)=[\lambda_-,\lambda_+]$. The result follows from classical entropy bounds. Otherwise, every real valued symmetric matrix $\Sigma$ can be written as $\Sigma=UDU^T$ where $D$ is the diagonal matrix containing the real eigenvalues of $\Sigma$ and $U$ is an orthonormal matrix. For $\Sigma_1=U_1\text{diag}(\lambda_{1,1},\dots,\lambda_{d,1}) U_1^T$ and  $\Sigma_2=U_2\text{diag}(\lambda_{1,2},\dots,\lambda_{d,2}) U_2^T$ we have
\begin{align*}
||\Sigma_1-\Sigma_2|| &\leq || U_1 (D_1-D_2) U_1^T|| + || (U_1-U_2) D_2 U_1^T|| + || U_2 D_2 (U_1-U_2)^T||\\
&\leq ||D_1-D_2||+2\lambda_+||U_1-U_2||\\
&= \max_{1\leq i\leq d} |\lambda_{i,1}-\lambda_{i,2}|+2\lambda_+||U_1-U_2||.
\end{align*}
Therefore
\[ N\left(\text{Sym}(\lambda_-,\lambda_+),||\cdot||,\epsilon\right)\leq N\left(B((\lambda_+-\lambda_-)/2),||\cdot||_{\infty},\epsilon_1\right) \times N\left(\text{ON}(d), ||\cdot||,\epsilon_2\right) \]
with $\epsilon=\epsilon_1+2\lambda_+\epsilon_2$. We have the classic bound
\[ N\left(B((\lambda_+-\lambda_-)/2),||\cdot||_{\infty},\epsilon_1\right)\leq \left(3\frac{\lambda_+-\lambda_-}{2\epsilon_1}\right)^d. \]
\begin{itemize}
\item For $d=2$, the orthonormal matrices are of the form 
\[ U_{\alpha,\theta}=\begin{pmatrix}
\cos(\theta) & -\alpha \sin(\theta)\\
\sin(\theta) & \alpha \cos(\theta)
\end{pmatrix},\theta\in[0,2\pi],\alpha\in\{-1,1\}.
\]
We have
\[ ||U_{\alpha,\theta}-U_{\alpha,\theta'}||^2 = 2 [1-\cos(\theta-\theta')] \leq (\theta-\theta')^2, \]
and therefore
\[ N\left(\text{ON}(2), ||\cdot||,\epsilon\right)\leq 2 \frac{3\pi}{\epsilon}=6\pi/\epsilon,
\]
where the factor 2 comes from the presence of $\epsilon$ for positively and negatively oriented basis. We obtain the final result for $\epsilon_1=2\epsilon/3$ and $\epsilon_2=\epsilon/6\lambda_+$.
\item We proceed similarly for $d=3$. Every orthonormal basis in dimension 3 can be written in the form
\[ U_{\epsilon,\theta,\beta,\gamma}:=\begin{pmatrix} \cos\theta & \cos\gamma \sin\theta & -\epsilon \sin \gamma \sin\theta\\
\sin \theta \cos\beta & -\cos\gamma \cos\theta \cos\beta +\sin\gamma \sin\beta & \epsilon(\sin\gamma \cos\theta \cos\beta + \cos \gamma \sin\beta)\\
\sin\theta\sin\beta & -\cos\gamma \cos\theta \sin\beta - \sin\gamma \cos\beta & \epsilon(\sin\gamma \cos\theta \sin\beta - \cos\gamma \cos\beta)
\end{pmatrix},\]
$\theta\in[0,2\pi],\beta\in[0,2\pi],\gamma\in[0,2\pi],\epsilon\in\{-1,1\}$. As before, one can check that we have
\begin{align*}
|| U_{\epsilon,\theta,\beta,\gamma} - U_{\epsilon,\theta',\beta,\gamma} || &\leq |\theta-\theta'|^2\\
|| U_{\epsilon,\theta,\beta,\gamma} - U_{\epsilon,\theta,\beta',\gamma} || &\leq |\beta-\beta'|^2\\
|| U_{\epsilon,\theta,\beta,\gamma} - U_{\epsilon,\theta,\beta,\gamma'} || &\leq |\theta-\theta'|^2.
\end{align*}
Therefore we have
\begin{equation}
\label{eq:on_3_induction}
N\left(\text{ON}(3), ||\cdot||,\epsilon\right)\leq \left(N\left([0,2\pi],|\cdot|,\epsilon/\sqrt{3}\right) \right)^3\leq 2\left( \frac{3\sqrt{3}\pi}{\epsilon} \right)^3,
\end{equation}
where the factor 2 comes from the presence of $\epsilon$ for positively and negatively oriented basis. We obtain the final result for $\epsilon_1=\epsilon/2$ and $\epsilon_2=\epsilon/4\lambda_+$.
\item For higher dimensions, we have the following lemma.
\begin{lemme}
\label{lem:on_entropy_induction}
For $d\geq 3$, we can build an $\epsilon$-net $ON(d)[\epsilon]$ of $ON(d)$ with respect to the operator norm such that
\[ |ON(d)[\epsilon]|\leq C \frac{\pi^{d(d-1)/2}}{e^{(d-1)(d-2)/4}} d^{(d-1)(d+2)/2} (d-1)^{(d-1)(d+1)/2} \epsilon^{-d(d-1)/2}, \forall d\geq 1, \]
with $C=\frac{e^{1/2}}{3^{1/2}2^3}$.
\end{lemme}
We obtain the final bound with $\epsilon_1 = \frac{2\epsilon}{d+1}$ and $\epsilon_2 = \frac{\epsilon}{2\lambda_+} \frac{d-1}{d+1}$.
\end{itemize}
\subsubsection{Proof of Lemma \ref{lem:on_entropy_induction}}
We prove this by induction. From (\ref{eq:on_3_induction}) we have the desired inequality for $d=3$ with $C_3=\frac{e^{1/2}}{3^{1/2} 2^3}$. Let $\epsilon$ be in $(0,1]$ and $d\geq 3$. Let us now assume that for $\lambda_1>0$ we have a $\lambda_1$-net $ON(d)[\lambda_1]$ with
\[ |ON(d)[\lambda_1]| \leq C \frac{\pi^{d(d-1)/2}}{e^{(d-1)(d-2)/4}} d^{(d-1)(d+2)/2} (d-1)^{(d-1)(d+1)/2} \lambda_1^{-d(d-1)/2} .\]
Let $U\in\mbb{R}^{d+1}$ be a unitary vector, i.e.  $U_1^2+\dots+U_{d+1}^2=1$. There is $\theta \in [0,2\pi]^d$ such that $U=f(\theta)$ with
\[ U_i = f_i(\theta) := \cos \theta_i \prod_{j\leq i} \sin \theta_j, \]
with the convention $\theta_{d+1}=0$ and that a product over an empty set of indices is equal to 1. We define applications $a_1,\dots,a_d,a_{d+1}$ by $a_1=id$ and
\[ a_i(\theta)=\left(\theta_1+\frac{\pi}{2},\dots,\theta_{i-1}+\frac{\pi}{2},\theta_i,\dots,\theta_d\right),\forall i\in\{2,\dots,d+1\}. \]
One can check that the set of vectors $A_1(\theta),\dots,A_{d+1}(\theta)\in\mbb{R}^{d+1}$, given by $A_i(\theta)=f(a_i(\theta))$ for $i$ in $\{1,2,\dots,d+1\}$, is an orthonormal basis of $\mbb{R}^d$. We take $n_j= \left\lceil \frac{\sqrt{d+1-j}}{\lambda_2} \right\rceil,\forall j\in\{1,2,\dots,d\}$ and we take
\[ \mscr{A}_{d+1}[\lambda_2] := \left\{ A(\psi_{i_1,\dots,i_d}); i_j\in\{1,2,\dots,n_j\}, j\in\{1,2,\dots,d\} \right\} \subset ON(d+1), \]
with 
\[ \psi_{i_1,\dots,i_d} = \left( \frac{\pi (2i_j-1)}{n_j} \right)_{1\leq j\leq d}. \]
\begin{lemme}
\label{lem:on_net}
The set
\[ O[\lambda_1,\lambda_2] := \left\{ A \begin{pmatrix}
    1& 0\\
    0& B                                                                                                                                                                                                     \end{pmatrix}; A\in \mscr{A}_{d+1}[\lambda_2], B \in ON(d)[\lambda_1] \right\}, \]
is a $\lambda_1 + \sqrt{d}\pi \lambda_2$-net of $ON(d+1)$ with respect to the operator norm.
\end{lemme}
One can easily check that we have the following bound
\[ |\mscr{A}_{d+1}[\lambda_2]| \leq \left(\frac{2}{\lambda_2}\right)^d \sqrt{ d!}. \]
Therefore, we have
\begin{align*}
|O[\lambda_1,\lambda_2]| &= |ON(d)[\lambda_1]| \times |\mscr{A}_{d+1}[\lambda_2]|\\
&\leq C \frac{\pi^{d(d-1)/2}}{e^{(d-1)(d-2)/4}} d^{(d-1)(d+2)/2} (d-1)^{(d-1)(d+1)/2} \lambda_1^{-d(d-1)/2} \times \left(\frac{2}{\lambda_2}\right)^d \sqrt{ d!}.
\end{align*}
For $\lambda_1=\epsilon\frac{d-1}{d+1}$ and $\lambda_2=\epsilon \frac{2}{\sqrt{d}\pi (d+1)}$, we get
\begin{align*}
&|O[\lambda_1,\lambda_2]|\\
&\leq C \frac{\pi^{d(d-1)/2}}{e^{(d-1)(d-2)/4}} d^{(d-1)(d+2)/2} (d-1)^{(d-1)(d+1)/2} \left( \frac{d+1}{d-1} \right)^{d(d+1)/2} \epsilon^{-d(d-1)/2}\\
&\times \sqrt{ d!}  \left( \sqrt{d}\pi (d+1) \right)^d \epsilon^{-d}\\
&= C (d-1)^{-(d+1)/2} d^{-1} \sqrt{ d!} e^{(d-1)/2} \frac{\pi^{d(d+1)/2}}{e^{(d-1)(d-2)/2}} (d+1)^{d(d+3)/2} d^{d(d+2)/2} \epsilon^{-d(d+1)/2}.
\end{align*}
We use the bound $n!\leq \sqrt{2\pi} n^{n+\frac{1}{2}} e^{-n} e^{\frac{1}{12n}}$ and we get
\begin{align*}
&|O[\lambda_1,\lambda_2]|\\
&\leq C (d-1)^{-(d+1)/2} d^{-1/2} \sqrt{(d-1)!} e^{(d-1)/2} \frac{\pi^{d(d+1)/2}}{e^{(d-1)(d-2)/2}} (d+1)^{d(d+3)/2} d^{d(d+2)/2} \epsilon^{-d(d+1)/2}\\
&\leq C (d-1)^{-3/4} d^{-1/2} (2\pi)^{1/4} e^{\frac{1}{24(d-1)}} \frac{\pi^{d(d+1)/2}}{e^{(d-1)(d-2)/2}} (d+1)^{d(d+3)/2} d^{d(d+2)/2} \epsilon^{-d(d+1)/2}.
\end{align*}
We have
\[ (d-1)^{-3/4} d^{-1/2} (2\pi)^{1/4} e^{\frac{1}{24(d-1)}} \leq 1 \]
for all $d\geq 3$. Therefore, we satisfy the desired property for $d+1$ with $ON[\epsilon]=O[\lambda_1,\lambda_2]$.
\subsubsection{Proof of Lemma \ref{lem:on_net}}
Let $C=(C_1\dots C_{d+1})$ be in $ON(d+1)$. There is $\theta$ in $[0,2\pi]^d$ such that $C_1=A_1(\theta)$. Let $B$ be the matrix in $ON(d)$ given by
\[ A(\theta)^T C = \begin{pmatrix}
            1& 0\\
            0& B
           \end{pmatrix}
.\]
For $\theta\in[0,2\pi]^d$ there exists $\psi_{i_1,\dots,i_d}$ such that
\[ \left| \theta_i - \frac{\pi (2i_j-1)}{n_j} \right| \leq \frac{\pi}{n_j} \leq \frac{\pi \lambda_2}{\sqrt{d+1-j}},\forall j\in\{1,\dots,d\}. \]
\begin{lemme}
\label{lem:approx_basis_theta}
We have
\[ ||A(\theta)-A(\theta+h)||_{op} \leq \sqrt{ \suml_{k=0}^{d-1} (d-k) h_{k+1}^2 }. \]
\end{lemme}
Therefore we have
\[ || A(\theta) - A(\psi_{i_1,\dots,i_d}) ||_{op} \leq d^{1/2} \pi \epsilon. \]
There exists $B'$ in $ON(d)[\lambda_1]$ such that $||B-B'||_{op} \leq \lambda_1$. We define $C'\in ON(d+1)$ by 
\[ C'=A(\psi_{i_1,\dots,i_d}) \begin{pmatrix}
    1& 0\\
    0& B'                                                                                                                                                                                                     \end{pmatrix}\in ON[\lambda_1,\lambda_2]. \]
Then we have
\begin{align*}
|| C-C' ||_{op} &\leq \left|\left| A(\theta) \begin{pmatrix}
            0& 0\\
            0& B-B'
           \end{pmatrix} \right|\right|_{op} + \left|\left| \left( A(\theta) - A(\psi_{i_1,\dots,i_d}) \right) \begin{pmatrix}
            1& 0\\
            0& B'
           \end{pmatrix} \right|\right|_{op}\\
           &\leq ||B-B'||_{op} + \left|\left| A(\theta) - A(\psi_{i_1,\dots,i_d}) \right|\right|_{op}\\
           &\leq \lambda_1 + d^{1/2} \pi \lambda_2.
\end{align*}
\subsubsection{Proof of Lemma \ref{lem:approx_basis_theta}}
For $\theta\in\mbb{R}^d$ and $h\in\mbb{R}^d$, we define $U_0=f(\theta)$ and
\[ U_i=f(\theta_1+h_1,\dots,\theta_i+h_i,\theta_{i+1},\dots,\theta_d),i\in\{1,\dots,d\}. \]
Similarly, we write $A^{(i)}=A(\theta^{(h,i)})$ with
\[ \theta^{(h,i)} = (\theta_1+h_1,\dots,\theta_i+h_i,\theta_{i+1},\dots,\theta_d), \]
for $i\in\{0,1,\dots,d\}$ and $j\in\{1,\dots,d+1\}$. It implies $A^{(0)}_1=U_0$ and $A^{(d)}_1=U_d$. We have
\begin{align*}
A^{(k)}_{ij} &= f_i(a_j(\theta^{(h,k)})) = \cos\left( a_j(\theta^{(h,k)}) \right) \prod_{l\leq i} \sin\left( a_j(\theta^{(h,k)}) \right)\\
&= \cos\left( \theta_i +\mathbbm{1}_{i<j}\frac{\pi}{2} + \mathbbm{1}_{l\leq i} h_i \right) \prod_{l\leq i} \sin\left( \theta_l + \mathbbm{1}_{l< j} \frac{\pi}{2} + \mathbbm{1}_{l\leq k} h_l \right),
\end{align*}
and therefore
\begin{align*}
A^{(k+1)}_{ij} - A^{(k)}_{ij} &=
\begin{cases}
&0 \text{  if  } i\leq k\\
&\prod\limits_{l\leq k} \sin\left( \theta_l +\mathbbm{1}_{l<j}\frac{\pi}{2} + h_l \right)\\
&\times \left[ \cos \left( \theta_{k+1} +\mathbbm{1}_{k+1<j}\frac{\pi}{2} + h_{k+1} \right) - \cos \left( \theta_{k+1} +\mathbbm{1}_{k+1<j}\frac{\pi}{2} \right) \right] \text{  if  } i=k+1\\
&\prod\limits_{\substack{l<i\\
l\neq k+1}} \sin\left( \theta_l +\mathbbm{1}_{l<j}\frac{\pi}{2} + \mathbbm{1}_{l\leq k} \right) \times \cos \left( \theta_i +\mathbbm{1}_{i<j}\frac{\pi}{2} \right)\\
&\times \left[ \sin \left( \theta_{k+1} +\mathbbm{1}_{k+1<j}\frac{\pi}{2} + h_{k+1} \right) - \sin \left( \theta_{k+1} +\mathbbm{1}_{k+1<j}\frac{\pi}{2} \right) \right] \text{  if  } i>k+1,
\end{cases}\\
&= 2 \sin\left(\frac{h_{k+1}}{2}\right) \prod\limits_{l\leq k} \sin\left( \theta_l +\mathbbm{1}_{l<j}\frac{\pi}{2} + h_l \right)\\
&\times
\begin{cases}
&0 \text{  if  } i\leq k\\
&- \sin\left(\theta_{k+1} + \mathbbm{1}_{k+1<j}\frac{\pi}{2} + \frac{h_{k+1}}{2}\right)  \text{  if  } i=k+1\\
&\cos \left(\theta_{k+1} + \mathbbm{1}_{k+1<j}\frac{\pi}{2} + \frac{h_{k+1}}{2}\right) \prod\limits_{k+1<l<i} \sin\left( \theta_l +\mathbbm{1}_{l<j}\frac{\pi}{2} \right)\\
&\times \cos \left( \theta_i +\mathbbm{1}_{i<j}\frac{\pi}{2} \right) \text{  if  } i>k+1.
\end{cases}
\end{align*}
We have ($k+1\leq d$, $k\geq 0$)
\begin{align*}
&||A^{(k+1)} - A^{(k)}||^2_F = \suml_{i,j} \left(A^{(k+1)}_{ij}-A^{(k)}_{ij}\right)^2\\
&= 4 \sin^2\left(\frac{h_{k+1}}{2}\right) \suml_{1\leq j\leq d+1}  \prod_{l\leq k} \sin^2\left( \theta_l + \mathbbm{1}_{l<j} \frac{\pi}{2} + h_l \right) \bigg[ \sin^2\left( \theta_{k+1} + \mathbbm{1}_{k+1<j} \frac{\pi}{2} + \frac{h_{k+1}}{2} \right)\\
&+ \cos^2\left( \theta_{k+1} + \mathbbm{1}_{k+1<j}\frac{\pi}{2} + \frac{h_{k+1}}{2} \right) \suml_{i=k+2}^{d+1} \prod_{k+1<l<i} \sin^2\left(\theta_l+\mathbbm{1}_{l<j} \frac{\pi}{2}\right) \cos^2\left(\theta_i+\mathbbm{1}_{i<j} \frac{\pi}{2} \right) \bigg]\\
&= 4 \sin^2\left(\frac{h_{k+1}}{2}\right) \suml_{1\leq j\leq d+1}  \prod_{l\leq k} \sin^2\left( \theta_l + \mathbbm{1}_{l<j} \frac{\pi}{2} + h_l \right)\\
&= 4 \sin^2\left(\frac{h_{k+1}}{2}\right) \bigg[ (d+1-k) \prod_{l\leq k} \sin^2\left( \theta_l + \mathbbm{1}_{l<j} \frac{\pi}{2} + h_l \right)\\
&+ \suml_{1\leq j\leq k}  \prod_{l\leq k} \sin^2\left( \theta_l + \mathbbm{1}_{l<j} \frac{\pi}{2} + h_l \right)\bigg]\\
&\leq 4 \sin^2\left(\frac{h_{k+1}}{2}\right) \bigg[ (d+1-k) \prod_{l\leq k} \cos^2\left( \theta_l + h_l \right) + 1 - \prod_{l\leq k} \cos^2\left( \theta_l + h_l \right) \bigg]\\
&\leq 4 \sin^2\left(\frac{h_{k+1}}{2}\right) (d-k) \prod_{l\leq k} \cos^2\left( \theta_l + h_l \right)\\
&\leq (d-k) h_{k+1}^2.
\end{align*}
Finally, with $||\cdot||_{op}\leq ||\cdot||_F$ we get
\begin{align*}
|| A^{(d)}-A^{(0)}||_{op} &\leq \suml_{k=0}^{d-1} || (A^{(k+1)}-A^{(k}))^T||_{op}\\
&\leq \suml_{k=0}^{d-1} (d-k) h^2_{k+1}.
\end{align*}
\section{Hidden Markov models}
This section gathers the proof of Theorems \ref{th:first_hmm}, \ref{th:hmm_log_concave}, \ref{th:hmm_log_concave_high}, \ref{th:ibragimov_exponential}, \ref{th:faster_rates}, Corollary \ref{coro:hmm_robust} and Proposition \ref{prop:exponential_family}, \ref{prop:regular_parametric_model}, \ref{prop:faster_rates}.
\subsection{Proof of Theorem \ref{th:first_hmm}}
\label{sec:proof_th_hmm}
The next result is proven in Section \ref{sec:proof_prop_rho_dimension_hmm} and gives a bound on the $\rho$-dimension function.
\begin{prop}
\label{prop:rho_dimension_hmm}
Under Assumption \ref{hyp:hmm_vc} and with $\delta(s)$ given by (\ref{eq:delta_s}, we can take
\[ D_{n(s,1)}\left( \mscr{M}_{\delta(s)} \right) = C L \overline{V} \left[ 1 + \log\left( \frac{K n(s,1)}{\overline{V}\wedge n(s,1)} \right) \right], \]
with $C=3 930$.
\end{prop}
With Theorem \ref{th:main_rho} we have
\begin{align*}
C \mbb{E} \left[ h^2\left( \overline{P}, \hat{P}_s \right) \right] &\leq  h^2\left( \overline{P}, \mscr{M}_{\delta} \right) + n^{-1} \suml_{i=1}^n h^2\left( P_i, \overline{P} \right)\\
&+ n^{-1} \suml_{b=1}^{s+1} \mbf{K}\left( \mbf{P}_{s,b}^* || \mbf{P}^{ind}_{s,b} \right) + (s+1) L \overline{V} \frac{ \log n}{n},
\end{align*}
for some positive constant $C$. The following result is proven in Proposition \ref{sec:proof_prop_approximation_hmm} and tells us how well $\mscr{M}_{\delta}$ approximates $\mscr{M}$.
\begin{prop}
\label{prop:approximation_hmm}
For $K\geq 2$, $w,v$ in $\mcal{W}_K$, $Q,R$ in $\mcal{T}_K$ and probability distributions $F_1,\dots,F_K,G_1,\dots,G_K$ on $\left(\mscr{Y},\mcal{Y}\right)$, we have 
\begin{align*}
h^2\left( P_{w,Q,F},P_{v,R,G}\right) &\leq h^2(w,v) + (L-1) \max_{k \in [K]} h^2\left( Q_{k\cdot}, R_{k\cdot} \right)\\
&+ L \max_{k\in[K]} h^2\left(F_k,G_k\right).
\end{align*}
\end{prop}
With Proposition \ref{prop:approximation_hmm} and inequality (B.5) in Lecestre \cite{robust_mixture} we have
\begin{equation}
\label{eq:approximation_delta}
h^2\left(P,\mscr{M}_{\delta}\right) \leq (K-1)L\delta + L \epsilon^2, \forall P\in \mscr{M}.
\end{equation}
With the choice of $\delta$ given in (\ref{eq:delta_s}) we get
\begin{align*}
C \mbb{E} \left[ h^2\left( \overline{P}, \hat{P}_s \right) \right] &\leq  h^2\left( \overline{P}, \mscr{M} \right) + n^{-1} \suml_{i=1}^n h^2\left( P_i, \overline{P} \right) + n^{-1} \suml_{b=1}^{s+1} \mbf{K}\left( \mbf{P}_{s,b}^* || \mbf{P}^{ind}_{s,b} \right)\\
&+ L \epsilon^2 + (s+1) L \overline{V} \frac{ \log n}{n},
\end{align*}
for some positive constant $C$. We now turn to the second bound in Theorem \ref{th:first_hmm}. The next result is proven later in Section \ref{sec:proof_lem_k_s_hmm}.
\begin{lemme}
\label{lem:k_s_hmm}
Under Assumption \ref{hyp:hmm_ergodic}, there are positive constants $C(Q^*)$ and $r(Q^*)$ that only depend on $Q^*$ such that
\[ n^{-1} \suml_{b=1}^{s+1} \mbf{K}\left(\mbf{P}^*_{s,b}||\mbf{P}^{ind}_{s,b}\right) \leq C(Q^*) e^{-r(Q^*) s},\forall s\geq L-1,\forall b\in[s+1], \]
and $h^2\left(P^*,P_i\right) \leq C(Q^*) e^{-r(Q^* ) i}$ for all $i\in[n]$.
\end{lemme}
In this situation, for $\overline{P}=P^*$ and $s\geq L-1$ we have
\begin{align*}
C \mbb{E} \left[ h^2\left( \overline{P}, \hat{P}_s \right) \right] &\leq  h^2\left( \overline{P}, \mscr{M} \right) + \frac{C(Q^*)}{n(e^{r(Q^*)}-1)} + C(Q^*) e^{-r(Q^*) s}\\
&+ L \epsilon^2 + (s+1) L \overline{V} \frac{ \log n}{n},
\end{align*}
for some positive constant. The condition on $s$ leads to the desired inequality.
\subsubsection{Proof of Proposition \ref{prop:rho_dimension_hmm}}
\label{sec:proof_prop_rho_dimension_hmm}
From Proposition A.1. \cite{robust_mixture}, we have
\[ D^{\mscr{H}_{\delta(s)}} \left( \bigotimes_{i=1}^{n(s,b)} P_i, Q^{\otimes n(s,b)} \right) \leq 545.3 \overline{V} \left[ 5.82  + \log\left( \frac{(K^L+1)^2}{\delta(s)^L}\right) + \log_+\left(\frac{n(s,b)}{\overline{V}}\right)\right].\]
\begin{itemize}
\item If $\overline{V}\leq n(s,1) (K-1)/K$, we have
\begin{align*}
\log\left( \frac{(K^L+1)^2}{\delta(s,b)^L} \right) + \log_+\left( \frac{n(s,b)}{\overline{V}} \right) &\leq \log\left( \frac{(K^L+1)^2 n(s,1)^L (K-1)^L}{\overline{V}^L} \frac{n(s,1)}{\overline{V}} \right)\\
&= \log\left( \frac{(K^L+1)^2 (K-1)^L}{K^{L+1}} \right) + \log\left( \frac{K^{L+1} n(s,1)^{L+1}}{\overline{V}^{L+1}} \right)\\
&= \log\left( \frac{(K^L+1)^2 (K^2-1)^L}{K^{L+1}(K+1)^L} \right) + (L+1) \log\left( \frac{K n(s,1)}{\overline{V}} \right).
\end{align*}
One can check that for $L\geq 2$, we have $\frac{(K^L+1)^2(K^2-1)^L}{K^{L+1}(K+1)^L} \leq K^{2L-1}$ for all $K\geq 1$. Therefore,
\begin{align*}
\log\left( \frac{(K^L+1)^2}{\delta(s)^L} \right) + \log_+\left( \frac{n(s,b)}{\overline{V}} \right) &\leq (2L-1) \log K + (L+1) \log\left( \frac{K n(s,1)}{\overline{V}} \right)\\
&\leq 3L \log\left( \frac{K N}{\overline{V}} \right) = 3L \log\left( \frac{K N}{\overline{V}\wedge N} \right).
\end{align*}
\item Otherwise $\overline{V}> n(s,1) (K-1)/K$ and $\log\left(\frac{K n(s,1)}{\overline{V}\wedge n(s,1)}\right) = \log K$. We have
\begin{align*}
\log\left( \frac{(K^L+1)^2}{\delta(s)^L} \right) + \log_+\left( \frac{n(s,b)}{\overline{V}}\right) &\leq \log\left( \frac{ (K^L+1)^2 K^L  n(s,1)}{\overline{V}}\right)\\
&= \log\left( \frac{K n(s,1)}{\overline{V}} \right) + (L-1)\log K + 2\log\left( 1+K^L \right)\\
&\leq 3L\log\left( \frac{K n(s,1)}{\overline{V}\wedge n(s,1)}\right) + 2\log(1+K^{-L})\\
&\leq 2\log 2 + 3L\log\left( \frac{K n(s,1)}{\overline{V}\wedge n(s,1)}\right). 
\end{align*}
\end{itemize}
\subsubsection{Proof of Proposition \ref{prop:approximation_hmm}}
\label{sec:proof_prop_approximation_hmm}
With Lemma B.3 \cite{robust_mixture}, we have
\[ h\left( P_{w,Q,F}, P_{v,R,G} \right) \leq h\left( wQ^{\bigcirc{L}}, vR^{\bigcirc{L}} \right) + \max_{k_1,\dots,k_L\in [K]^L} h\left( \bigotimes_{l=1}^L F_{k_l}, \bigotimes_{l=1}^L G_{k_l} \right), \]
with 
\begin{equation}
\label{eq:w_Q_L}
wQ^{\bigcirc{L}}(k_1,\dots,k_L) = w_{k_1} Q_{k_1,k_2} \dots Q_{k_{L-1},k_L},\forall k_1,\dots,k_L\in[K]. 
\end{equation}
Let $\rho$ denote the Hellinger affinity defined by $\rho=1-h^2$ For $\rho_-=\min_{k\in[K]} \rho\left(Q_{k,\cdot},R_{k,\cdot}\right)$, we have
\begin{align*}
h^2\left( wQ^{\bigcirc{L}}, vR^{\bigcirc{L}} \right) &= 1 - \rho\left( wQ^{\bigcirc{L}}, vR^{\bigcirc{L}} \right)\\
&= 1 - \suml_{k_1,\dots,k_L} \sqrt{ w_{k_1} v_{k_1} Q_{k_1,k_2} R_{k_1,k_2} \dots Q_{k_{L-1},k_L} R_{k_{L-1},k_L} }\\
&= 1 - \suml_{k_1,\dots,k_{L-1}} \sqrt{ w_{k_1} v_{k_1} Q_{k_1,k_2} R_{k_1,k_2} \dots Q_{k_{L-2},k_{L-1}} R_{k_{L-2},k_{L-1}}} \rho\left(Q_{k_{L-1},\cdot},R_{k_{L-1},\cdot}\right)\\
&\leq 1 - \rho_- \suml_{k_1,\dots,k_{L-1}} \sqrt{ w_{k_1} v_{k_1} Q_{k_1,k_2} R_{k_1,k_2} \dots Q_{k_{L-2},k_{L-1}} R_{k_{L-2},k_{L-1}}}.
\end{align*}
By induction we get
\[ h^2\left( wQ^{\bigcirc{L}}, vR^{\bigcirc{L}} \right) \leq 1 - \rho_-^{L-1} \rho\left( w,v \right) \leq h^2(w,v) + (L-1) \max_{k\in[K]} h^2\left(Q_{k,\cdot},R_{k,\cdot}\right). \]
We also have
\begin{align*}
h^2\left( \bigotimes_{l=1}^L F_{k_l}, \bigotimes_{l=1}^L G_{k_l} \right) &= 1 - \rho\left( \bigotimes_{l=1}^L F_{k_l}, \bigotimes_{l=1}^L G_{k_l} \right)\\
&= 1 - \prod_{l=1}^L \rho\left( F_{k_l}, G_{k_l} \right) \leq \suml_{l=1}^L h^2\left( F_{k_l}, G_{k_l} \right),
\end{align*}
which allows to conclude the proof.
\subsubsection{Proof of Lemma \ref{lem:k_s_hmm}}
\label{sec:proof_lem_k_s_hmm}
Let $s$ not be smaller than $L-1$ and $b$ be in $[s+1]$. Since $\left(Y_i,H_i\right)_{1\leq i\leq N}$ is a hidden Markov model, we have that
\[ \left(X^{(s,b)}_i,H^{(L,s,b)}_i\right)_{1\leq i\leq n} \]
is also a hidden Markov model, with
\[ X^{(s,b)}_i = X_{b+(i-1)(s+1)} \text{  and  } H^{(L,s,b)}_i = \left(H_{b+(i-1)(s+1)},\dots,H_{b+(i-1)(s+1)+L-1}\right). \]
From Lemma \ref{lem:divergence_hmm}, we have
\[ \mbf{K}\left( \mbf{P}^*_{s,b} || \mbf{P}^{ind}_{s,b} \right) \leq \suml_{i=1}^{n(s,b)-1} \mbf{K}\left( \mcal{L}\left(H^{(L,s,b)}_i,H^{(L,s,b)}_{i+1} \right) || \mcal{L}\left(H^{(L,s,b)}_i\right) \otimes \mcal{L}\left(H^{(L,s,b)}_{i+1} \right) \right). \]
We can use the following result to bound the terms in the sum on the right-hand side of the inequality.
\begin{lemme}
\label{lem:reverse_pinsker}
Let $A$ and $B$ be random variables taking values in the finite sets $\mscr{A}$ and $\mscr{B}$ respectively. We have
\begin{align*}
\mbf{K}\left(\mcal{L}(A,B)||\mcal{L}(A)\otimes\mcal{L}(B)\right) &\leq  2 \suml_{a\in\mscr{A}} d_{TV}\left( \mcal{L}(B|A=a), \mcal{L}(B) \right). 
\end{align*}
\end{lemme}
For $k_1,\dots,k_{2L}\in[K^*]$, we have
\begin{align*}
&\mbb{P}\left( H^{(L,s,b)}_{i+1} = (k_{L+1},\dots,k_{2L}) | H^{(L,s,b)}_i=(k_1,\dots,k_L) \right)\\
&= Q^*_{k_{2L-1},k_{2L}} \dots Q^*_{k_{L+1},k_{L+2}} (Q^*)^{s+2-L}_{k_L,k_{L+1}}
\end{align*}
Therefore, we have
\[ \mbf{K}\left( \mbf{P}^*_{s,b} || \mbf{P}^{ind}_{s,b} \right) \leq 2 \suml_{i=1}^{n(s,b)-1} \suml_{k\in [K^*]} d_{TV}\left( (Q^*)^{s+2-L}_{k,\cdot}, \nu_i Q^{s+2-L} \right), \]
where $\nu_i=w^* (Q^*)^{b+(i-1)(s+1)+L-2}$ is the distribution of $H_{b+(i-1)(s+1)+L-1}$. Since $Q^*$ is irreducible and aperiodic, there exists a unique invariant probability $\pi^*$ and there are positive constants $C(Q^*)$ and $r(Q^*)$ such that
\[ d_{TV}\left( (Q^*)^t_{k,\cdot}, \pi^* \right) \leq C(Q^*) e^{-r(Q^*) t}, \forall k\in[K^*],\forall t\geq1. \]
Combining the different inequalities we get
\[ \mbf{K}\left( \mbf{P}^*_{s,b} || \mbf{P}^{ind}_{s,b} \right) \leq 4 K^* (n(s,b)-1) C(Q^*) e^{-r(Q^*) (s+1)}. \]
We have
\[ h^2\left(P^*,P_i\right) \leq d_{TV}\left(P^*,P_i\right) = d_{TV}\left( \pi^*, w^* (Q^*)^{i-1} \right) \leq C(Q^*) e^{-r(Q^*) (i-1)}. \]
\subsubsection{Proof of Lemma \ref{lem:reverse_pinsker}}
We denote by $(\mscr{A}\times \mscr{B})^+$ the set $\left\{ (a,b)\in\mscr{A}\times\mscr{B}; \mbb{P}(A=a,B=b)>0 \right\}$. We have
\begin{align*}
\mbf{K}\left(\mcal{L}(A,B)||\mcal{L}(A)\otimes\mcal{L}(B)\right) &= \suml_{(a,b)\in(\mscr{A}\times\mscr{B})^+} \mbb{P}\left(A=a,B=b\right) \log\left( \frac{ \mbb{P}\left(A=a,B=b\right) }{ \mbb{P}\left(A=a\right) \mbb{P}\left(B=b\right) } \right)\\
&\leq \suml_{(a,b)\in(\mscr{A}\times\mscr{B})^+} \mbb{P}\left(A=a,B=b\right) \left( \frac{ \mbb{P}\left(A=a,B=b\right) }{ \mbb{P}\left(A=a\right) \mbb{P}\left(B=b\right) } -1 \right)\\
&= \suml_{(a,b)\in(\mscr{A}\times\mscr{B})^+} \frac{ \left( \mbb{P}\left(A=a,B=b\right)  - \mbb{P}\left(A=a\right) \mbb{P}\left(B=b\right) \right)^2 }{ \mbb{P}\left(A=a\right) \mbb{P}\left(B=b\right) }.
\end{align*}
For $(a,b)\in(\mscr{A}\times\mscr{B})^+$,
\begin{align*}
&\frac{ \left( \mbb{P}\left(A=a,B=b\right)  - \mbb{P}\left(A=a\right) \mbb{P}\left(B=b\right) \right)^2 }{ \mbb{P}\left(A=a\right) \mbb{P}\left(B=b\right) }\\
&= \left| \mbb{P}\left(A=a|B=b\right)  - \mbb{P}\left(A=a\right) \right| \times \left| \mbb{P}\left(B=b|A=a\right)  - \mbb{P}\left(B=b\right) \right|\\
&\leq \left| \mbb{P}\left(B=b|A=a\right)  - \mbb{P}\left(B=b\right) \right|.
\end{align*}
Finally, we get
\begin{align*}
\mbf{K}\left(\mcal{L}(A,B)||\mcal{L}(A)\otimes\mcal{L}(B)\right) &\leq  \suml_{a\in\mscr{A}} 2 d_{TV}\left( \mcal{L}\left(B|A=a\right), \mcal{L}\left(B\right) \right).
\end{align*}
\subsection{Proof of Corollary \ref{coro:hmm_robust}}
\label{sec:proof_coro_robust_hmm}
We have
\[ \mbb{P}\left( X_i = (Y'_i,\ldots,Y'_{i+L-1}) \right) \geq \mbb{P}\left( E_i=\dots=E_{i+L-1}=1 \right) = p_i p_{i+1} \dots p_{i+L-1}, \]
and with the convexity of the squared Hellinger distance
\begin{align*}
h^2\left( P_i, P^* \right) &\leq p_i p_{i+1} \dots p_{i+L-1} h^2\left( P'_i,P^* \right) + (1-p_i p_{i+1} \dots p_{i+L-1})\\
&\leq h^2\left( P'_i,P^* \right) + (1-p_i) + \dots + (1-p_{i+L-1}),
\end{align*}
where $P'_i=\mcal{L}(Y'_i,\dots,Y'_{i+L-1})$. One can check that $n\geq 1+N/2$ with our conditions on $L$. With Theorem \ref{th:first_hmm}, Lemma \ref{lem:robust_kl} and Lemma \ref{lem:k_s_hmm} we have
\begin{align*}
C \mbb{E} \left[ h^2\left( P^*, \hat{P}_s \right) \right] &\leq h^2\left( P^*, \mscr{M} \right) + \frac{C(Q^*)}{n(e^{r(Q^*)}-1)} + \frac{L}{N} \suml_{i=1}^N (1-p_i)\\
&+  e^{-r(Q^*)s} + L\epsilon^2 + (s+1) L \overline{V} \frac{\log n}{n},
\end{align*}
for some positive constant $C$ and $s\geq L-1$.
\subsection{Proof of Theorems \ref{th:hmm_log_concave} and \ref{th:hmm_log_concave_high}}
\label{sec:proof_th_hmm_log_concave}
With (\ref{eq:v_hmm_log_concave}) and Theorem \ref{th:first_hmm}, we have
\begin{align*}
C\mbb{E} \left[ h^2 \left( \overline{P}, \hat{P}_s \right) \right] &\leq h^2\left( \overline{P}, \mscr{M} \right) + n^{-1} \suml_{i=1}^n h^2\left( P_i,\overline{P} \right)\\
&+ n^{-1} \suml_{b=1}^{s+1} \mbf{K}\left( \mbf{P}^*_{s,b} || \mbf{P}^{ind}_{s,b} \right)\\
&+ L \epsilon^2 + (s+1) L^2 K^L \log(2|\mscr{F}_{\lambda_-,\lambda_+,M}|[\epsilon]) \frac{\log n}{n}.
\end{align*}
We can simply follow the proof of Theorems \ref{th:kolmogorov} and \ref{th:kolmogorov_high} to conclude.
\subsection{Proof of Proposition \ref{prop:exponential_family}}
\label{sec:proof_vc_exponential}
The proof relies on the following lemma.
\begin{lemme}
\label{lem:vc_exponential_family}
The set $\mcal{A}$ of probability density functions, defined by
\[ \mcal{A} = \left\{ (x_1,\dots,x_L) \mapsto q_1(x_1) \dots q_L(x_L); q_i\in\mcal{E}\left(\overline{\Theta}_i,\eta_i,T_i,d_i,B_i\right),\forall i \in\{1,\dots,L\} \right\}, \]
is VC-subgraph with VC-index $3+d_1+\dots+d_L$.
\end{lemme}
As $L\geq 2$ and $\max\limits_{1\leq k\leq K} d_k\geq 2$, Assumption \ref{hyp:hmm_vc} is met with
\[ \overline{V} = 3 K^L + K^{L-1} L \suml_{k=1}^K d_k \leq K^L \left( 3 + L \max\limits_{1\leq k\leq K} d_k \right). \]
\subsubsection{Proof of Lemma \ref{lem:vc_exponential_family}}
We have
\begin{align*}
\mcal{A} &= \left\{ (x_1,\dots,x_L) \mapsto f_{\theta_1}(x_1) \dots f_{\theta_L}(x_L); \theta_i\in\overline{\Theta}_i,\forall i \in\{1,\dots,L\}\right\}\\
&= \exp \circ \left\{ (x_1,\dots,x_L) \mapsto \suml_{i=1}^L \langle \eta_i(\theta_i), T_i(x_i)\rangle + A_i(\theta_i) + B_i(x_i), \forall i \in\{1,\dots,L\} \right\}\\
&\subset \exp \circ \left( V + B \right)
\end{align*}
with $B: (x_1,\dots,x_L) \mapsto B_i(x_1) + \dots + B_i(x_L)$ and
\[V = \left\{ (x_1,\dots,x_L) \mapsto A + \suml_{i=1}^K \langle \eta_i, T_i(x_i)\rangle; \eta_i \in\mbb{R}^d, \forall i \in\{1,\dots,L\}, A\in\mbb{R} \right\}.\]
The set $V$ is a vector space of dimension $1+d_1+\dots+d_L$ and $\exp$ is monotone, therefore, from Proposition 42-(i,ii) \cite{baraudinventiones} and Lemma 2.6.15 \cite{VanDerVaart} and Lemma 2.6.18-(v) \cite{VanDerVaart}, the class of functions $\mcal{A}$ is VC-subgraph with $VC$-index $V(\mcal{A})\leq 3+d_1+\dots+d_L$.
\subsection{Proof of Proposition \ref{prop:regular_parametric_model}}
\label{sec:proof_prop_regular}
We first need the following lemma to apply results of regular parametric models.
\begin{lemme}
\label{lem:regular_statistical_experiment}
Under Assumption \ref{hyp:regular_parametric_model}, our model is regular, i.e.
\begin{itemize}
 \item $\phi\mapsto p(\mbf{x};\phi)$ is continuous for all $\mbf{x}$,
 \item it is differentiable for all $\mbf{x}$,
 \item and the information matrix function
 \[ I : \phi \mapsto I(\phi)=\int_{\mscr{X}^L} \partial_{\phi}p(\mbf{x};\phi) \left( \partial_{\phi}p(\mbf{x};\phi) \right)^T \frac{\mu(\mbf{x})}{p(\mbf{x};\phi)}  \]
 is well-defined and continuous.
\end{itemize}
\end{lemme}
We can now apply results of Ibragimov and Has’minski\u{\i} \cite{Ibragimov}, in particular (7.20) which is a consequence of Theorem 7.6. Let $\kappa$ be a compact subset of $\overline{\Phi}$ such that $\overline{\Phi}$ belongs to the interior of $\kappa$. There is a positive constants $a(\kappa)$ such that
\[ \forall \phi\in \kappa, h^2\left(P_{\phi},P_{\overline{\phi}}\right) \geq a(\kappa) \frac{ ||\phi-\overline{\phi}||^2}{1+||\phi-\overline{\phi}||^2}\geq  \frac{a(\kappa)}{1+b(\kappa)} ||\phi-\overline{\phi}||^2, 
\]
with $b(\kappa)=\max\limits_{\phi\in\kappa} ||\phi-\overline{\phi}||^2$. We know that $c(\kappa):=\inf\limits_{\phi\in \overline{\Phi}\backslash \kappa} h^2\left(P_{\phi},P_{\overline{\phi}} \right)$ is positive. Therefore, there exist a positive constant $C(\overline{\phi})$ such that
\begin{align*}
\forall \phi\in \overline{\Phi}, h^2\left( P_{\phi}, P_{\overline{\phi}} \right) &\geq \mathbbm{1}_{\phi\in\kappa} \frac{a(\kappa)}{1+b(\kappa)} ||\phi-\overline{\phi}||^2 + \mathbbm{1}_{\phi\in \overline{\Phi}\backslash \kappa} c(\kappa)\\
&\geq C(\overline{\phi}) \left[ \left|\left|\overline{w}-w\right|\right|^2 + \left|\left|\overline{Q}-Q\right|\right|^2 + \suml_{k=1}^K \left|\left|\overline{\theta}-\theta\right|\right|^2 \wedge 1 \right].
\end{align*}
\subsubsection{Proof of Lemma \ref{lem:regular_statistical_experiment}}
For $k_1,\dots,k_L\in[K]$ we have
\begin{equation}
\label{eq:p_beta_lower_bound}
p(\mbf{x};\phi) \geq w_{k_1} Q_{k_1,k_2} \dots Q_{k_{L-1},k_L} \prod_{l=1}^L f_{\theta_{k_l}}(x_l).
\end{equation}
\begin{itemize}
\item Since $\eta_k$ and $A_k$ are continuous for all $k$ in $[K]$, then the applications $\theta_k\mapsto f_{\theta_k}(x)$ are continuous for all $x\in\mscr{X}$ and so is $\phi\mapsto p(\mbf{x};\phi)$ for all $\mbf{x}\in\mscr{Y}^L$.
\item The function $u\mapsto p(\mbf{x};u)$ is differentiable at the point $u=\phi$ for all $\mbf{x}\in\mscr{Y}^L$ since $A_k$ and $\eta_k$ are differentiable for all $k\in[K]$. For all $\overline{k}\in[K]$ and $j\in[e_k]$,
\begin{align}
\partial_{\theta_{\overline{k},j}} p(\mbf{x};\phi) &= \suml_{k_1,\dots,k_L} w_{k_1} Q_{k_1,k_2} \dots Q_{k_{L-1},k_L} \suml_{l=1}^L \mathbbm{1}_{k_l=\overline{k}} \left( \prod_{i\neq l} f_{\theta_{k_j}}(x_j) \right) \partial_{\theta_{\overline{k},j}} f_{\theta_{\overline{k}}}(x_l)\nonumber\\
&= \suml_{k_1,\dots,k_L} w_{k_1} Q_{k_1,k_2} \dots Q_{k_{L-1},k_L} \prod_{i=1}^L f_{\theta_{k_i}}(x_i)\nonumber\\
&\times \suml_{l=1}^L \mathbbm{1}_{k_l=\overline{k}} \left[ \langle \partial_{\theta_{\overline{k},j}} \eta_{\overline{k}}(\theta_{\overline{k}}), T_{\overline{k}}(x_l) \rangle + \partial_{\theta_{\overline{k},j}} A_{\overline{k}} (\theta_{\overline{k}}) \right].\label{eq:partial_theta}
\end{align}
For $\overline{k}\in[K-1]$ and $k'\in[K]$ we have
\begin{align}
\partial_{w_{\overline{k}}} p(\mbf{x};\phi) &= \suml_{k_2,\dots,k_L} Q_{\overline{k},k_2} \dots Q_{k_{L-1},k_L} f_{\theta_{\overline{k}}}(x_1) \prod_{l=2}^L f_{\theta_{k_l}}(x_l)\nonumber\\
&- \suml_{k_2,\dots,k_L} Q_{K,k_2} \dots Q_{k_{L-1},k_L} f_{\theta_K}(x_1) \prod_{l=2}^L f_{\theta_{k_l}}(x_l)\label{eq:partial_w}
\end{align}
and
\begin{align}
\partial_{Q_{k',\overline{k}}} p(\mbf{x};\phi) &= \suml_{k_1,k_2,\dots,k_L} w_{k_1} \partial_{Q_{k',\overline{k}}} \left[  Q_{k_1,k_2} \dots Q_{k_{L-1},k_L} \right] \prod_{l=1}^L f_{\theta_{k_l}}(x_l)\nonumber\\
&= \suml_{k_1,k_2,\dots,k_L} w_{k_1} \prod_{i=1}^L f_{k_i,\theta_{k_i}}(x_i) \suml_{l=2}^L \left[ \mathbbm{1}_{(k',\overline{k})=(k_{l-1},k_l)} - \mathbbm{1}_{(k',K)=(k_{l-1},k_l)} \right] \prod_{\substack{2\leq j\leq L,\\j\neq l}} Q_{k_{j-1},k_j}. \label{eq:partial_Q}
\end{align}
\end{itemize}
Since $A_k$ and $\eta_k$ are $\mcal{C}^1$, we just need to check that the functions
\begin{align}
\label{eq:t_2}
\phi &\mapsto \int_{\mscr{Y}^L} T_{\overline{k},j}(x_i) T_{\overline{k}',j'}(x_{i'})  \prod_{l=1}^L f_{\theta_{k_l}}(x_l) f_{\theta_{k'_l}}(x_l)  \frac{\mu(d\mbf{x})}{p(\mbf{x};\phi)},\\
\label{eq:t_1}
\phi &\mapsto \int_{\mscr{Y}^L} T_{\overline{k},j}(x_i) \prod_{l=1}^L f_{\theta_{k_l}}(x_l) f_{\theta_{k'_l}}(x_l)  \frac{\mu(d\mbf{x})}{p(\mbf{x};\phi)},\\
\label{eq:t_0}
\phi &\mapsto \int_{\mscr{Y}^L} \prod_{l=1}^L f_{\theta_{k_l}}(x_l) f_{\theta_{k'_l}}(x_l) \frac{\mu(d\mbf{x})}{p(\mbf{x};\phi)},
\end{align}
are well-defined and continuous for all $k_1,k'_1,\ldots,k_L,k'_L,\overline{k},\overline{k}\in[K],j\in[d_{\overline{k}}],j'\in[d_{\overline{k}'}],i,i'\in[L]$, where 
\[ T_k(x)=(T_{k,1}(x),\ldots,T_{k,d_k}(x))\in\mbb{R}^{d_k},\forall x\in\mscr{Y}. \]
We deal with integrability in the first time and then look at continuity, using (\ref{eq:p_beta_lower_bound}) repeatedly.
\begin{itemize}
\item We have
\begin{align*}
0 &\leq \int_{\mscr{Y}^L} \prod_{l=1}^L f_{\theta_{k_l}}(x_l) f_{\theta_{k'_l}}(x_l) \frac{\mu(d\mbf{x})}{p(\mbf{x};\phi)}\\
&\leq \left( w_{k_1} Q_{k_1,k_2} \dots Q_{k_{L-1},k_L} \right)^{-1} \int_{\mscr{Y}^L} \prod_{l=1}^L f_{\theta_{k'_j}}(x_j) \mu(d\mbf{x})\\
&= \left( w_{k_1} Q_{k_1,k_2} \dots Q_{k_{L-1},k_L} \right)^{-1} <\infty,
\end{align*}
and (\ref{eq:t_0}) is well defined. Similarly
\begin{align*}
0 &\leq \int_{\mscr{Y}^L} \left| T_{\overline{k},j}(x_i) \right| \prod_{l=1}^L f_{\theta_{k_l}}(x_l) f_{\theta_{k'_l}}(x_l) \frac{\mu(d\mbf{x})}{p(\mbf{x};\phi)}\\
&\leq \left( w_{k'_1} Q_{k'_1,k'_2} \dots Q_{k'_{L-1},k'_L} \right)^{-1} \int_{\mscr{Y}}  \left| T_{\overline{k},j}(x_i) \right| f_{\theta_{k_i}}(x_i) \nu(d x_i)\\
&\leq \left( w_{k'_1} Q_{k'_1,k'_2} \dots Q_{k'_{L-1},k'_L} \right)^{-1} \sqrt{ \int_{\mscr{Y}}  \left| T_{\overline{k},j}(x_i) \right|^2 f_{\theta_{k_i}}(x_i) \nu(d x_i) } < \infty,
\end{align*}
and (\ref{eq:t_1}) is well defined. Finally
\begin{align*}
0 &\leq \int_{\mscr{Y}^L} \left| T_{\overline{k},j}(x_i) T_{\overline{k}',j'}(x_{i'}) \right| \prod_{l=1}^L f_{\theta_{k_l}}(x_l) f_{\theta_{k'_l}}(x_l) \frac{\mu(d\mbf{x})}{p(\mbf{x};\phi)}\\
&\leq \left( w_{k_1} w_{k'_1} Q_{k_1,k_2} Q_{k'_1,k'_2} \dots Q_{k_{L-1},k_L} Q_{k'_{L-1},k'_L} \right)^{-1/2}\\
&\times \int_{\mscr{Y}^L}  \left| T_{\overline{k},j}(x_i) T_{\overline{k}',j'}(x_{i'}) \right| \sqrt{ \prod_{l=1}^L f_{\theta_{k_l}}(x_l) f_{\theta_{k'_l}}(x_l) } \mu(d\mbf{x})\\
&\leq \left( w_{k_1} w_{k'_1} Q_{k_1,k_2} Q_{k'_1,k'_2} \dots Q_{k_{L-1},k_L} Q_{k'_{L-1},k'_L} \right)^{-1/2}\\
&\times \sqrt{ \int_{\mscr{Y}} \left| T_{\overline{k},j}(x_i) \right|^2  f_{\theta_{k_i}}(x_i) \nu(dx_i) } \sqrt{ \int_{\mscr{Y}} \left| T_{\overline{k}',j'}(x_{i'}) \right|^2  f_{\theta_{k'_{i'}}}(x_{i'}) \nu(dx_{i'}) } < \infty,
\end{align*}
and (\ref{eq:t_2}) is well defined. The Fisher information matrix $I(\phi)$ is well-defined for all $\phi$. We now turn to continuity.
\item  We have
\begin{align*}
&\left| \frac{ \prod_{l=1}^L f_{\theta_{k_l}}(x_l) f_{\theta_{k'_l}}(x_l) }{p(\mbf{x};\phi)} - \frac{ \prod_{l=1}^L f_{k_l,\theta'_{k_l}}(x_l) f_{k'_l,\theta'_{k'_l}}(x_l) }{p(\mbf{x};\phi')} \right|\\
&\leq \frac{\prod_{l=1}^L f_{\theta_{k'_l}}(x_l)}{p(\mbf{x};\phi)} \left| \prod_{l=1}^L f_{\theta_{k_l}}(x_l) - \prod_{l=1}^L f_{\theta'_{k_l}}(x_l) \right|\\
&+ \prod_{l=1}^L f_{\theta'_{k_l}}(x_l) \prod_{l=1}^L f_{\theta_{k'_l}}(x_l) \left| \frac{1}{p(\mbf{x};\phi)} - \frac{1}{p(\mbf{x};\phi')} \right|\\
&+ \frac{\prod_{l=1}^L f_{\theta'_{k_l}}(x_l)}{p(\mbf{x};\phi')} \left| \prod_{l=1}^L f_{\theta_{k'_l}}(x_l) - \prod_{l=1}^L f_{\theta'_{k'_l}}(x_l) \right| \\
&\leq \frac{ \left| \prod_{l=1}^L f_{\theta_{k_l}}(x_l) - \prod_{l=1}^L f_{\theta'_{k_l}}(x_l) \right| }{ w_{k'_1} Q_{k'_1,k'_2} \dots Q_{k'_{L-1},k'_L} } \\
&+ \frac{ \left| p(\mbf{x};\phi) - p(\mbf{x};\phi') \right| }{ w'_{k_1} w_{k'_1} Q'_{k_1,k_2} Q_{k'_1,k'_2} \dots Q'_{k_{L-1},k_L} Q_{k'_{L-1},k'_L} }\\
&+ \frac{ \left| \prod_{l=1}^L f_{\theta_{k'_l}}(x_l) - \prod_{l=1}^L f_{\theta'_{k'_l}}(x_l) \right| }{ w'_{k_1} Q'_{k_1,k_2} \dots Q'_{k_{L-1},k_L} }.
\end{align*}
Therefore,
\begin{align*}
&\left| \int_{\mscr{Y}^L} \frac{ \prod_{l=1}^L f_{\theta_{k_l}}(x_l) f_{\theta_{k'_l}}(x_l) }{p(\mbf{x};\phi)} \mu(d\mbf{x}) - \int_{\mscr{Y}^L} \frac{ \prod_{l=1}^L f_{\theta'_{k_l}}(x_l) f_{\theta'_{k'_l}}(x_l) }{p(\mbf{x};\phi')} \mu(d\mbf{x}) \right|\\
&\leq \frac{2 d_{TV}\left( \bigotimes_{l=1}^L F_{\theta_{k_l}}, \bigotimes_{l=1}^L F_{\theta'_{k_l}} \right) }{ w_{k'_1} Q_{k'_1,k'_2} \dots Q_{k'_{L-1},k'_L} } \\
&+ \frac{ 2 d_{TV}\left( P_{\phi}, P_{\phi'} \right) }{ w'_{k_1} w_{k'_1} Q'_{k_1,k_2} Q_{k'_1,k'_2} \dots Q'_{k_{L-1},k_L} Q_{k'_{L-1},k'_L} }\\
&+ \frac{2 d_{TV}\left( \bigotimes_{l=1}^L F_{\theta_{k'_l}}, \bigotimes_{l=1}^L F_{\theta'_{k'_l}} \right) }{ w'_{k_1} Q'_{k_1,k_2} \dots Q'_{k_{L-1},k_L} }.
\end{align*}
Since convergence with respect to the total variation distance and to the Hellinger distance are equivalent, we get continuity of (\ref{eq:t_0}) with Proposition \ref{prop:approximation_hmm}. Similarly, we have
\begin{align*}
&\left| \int_{\mscr{Y}^L} \frac{ T_{\overline{k},j}(x_i) \prod_{l=1}^L f_{\theta_{k_l}}(x_l) f_{\theta_{k'_l}}(x_l) }{ p(\mbf{x};\phi) } \mu(d\mbf{x}) - \int_{\mscr{Y}^L} \frac{ T_{\overline{k},j}(x_i) \prod_{l=1}^L f_{\theta'_{k_l}}(x_l) f_{\theta'_{k'_l}}(x_l) }{ p(\mbf{x};\phi') } \mu(d\mbf{x}) \right|\\
&\leq \frac{ \int_{\mscr{Y}^L} |T_{\overline{k},j}(x_i)|  \left| \prod_{l=1}^L f_{\theta_{k_l}}(x_l) - \prod_{l=1}^L f_{\theta'_{k_l}}(x_l) \right| \mu(d\mbf{x}) }{ w_{k'_1} Q_{k'_1,k'_2} \dots Q_{k'_{L-1},k'_L} } \\
&+ \frac{ \int_{\mscr{Y}^L} |T_{\overline{k},j}(x_i)| \left| p(\mbf{x};{\phi}) - p(\mbf{x};{\phi'}) \right| \mu(d\mbf{x}) }{ w'_{k_1} w_{k'_1} Q'_{k_1,k_2} Q_{k'_1,k'_2} \dots Q'_{k_{L-1};k_L} Q_{k'_{L-1},k'_L} }\\
&+ \frac{ \int |T_{k_l}(x_l)| \left| \prod_{i=l}^L f_{\theta_{k'_l}}(x_l) - \prod_{i=1}^L f_{k'_l,\theta'_{k'_l}}(x_l) \right| \mu(d\mbf{x}) }{ w'_{k_1} Q'_{k_1,k_2} \dots Q'_{k_{L-1},k_L} }.
\end{align*}
We have
\begin{align*}
&\int_{\mscr{Y}^L} |T_{\overline{k},j}(x_i)| \left| p(\mbf{x};{\phi}) - p(\mbf{x};{\phi'}) \right| \mu(d\mbf{x})\\
&\leq \suml_{1\leq k_1,\ldots,k_L\leq K}  \int_{\mscr{Y}^L} |T_{\overline{k},j}(x_i)| \left| \prod_{l=1}^L f_{\theta_{k_l}}(x_l) - \prod_{l=1}^L f_{\theta'_{k_l}}(x_l) \right| \mu(d\mbf{x})
\end{align*}
and
\begin{align*}
&\int_{\mscr{Y}^L} |T_{\overline{k},j}(x_i)|  \left| \prod_{l=1}^L f_{\theta_{k_l}}(x_l) - \prod_{l=1}^L f_{\theta'_{k_l}}(x_l) \right| \mu(d\mbf{x})\\
&\leq \int_{\mscr{Y}^L} |T_{\overline{k},j}(x_i)|  \left| f_{\theta_{k_i}}(x_i) - f_{\theta'_{k_i}}(x_i) \right| \nu(d x_i)\\
&+ 2 \int_{\mscr{Y}} |T_{\overline{k},j}(x_i)| f_{\theta_{k_i}}(x_i) \nu(dx_i) \times \suml_{l<i}  d_{TV} \left( F_{\theta_{k_l}}, F_{\theta'_{k_l}} \right)\\
&+ 2 \int_{\mscr{Y}} |T_{\overline{k},j}(x_i)| f_{\theta'_{k_i}}(x_i) \nu(dx_i) \times \suml_{l>i} d_{TV} \left( F_{\theta_{k_l}}, F_{\theta'_{k_l}} \right).
\end{align*}
As
\begin{align*}
&\int_{\mscr{Y}} |T_{\overline{k},j}(x)| \left| f_{\theta_k}(x) - f_{\theta'_k}(x) \right| \nu(dx)\\
&\leq \sqrt{ \int_{\mscr{Y}} |T_{\overline{k},j}(x)|^2 \left| f_{\theta_k}(x) - f_{\theta'_k}(x) \right| \nu(dx)} \times \sqrt{ 2 d_{TV}\left( F_{\theta_k}, F_{\theta'_k} \right) } \xrightarrow[\theta'_k\rightarrow \theta_k]{} 0.
\end{align*}
for all $k\in[K]$ and $\theta_k\in\Theta_k$, we get continuity of (\ref{eq:t_1}). Similarly, we only need
\[ \int_{\mscr{Y}} |T_{\overline{k},j}(x)|^2 \left| f_{\theta_k}(x) - f_{\theta'_k}(x) \right| \nu(dx) \xrightarrow[\theta'_k\rightarrow \theta_k]{} 0   \]
to obtain the continuity of (\ref{eq:t_2}). 
\end{itemize}
\subsection{Proof of Theorem \ref{th:ibragimov_exponential}}
\label{sec:proof_th_ibragimov_exponential}
We start the proof with two lemmas that ensure we fit into the framework of Proposition \ref{prop:regular_parametric_model}.
\begin{lemme}
\label{lem:definite_positive_information} 
The information matrix $I(\phi)$ is definite positive for all $\phi$ in $\overline{\Phi}$. 
\end{lemme}
\begin{lemme}
\label{lem:identifiable_suite}
Let $(\phi_n)_{n\in\mbb{N}}$ be a sequence in $\overline{\Phi}$. If $\lim\limits_{n\rightarrow \infty} h\left( P_{\phi_n},P_{\overline{\phi}}\right)=0$, then we have $\lim\limits_{n\rightarrow \infty}\phi_n=\overline{\phi}$.
\end{lemme}
One can see that Lemma \ref{lem:identifiable_suite} implies that
$\inf\limits_{\substack{||\phi-\overline{\phi}||\geq a\\\phi\in \overline{\Phi}}} h^2\left(P_{\phi},P_{\overline{\phi}}\right)>0$ for all $a>0$. Therefore we can apply Proposition \ref{prop:regular_parametric_model}. From Proposition \ref{prop:exponential_family}, we get $\overline{V}\leq (3+L)K^L=5 K^3$.
\subsubsection{Proof of Lemma \ref{lem:definite_positive_information}}
For $\mbf{k}=(k_1,\dots,k_L)\in[K]^L$, the notation $wQ^{\bigcirc L}(\mbf{k})$ is defined by (\ref{eq:w_Q_L}). Following Theorem 1 of Meijer \& Ypma \cite{meijer}, we have
\begin{align*}
\det(I(\phi))=0 &\Leftrightarrow \exists \lambda\neq 0, \suml_i \lambda_i \partial_{\phi_i} p(\mbf{x};\phi) = 0 \text{ for } \mu\text{-almost all }\mbf{x}.  
\end{align*}
We can use (\ref{eq:partial_theta}), (\ref{eq:partial_w}) and (\ref{eq:partial_Q}) to get
\begin{align*}
0 &= \suml_{\mbf{k}\in[K]^L} wQ^{\bigcirc L}(\mbf{k}) \prod_{l=1}^L f_{\theta_{k_l}}(x_l) \suml_{l=1}^L \suml_{j=1}^{e_{k_l}} \lambda_{\theta_{k_l,j}} \left[ \langle \partial_{\theta_{k_l,j}} \eta_{k_l}(\theta_{k_l}), T_{k_l}(x_l) \rangle + \partial_{\theta_{k_l,j}} A_{k_l}(\theta_{k_l}) \right]\\
&+ \suml_{k_1=1}^{K-1} \lambda_{w_{k_1}} \left[ f_{\theta_{k_1}}(x_1) - f_{\theta_K}(x_1) \right] \suml_{k_2,\dots,k_L}  \frac{ wQ^{\bigcirc L}(\mbf{k}) }{ w_{k_1} } \prod_{i=2}^L f_{\theta_{k_i}}(x_i)\\
&+ \suml_{l=2}^L \suml_{k_l=1}^{K-1} \suml_{k_1,\dots,k_{l-1},k_{l+1},\dots,k_L} \lambda_{Q_{k_{l-1},k_l}} \frac{ wQ^{\bigcirc L}(\mbf{k}) }{ Q_{k_{l-1},k_l} } \left[ f_{\theta_{k_l}}(x_l) -  f_{\theta_K}(x_l)\right] \prod_{i\neq l} f_{\theta_k}(x_i), 
\end{align*}
for almost all $x$. If we apply it to exponential distributions, we get
\begin{align}
0 &= - \suml_{\mbf{k}\in[K]^L} wQ^{\bigcirc L}( \mbf{k} ) \theta_{k_1}\dots \theta_{k_L} e^{-\theta_{k_1}x_1 -\dots-\theta_{k_L}x_L} \left(\suml_{l=1}^L \lambda_{\theta_{k_l}} x_l  \right)\label{eq:aux_information_exp}\\
&+ \suml_{k_1=1}^{K-1} \lambda_{w_{k_1}}  \suml_{k_2,\dots,k_L} \frac{ wQ^{\bigcirc L}(\mbf{k}) }{ w_{k_1} } \theta_{k_1}\dots\theta_{k_L} e^{-\theta_{k_1}x_1-\dots-\theta_{k_L}x_L}\nonumber\\
&- \suml_{k_1=1}^{K-1} \lambda_{w_{k_1}} \suml_{k_2,\dots,k_L} \frac{ wQ^{\bigcirc L}(\mbf{k}) }{ w_{k_1} } \theta_K \theta_{k_2}\dots\theta_{k_L} e^{-\theta_K x_1-\dots-\theta_{k_L}x_L}\nonumber\\
&+ \suml_{l=2}^L \suml_{k_l=1}^{K-1} \suml_{k_i;i\neq l}  \lambda_{Q_{k_{l-1},k_l}} \frac{ wQ^{\bigcirc L}(\mbf{k}) }{ Q_{k_{l-1},k_l} } \theta_{k_1}\dots\theta_{k_L} e^{-\theta_{k_1}x_1-\dots-\theta_{k_l}x_l}\nonumber\\
&- \suml_{l=2}^L \suml_{k_l=1}^{K-1} \suml_{k_i;i\neq l}  \lambda_{Q_{k_{l-1},k_l}} \frac{ wQ^{\bigcirc L}(\mbf{k}) }{ Q_{k_{l-1},k_l} } \theta_{k_1}\dots \theta_{k_{l-1}} \theta_K \theta_{k_{l+1}} \dots \theta_{k_L} e^{-\theta_{k_1}x_1-\dots-\theta_K x_l-\dots-\theta_{k_l}x_l}.\nonumber 
\end{align}
As $\theta_1>\dots>\theta_K$, we can identify the coefficients for each $x\mapsto e^{-\theta_{k_1}x_1-\dots-\theta_{k_L}x_L}$. For $\mbf{k}\in[K-1]^L$, we get
\begin{align*}
&0 = -wQ^{\bigcirc L}(\mbf{k}) \theta_{k_1}\dots\theta_{k_L} \left( \suml_{l=1}^L \lambda_{\theta_{k_l}} x_l \right) + \lambda_{w_{k_1}} \frac{ wQ^{\bigcirc L}(\mbf{k}) }{ w_{k_1} } \theta_{k_1}\dots\theta_{k_L}\\
&+ \suml_{l=2}^L \lambda_{Q_{k_{l-1},k_l}} \frac{ wQ^{\bigcirc L}(\mbf{k}) }{ Q_{k_{l-1},k_l} } \theta_{k_1}\dots\theta_{k_L} \text{ for almost all }\mbf{x}\\
&\Rightarrow 0 = \lambda_{\theta_{k_1}} = \dots = \lambda_{\theta_{k_L}} = \frac{\lambda_{w_{k_1}}}{w_{k_1}} + \suml_{l=2}^L \frac{ \lambda_{Q_{k_{l-1},k_l}} }{ Q_{k_{l-1},k_l}}.
\end{align*}
This implies $\lambda_{\theta_k}=0$ for all $k\in[K-1]$ and there are quantities $\lambda^*_w$ and $\lambda^*_Q$ such that $\frac{\lambda_{w_k}}{w_k} = \lambda^*_k$ for all $k\in[K-1]$ and $\frac{\lambda_{ Q_{k_1,k_2}} }{ Q_{k_1,k_2} }= \lambda^*_Q$ for $k_1,k_2\in[K-1]$ and $\lambda^*_w + (L-1)\lambda^*_Q=0$. Therefore, (\ref{eq:aux_information_exp}) becomes
\begin{align}
0 &= \lambda^*_w \suml_{k_1=1}^{K-1} \suml_{k_2,\dots,k_L} wQ^{\bigcirc L}(\mbf{k}) \theta_{k_1}\dots\theta_{k_L} e^{-\theta_{k_1}x_1-\dots-\theta_{k_L}x_L}\label{eq:aux_information_exp_2}\\
&- \lambda^*_w \suml_{k_2,\dots,k_L} \left( \suml_{k_1=1}^{K-1} wQ^{\bigcirc L}(\mbf{k})  \right) \theta_K \theta_{k_2}\dots\theta_{k_L} e^{-\theta_K x_1-\dots-\theta_{k_L}x_L}\nonumber\\
&+ \suml_{l=2}^L \suml_{k_l=1}^{K-1} \suml_{k_i:i\neq l}  \lambda_{Q_{k_{l-1},k_l}} \frac{ wQ^{\bigcirc L}(\mbf{k}) }{ Q_{k_{l-1},k_l} } \theta_{k_1}\dots\theta_{k_L} e^{-\theta_{k_1}x_1-\dots-\theta_{k_l}x_l}\nonumber\\
&- \suml_{l=2}^L \suml_{k_l=1}^{K-1} \suml_{k_i:i\neq l}  \lambda_{ Q_{k_{l-1},k_l} } \frac{ wQ^{\bigcirc L}(\mbf{k}) }{ Q_{k_{l-1},k_l} } \theta_{k_1}\dots\theta_K\dots \theta_{k_L} e^{-\theta_{k_1}x_1-\dots-\theta_K x_l-\dots-\theta_{k_l}x_l}.\nonumber 
\end{align}
For $k_2,\dots,k_L\in[K-1]^{L-1}$, we write $\mbf{k}'=(K,k_2,\dots,k_L)$ and with identification with respect to $\mbf{x}\mapsto e^{-\theta_K x_1 - \theta_{k_2} x_2 -\dots-\theta_{k_L}x_L}$ we have
\begin{align*}
&0 = - \lambda^*_w \left(\suml_{k_1=1}^{K-1} wQ^{\bigcirc L}(\mbf{k}) \right) \theta_K \theta_{k_2}\dots\theta_{k_L} + \lambda_{Q_{K,k_2}} \frac{ wQ^{\bigcirc L}(\mbf{k}') }{ Q_{K,k_2} } \theta_K \theta_{k_2} \dots \theta_{k_L}\\
\Rightarrow& \lambda^*_w \left(\suml_{k_1=1}^{K-1} w_{k_1} Q_{k_1,k_2} \right) = \frac{ \lambda_{Q_{K,k_2}} }{ Q_{K,k_2} } w_K Q_{K,k_2}.
\end{align*}
For $k\in[K-1]$, 
\begin{equation}
\label{eq:regular_parametric_model_aux_1}
\frac{ \lambda_{Q_{K,k}} }{ Q_{K,k} } =  \lambda^*_w \beta_k \text{  with  } \beta_k=\frac{ \suml_{k'=1}^{K-1} w_{k'} Q_{k',k} }{ w_K Q_{K,k} }.
\end{equation}
Finally (\ref{eq:aux_information_exp_2}) becomes
\begin{align*}
0 &= \lambda^*_w \suml_{k_1=1}^{K-1} \suml_{k_2,\dots,k_L} wQ^{\bigcirc L}(\mbf{k}) \theta_{k_1}\dots\theta_{k_L} e^{-\theta_{k_1}x_1-\dots-\theta_{k_L}x_L}\nonumber\\
&- \lambda^*_w \suml_{k_2,\dots,k_L} \left( \suml_{k_1=1}^{K-1} wQ^{\bigcirc L}(\mbf{k}) \right) \theta_K \theta_{k_2}\dots\theta_{k_L} e^{-\theta_K x_1-\dots-\theta_{k_L}x_L}\\
&+ \lambda^*_Q \suml_{l=2}^L \suml_{k_{l-1},k_l\in[K-1]} \suml_{\substack{k_i\in[K];\\i\notin\{l-1,l\}}}  wQ^{\bigcirc L}(\mbf{k}) \theta_{k_1}\dots\theta_{k_L} e^{-\theta_{k_1}x_1-\dots-\theta_{k_L}x_L}\\
&+ \lambda^*_w \suml_{l=2}^L \suml_{k_l\in[K-1]} \suml_{\substack{k_i\in[K];\\i\neq l}}  \beta_{k_l} wQ^{\bigcirc L}(\mbf{k}) \theta_{k_1}\dots \theta_{k_{l-2}} \theta_K \theta_{k_l} \dots \theta_{k_L} e^{-\theta_{k_1}x_1-\dots-\theta_Kx_{l-1}-\theta_{k_l}x_l-\dots -\theta_{k_L}x_L}\\
&- \lambda^*_Q \suml_{l=2}^L \suml_{ k_{l-1},k_l\in[K-1]} \suml_{\substack{k_i\in[K];\\i\notin\{l-1,l\}}}  wQ^{\bigcirc L}(\mbf{k}) \theta_{k_1}\dots\theta_{k_{l-1}}\theta_K \theta_{k_{l+1}}\dots \theta_{k_L} e^{-\theta_{k_1}x_1-\dots-\theta_{k_{l-1}}x_{l-1}-\theta_K x_l-\dots-\theta_{k_L}x_L}\\
&- \lambda^*_w \suml_{l=2}^L \suml_{k_l=1}^{K-1} \suml_{\substack{k_i\in[K];\\i\neq l}}  \beta_{k_l} wQ^{\bigcirc L}(\mbf{k}) \theta_{k_1}\dots\theta_{k_{l-2}}\theta_K\theta_K\theta_{k_{l+1}}\dots \theta_{k_L} e^{-\theta_{k_1}x_1-\dots-\theta_Kx_{l-1}-\theta_K x_l-\dots-\theta_{k_L}x_L}. 
\end{align*}
Identification with respect to $\mbf{x}\mapsto e^{-\theta_Kx_1\dots-\theta_Kx_K}$ gives
\begin{align*}
&0 = - \lambda^*_w \left( \suml_{k=1}^{K-1} w_{k_1} \right) Q_{K,K}^{L-1} - \lambda^*_w \suml_{l=2}^{L-1} \suml_{k_l=1}^{K-1} \beta_{k_l} w_K Q_{K,K}^{L-3} Q_{k_l,K} Q_{K,k_l}) - \lambda^*_w \suml_{k_L=1}^{K-1} \beta_{k_L} w_K Q_{K,K}^{L-2} Q_{K,k_L}\\
\Rightarrow & 0= \lambda^*_w \left[ (1-w_K) Q_{K,K}^2 + (L-2)  \suml_{k_2=1}^{K-1} w_K \beta_{k_2} Q_{k_2,K} Q_{K,k_2} + Q_{K,K} \suml_{k_2=1}^{K-1} w_K \beta_{k_2} Q_{K,k_2}  \right]\\
\Rightarrow & 0= \lambda^*_w \left[ (1-w_K) Q_{K,K}^2 + (L-2)  \suml_{k_2=1}^{K-1} \left( \suml_{k_1}^{K-1} w_{k_1} Q_{k_1,k_2} \right) Q_{k_2,K} + Q_{K,K} \suml_{k_2=1}^{K-1} \suml_{k_1=1}^{K-1} w_{k_1} Q_{k_1,k_2} \right],
\end{align*}
where the last inequality comes from the definition of $\beta_k$. One can notice the quantity between the brackets is positive as a consequence of the definition of $O_K$. Therefore, we necessarily have $\lambda^*_w=0$ and consequently $\lambda^*_Q=\lambda{K,1}=\dots=\lambda_{K,K-1}=0$ which means $\lambda=0$ and therefore the information matrix is definite positive.
\subsubsection{Proof of Lemma \ref{lem:identifiable_suite}}
The parameters $w_k$ and $Q_{k,k'}$ are bounded so we can assume the sequences $w_{k,n}$ and $Q_{{k,k'}_n}$ are converging, with respective limits $w^*_k$ and $Q^*_{k,k'}$, even if it means extracting a subsequence. For other parameters, it is always possible to extract a subsequence $\phi_{\psi(n)}$ such that for all $k$ in $[K]$, we have $\theta_{k,\psi(n)} \xrightarrow[n\rightarrow \infty]{} \theta^*_k \in [0,\infty]$. We can deduce from the definition of $\overline{\Phi}$ that $\theta^*_1\geq \theta^*_2 \geq \dots\geq\theta^*_K$. Let us consider the following cases, dropping the dependency on $\psi$ in the notation.
\begin{itemize}
\item If $\theta^*_k=+\infty$, we have $\theta_{k,n} e^{-\theta_{k,n}x}\cdot dx \xrightarrow[n\rightarrow \infty]{\mbb{P}} \text{Dirac}(0)$. Since $\lim_{n\rightarrow\infty} h\left(P_{\phi,n},P_{\overline{\phi}}\right)$, we get that $w^*_{k_1} Q^*_{k_1,k_2} \dots Q_{k_{L-1},L} = 0$ if $k$ appears in $k_1,k_2,\dots,k_L$.
\item If $\theta^*_k=0$. We have
\[ P_{\overline{\phi}}\left( [\theta_{k,n}^{-1},+\infty)^L \right) \leq (e^{- \overline{\theta}_K/\theta_{k,n} } )^L \xrightarrow[n\rightarrow \infty]{} 0, \]
and
\[ P_{\phi_n}\left( [\theta_{k,n},+\infty)^L \right) \geq w_{k_n} Q_{k_n,k_n}^{L-1} e^{-L}. \]
Since $\lim\limits_{n\rightarrow \infty} h\left( P_{\phi_n},P_{\overline{\phi}}\right)=0$, we must have $w_k^* (Q^*_{k,k})^{L-1}=0$.
\end{itemize}
This proves that $P_{\phi_n}$ converges to
\[ P_{\infty}(d\mbf{x}) = \suml_{ k_1,\dots,k_L \in [K]^+ } w^*_{k_1} Q^*_{k_1,k_2} \dots Q^*_{k_{L-1},k_L} \theta^*_{k_l} \prod_{l_=1}^L e^{-\theta^*_{k_l}x_l} dx_1 \dots dx_L , \]
with $[K]^+=\left\{ k \in [K]; \theta^*_k \in(0,\infty) \right\}$, and necessarily $P_{\infty}=P_{\overline{\phi}}$. We can easily identify the different parameters which implies that $(w^*,Q^*,\theta^*)$ and $(\overline{w},\overline{Q},\overline{\theta})$ are equal up to a permutation $\sigma$ on $[K]$. The ordering of the $\overline{\theta}_k$ and the $\theta^*_k$ ensures that this equality is true, not even up to a permutation.
\subsection{Proof of Theorem \ref{th:faster_rates}}
\label{sec:proof_th_faster_rates}
We just need to check that we satisfy Assumption \ref{hyp:hmm_vc}. Then we can combine Proposition \ref{prop:faster_rates} and ??. We use Definition 41 \cite{baraudinventiones} that allows to consider functions taking values in $(-\infty,+\infty]$. From Lemma 2.6.15 \cite{VanDerVaart}, we have that 
\[ \left\{ \mbf{x}\mapsto (x_1-z_1)(x_2-z_2); z_1,z_2\in\mbb{R} \right\} \subset \left\{ \mbf{x} \mapsto ax_1 + bx_2 + x_1x_2 + c; a,b,c\in\mbb{R} \right\} \]
is VC-subgraph with VC-dimension smaller than or equal to 4. With Proposition 42-$(v)$ \cite{baraudinventiones}, we get that $\left\{ \mbf{x}\mapsto |x_1-z_1| \cdot |x_2-z_2|; z_1,z_2\in\mbb{R} \right\}$ is VC-subgraph with VC-dimension not larger than 37.608. We now need the following result.
\begin{lemme}
\label{lem:vc_set_to_function}
If $\mscr{A}\subset\mcal{P}\left(\mscr{X}\right)$ is a VC-class with dimension $V$, then $\mscr{F}_{\mscr{A},a} := \left\{ p_{A,a}; A\in\mscr{A} \right\}$ is VC-subgraph with dimension $V$ for any $a$ in $\mbb{R}$ where
\[ p_{A,a}(x) := \left\{ \begin{array}{l}
a \text{  if } x\in A,\\
+ \infty \text{  otherwise}.
                    \end{array}
 \right.\] 
\end{lemme}
Since $\mscr{C} := \left\{ C_{z_1,z_2}:=[z_1\pm 1]\times[z_2\pm 1]; z_1,z_2\in\mbb{R} \right\}$ is VC with VC-dimension $4$, we get that $\mscr{F}_{\mscr{C},0}$ is VC-subgraph with VC-dimension $4$. We can apply Proposition 42-$(v)$ \cite{baraudinventiones} one more time which implies that $\mscr{G} = \left\{ \mbf{x} \mapsto g_{z_1,z_2}(\mbf{x}) ; z_1,z_2\in\mbb{R} \right\}$ is VC-subgraph with dimension at most $4.701(37.608+4)\leq 196$, with
\begin{align*}
g_{z_1,z_2} (\mbf{x}) &:= p_{C_{z_1,z_2},0} \vee |x_1-z_1| \cdot |x_2-z_2|\\
&= \left\{ \begin{array}{l}
|x_1-z_1|\cdot|x_2-z_2| \text{  if } x\in[z_1\pm 1]\times[z_2\pm 1],\\
+ \infty \text{  otherwise}.
                    \end{array}\right.
\end{align*}
We need another lemma before we have a bound on the VC-dimension of
\[ \mscr{S}_{\alpha,2} := \left\{ \mbf{x} \mapsto f_{\alpha}(x_1-z_1) f_{\alpha}(x_2-z_2) = \frac{(1-\alpha)^2}{4} \frac{1}{g_{z_1,z_2}^{\alpha}(\mbf{x})}; z_1,z_2\in\mbb{R} \right\}. \]
\begin{lemme}
\label{lem:vc_inverse_infty}
Let $\mscr{G}$ be a set of functions $\mscr{X}\rightarrow[0,\infty]$. If $\mscr{G}$ is VC-subgraph with VC-dimension at most $V$, then $\mscr{G}^{-1} := \left\{ \frac{1}{g}; g\in\mscr{G} \right\}$ is VC-subgraph with VC-dimension at most $V$, with the convention $1/0=+\infty$ and $1/+\infty=0$.
\end{lemme}
Combining this lemma with Proposition 42-$(ii)$ \cite{baraudinventiones}, we get that $\mscr{S}_{\alpha,2}$ is VC-subgraph with VC-dimension at most 196. This proves that we satisfy Assumption \ref{hyp:hmm_vc} with
\[ \overline{V}=4\times 196=784.\]
\subsubsection{Proof of Lemma \ref{lem:vc_set_to_function}}
Assume that $\mscr{F}_{\mscr{A}}$ has VC-dimension larger than $V$. Therefore, there is $(x_i,u_i)_{i\in[V+1]}\in \left(\mscr{X}\times \mbb{R}\right)^{[V+1]}$ such that for each $I\subset [V+1]$ we can find $A_I$ in $\mscr{A}$ such that $i\in I \Leftrightarrow f_{A_I}(x_i) > u_i$. Necessarily, we have $u_i\geq a$ for all $i\in [V+1]$ and therefore $i\in I \Leftrightarrow x_i \notin A_I$. Therefore, $\mscr{A}$ can shatter $(u_i)_{i_in[V+1]}$ which contradicts the fact that its VC-dimension is at most $V$.
\subsubsection{Proof of Lemma \ref{lem:vc_inverse_infty}}
We adapt the proof of Lemma 2.6.18 \cite{VanDerVaart}. Let $(x_i,u_i)_{i\in [n]}\in \left(\mscr{X}\times \mbb{R} \right)^n$ be such that for each $I\subset[n]$, we have $g_I\in\mscr{G}$ such that
\[ i\in I \Leftrightarrow \frac{1}{g_I(x_i)} > u_i. \]
For all $i\in[n]$, we necessarily have $u_i\geq 0$ and we define $a_i := \max\{ g_J(x_i); \frac{1}{g_J(x_i)} > u_i \}$. One can check that we have
\begin{align*}
g_I(x_i) > a_i \Leftrightarrow \frac{1}{g_I(x_i)} \leq u_i.
\end{align*}
Therefore $\mscr{G}$ shatters $\left(x_i,a_i\right)_{i\in[n]}\in\left(\mscr{X}\times\mbb{R}\right)^n$ which implies $n\leq V$.
\subsection{Proof of Proposition \ref{prop:faster_rates}}
\label{sec:proof_faster_rates_prop}
For $\pi=(\pi_{11},\pi_{12},\pi_{21},\pi_{22})\in\mcal{W}_4$ and $z\in\mbb{R}$ we write
\[ p_{\pi,z} := \pi_{11} f_{\alpha} \otimes f_{\alpha} + \pi_{12} f_{\alpha} \otimes f_{\alpha}(\cdot-z) + \pi_{21} f_{\alpha}(\cdot-z) \otimes f_{\alpha} + \pi_{22} f_{\alpha}(\cdot-z) \otimes f_{\alpha}(\cdot-z). \]
We define $\pi^*\in\mcal{W}_4$ by $\pi^*_{11}=w^*(1-q^*_{12})$, $\pi^*_{12}=w^*q^*_{12}$ and $\pi^*_{21}=(1-w^*)q^*_{21}$. We also define $g:\mcal{W}_4\times\mbb{R}\rightarrow \mbb{R}$ by
\[ g(\pi,z) = 2 h^2\left(P_{\pi^*,z^*},P_{\pi,z}\right) = \int_{\mbb{R}^2} a_{\pi,z}^2(x_1,x_2) dx, \]
with $a_{\pi,z}:\mbb{R}^2\rightarrow \mbb{R}$ defined by $a_{\pi,z}(x_1,x_2) = |\sqrt{p_{\pi,z}} - \sqrt{p_{\pi^*,z^*}}|$. We will drop the dependence on $\pi$ and $z$, and just write $a=a_{\pi,z}$. Without loss of generality we can assume $z^*>0$ as we have $h^2(P_{\pi,-z},P_{\pi^*,-z^*})=h^2(P_{\pi,z},P_{\pi^*,z^*})$. We define the set of parameters
\[ \mscr{Y} = \left\{ (\pi,z)\in\mcal{W}_4\times \mbb{R}; z\in\left( \frac{z^*}{2} \vee z^*-\beta^*, z^*+\beta^* \right) \right\}, \]
where $\beta^*\in(0,1]$ is set in the proof of Lemma \ref{lem:in_the_set} which proves the desired inequality on $\mscr{Y}$.
\begin{lemme}
\label{lem:in_the_set}
There is a positive constant $C(\alpha,z^*,\pi^*)$ such that
\[ g(\pi,z) \geq C(\alpha,z^*,\pi^*) \left[ \left(\pi^*_{11}-\pi_{11}\right)^2 + \left(\pi^*_{12}-\pi_{12}\right)^2 + \left(\pi^*_{21}-\pi_{21}\right)^2 + |z-z^*|^{1-\alpha} \right],\]
for all $(\pi,z)$ in $\mscr{Y}$. 
\end{lemme}
We also get that $g$ is lower bounded out of $\mscr{Y}$ with the following lemma.
\begin{lemme}
\label{lem:out_the set}
There is a positive constant $C(\alpha,z^*,\pi^*_{22})$ such that
\[ g(\pi,z) \geq C(\alpha,z^*,\pi^*_{22}), \forall (\pi,z)\notin \mscr{Y}. \]
\end{lemme}
One can check that we have $|z-z^*|^{1-\alpha}=\left(|z-z^*|\wedge 1\right)^{1-\alpha}$ for $(\pi,z)\in\mscr{Y}$. And since $\left(\pi^*_{11}-\pi_{11}\right)^2 + \left(\pi^*_{12}-\pi_{12}\right)^2 + \left(\pi^*_{21}-\pi_{21}\right)^2 + (|z-z^*|\wedge 1)^{1-\alpha}\leq 3$ for all $\pi$ and all $z$, there is a positive constant $C(\alpha,z^*,\pi^*)$ such that
\begin{align*}
g(\pi,z) &\geq C(\alpha,z^*,\pi^*) \left[ \left(\pi^*_{11}-\pi_{11}\right)^2 + \left(\pi^*_{12}-\pi_{12}\right)^2 + \left(\pi^*_{21}-\pi_{21}\right)^2 + (|z-z^*|\wedge 1)^{1-\alpha} \right], 
\end{align*}
for all $\pi,z$. We now relate the distance to $\pi^*$ to the distance to $(w^*,q^*)$ with the following result.
\begin{lemme}
\label{lem:pi_to_w_and_q}
For $w,q_{12},q_{21}\in[0,1]$ we have
\begin{align*}
&(\pi_{11}-\pi^*_{11})^2 + (\pi_{12}-\pi^*_{12})^2 + (\pi_{21} -\pi^*_{21})^2\\
&\geq \max\left( \frac{1}{2} (w-w^*)^2, \frac{(1-w^*)^2}{3} \left(q^*_{21}-q_{21}\right)^2, \left(w^*\right)^2 \left( q_{12} - q^*_{12} \right)^2 \right).
\end{align*}
\end{lemme}
This last result allows to conclude the proof of Proposition \ref{prop:faster_rates}.\par
\subsubsection{Proof of Lemma \ref{lem:in_the_set}}
We will repeatedly use the following inequality
\begin{equation}
\label{eq:taylor}
\forall x,y>0, \left|  x^{1-\gamma} - y^{1-\gamma} \right| \geq \frac{(1-\gamma)|x-y|}{ (x\vee y)^{\gamma} }.
\end{equation}
Let $(\pi,z)$ be in $\mscr{Y}$. Our goal is to lower bound $a$ on subsets of $\mscr{Y}$ in a way that makes appearing the difference between some parameters. Inequalities (\ref{eq:I11}), (\ref{eq:I22}), (\ref{eq:I12}) and (\ref{eq:I21}) will be proved later.
\begin{itemize}
 \item For $I_{11}=[-1,b)^2$ with $b=(z^*\wedge z\wedge 1)-1$, we have
\begin{equation}
\int_{I_{11}} a(x_1,x_2)^2 dx_1 dx_2 \geq \frac{ (1-\alpha) \left( 1\wedge  |z^*|/2 \right)^2 }{16} \left( \pi^*_{11} - \pi_{11} \right)^2.\label{eq:I11}
\end{equation}
\item For
\[ I_{22} = \left\{ \begin{array}{ll} (z^*, z^* +1) \times \left(z^*, z^* + (1-\alpha)^{2/\alpha} (\pi^*_{22})^{1/\alpha} |z-z^*| \right) &\text{if  } z^*\geq z, \\
\left( \frac{z^*}{2}\vee (z^*-1), z^* \right) \times \left(z^*, z^* + \frac{ (1-\alpha) (\pi^*_{22})^{1/\alpha} }{ (1-\alpha) \left(2 (\pi^*_{22})^{1/\alpha}+1\right) + 2 } |z-z^*| \right) &\text{otherwise},
\end{array}\right. \]
we have
\begin{equation}
\int_{I_{22}} a^2(x_1,x_2) dx \geq \frac{\alpha^2}{4^3} \left( \frac{3-\alpha}{2-\alpha} \right)^2 \left( \frac{1-\alpha}{5-3\alpha} \right)^{1-\alpha} (\pi^*_{22})^{1/\alpha} \left( 1 \wedge |z^*|/2 \right)^{1-\alpha}    |z-z^*|^{1-\alpha}. \label{eq:I22}
\end{equation}
\item Let $\beta\in(0,1]$. For
\begin{align*}
I_{12} &:= (-1,-(1-z\wedge z^*)_+)\times (z\vee z^*+b_-,z\vee z^*+b_+),\\
I_{21} &:= (z\vee z^*+b_-,z\vee z^*+b_+) \times (-1,-(1-z\wedge z^*)_+),
\end{align*}
with
\begin{align}                                        
b_+ &= \mathbbm{1}_{z\vee z^*\geq \beta}(1-|z-z^*|) + \mathbbm{1}_{z\vee z^*<\beta} \frac{z\vee z^* (1-\beta)}{\beta}\nonumber\\
&\geq \mathbbm{1}_{z^*\geq\beta}(1-\beta) + \mathbbm{1}_{z^*<\beta} \frac{z^* (1-\beta)}{\beta} = \left(1\wedge |z^*|/\beta\right) (1-\beta)\label{eq:b_plus}
\end{align}
and $b_-=b_+\delta$, $\delta\in (0,1)$. We have
\begin{align}
\int_{I_{12}} a^2(x_1,x_2) dx &\geq (\pi^*_{12}-\pi_{12})^2 \frac{ (1-\alpha)^2 \left(1\wedge |z^*|/2 \right) }{ 8^2 } \left(b_+\right)^{1-\alpha} (1-\delta) \mathbbm{1}_{\Omega_{12}},\label{eq:I12}\\
\int_{I_{21}} a^2(x_1,x_2) dx &\geq (\pi^*_{21}-\pi_{21})^2 \frac{ (1-\alpha)^2 \left(1\wedge |z^*|/2 \right) }{ 8^2 } \left(b_+\right)^{1-\alpha} (1-\delta) \mathbbm{1}_{\Omega_{21}},\label{eq:I21}
\end{align}
with
\begin{align*}
I_{12} &:= \left\{ |\pi^*_{12}-\pi_{12}| \geq 2 \left[ \frac{\alpha |z-z^*| }{ \delta b_+ } + |\pi_{11}-\pi^*_{11}| (1-\beta)^{\alpha} \right] \right\},\\
I_{21} &:= \left\{ |\pi^*_{21}-\pi_{21}| \geq 2 \left[ \frac{\alpha |z-z^*| }{ \delta b_+ } + |\pi_{11}-\pi^*_{11}| (1-\beta)^{\alpha} \right] \right\}.
\end{align*}
\end{itemize}
Combining (\ref{eq:I11}), (\ref{eq:I22}), (\ref{eq:I12}) and (\ref{eq:I21}), we have
\begin{align*}
\int a^2(x_1,x_2) dx &\geq \left( \pi^*_{11} - \pi_{11} \right)^2  \frac{ (1-\alpha)^2 \left( 1\wedge |z^*|/2 \right)^2 }{16} \\
&+ |z-z^*|^{1-\alpha} \frac{\alpha^2}{4^3} \left( \frac{3-\alpha}{2-\alpha} \right)^2 \left( \frac{1-\alpha}{5-3\alpha} \right)^{1-\alpha} (\pi^*_{22})^{1/\alpha} \left( 1 \wedge |z^*|/2 \right)^{1-\alpha}\\
&+ (\pi^*_{12}-\pi_{12})^2 \frac{ (1-\alpha)^2 \left(1 \wedge |z^*|/2 \right) }{ 8^2 } \left(b_+\right)^{1-\alpha} (1-\delta) \mathbbm{1}_{\Omega_{12}}\\
&+ (\pi^*_{21}-\pi_{21})^2 \frac{ (1-\alpha) \left(1\wedge |z^*|/2 \right) }{ 8^2 } \left(b_+\right)^{1-\alpha} (1-\delta) \mathbbm{1}_{\Omega_{21}},
\end{align*}
for $(\pi,z)\in\mscr{Y}$. Then we can apply the following lemma.
\begin{lemme}
\label{lem:racor}
Let $g,A_1,A_2,A_3,B$ be functions $\Theta \rightarrow \mbb{R}$ and $D_1,D_{2,3},D_B,C_A,C_B$ be positive constants such that 
\[ \forall \theta\in\Theta, g(\theta) \geq D_1 A_1^2(\theta) + D_{2,3} \left( A_2^2(\theta) \mathbbm{1}_{\Omega_2} + A_3^2(\theta) \mathbbm{1}_{\Omega_3} \right) + D_B(\theta) B^{1-\alpha},\]
where $\Omega_2$ and $\Omega_3$ are subsets of $\Theta$ given by
\[ \Omega_i := \left\{ \theta\in\Theta; A_i(\theta) \geq C_A A_1(\theta) + C_B B(\theta) \right\}. \]
Then we have
\[ g(\theta) \geq \min\left( \frac{D_B}{1+4C_B^2}, \frac{D_1}{1+4C_A^2}, D_{2,3} \right) \left[ A_1^2(\theta) + A_2^2(\theta) + A_3^2(\theta) + B^{1-\alpha}(\theta) \right], \]
for all $\theta$ in $\Theta$. 
\end{lemme}
In our situation, we get
\begin{align*}
\int a^2(x_1,x_2) dx \geq C(\alpha,z^*,\pi^*) \left[ \left(\pi^*_{11}-\pi_{11}\right)^2 + \left(\pi^*_{12}-\pi_{12}\right)^2 + \left(\pi^*_{21}-\pi_{21}\right)^2 + |z-z^*|^{1-\alpha} \right]
\end{align*}
with
\begin{align*}
 C(\alpha,z^*,\pi^*) = \min &\left( \frac{ \frac{\alpha^2}{4^3} \left(\frac{3-\alpha}{2-\alpha}\right)^2 \left( \frac{1-\alpha}{5-3\alpha} \right)^{1-\alpha} \left( \pi^*_{22} \right)^{1/\alpha} (1\wedge |z^*|/2)^{1-\alpha} }{ 1 + 4^2 \frac{\alpha^2}{\delta^2 b_+^2} }, \frac{ \frac{(1-\alpha)^2 (1\wedge |z^*|/2)^2}{4^2} }{ 1+4 (1-\beta)^{2\alpha} },\right.\\
 &\left. \frac{(1-\alpha)(1\wedge|z^*|/2) }{8^2} \left(b_+\right)^{1-\alpha} (1-\delta) \right)\\
 \geq \min &\left( \frac{ \frac{\alpha^2}{4^3} \left(\frac{3-\alpha}{2-\alpha}\right)^2 \left( \frac{1-\alpha}{5-3\alpha} \right)^{1-\alpha} \left( \pi^*_{22} \right)^{1/\alpha} (1\wedge |z^*|/2)^{1-\alpha} }{ 1 + 4^2 \frac{\alpha^2}{\delta^2 (1\wedge |z^*|/2)^2 (1-\beta)^2} }, \frac{ (1-\alpha)^2 (1\wedge |z^*|/2)^2}{ 4^2 \left( 1+4 (1-\beta)^{2\alpha} \right) }, \right.\\
 &\left. \frac{(1-\alpha)(1\wedge|z^*|/2) }{8^2} (1\wedge |z^*|/2)^{1-\alpha} (1-\beta)^{1-\alpha} (1-\delta) \right)>0.
\end{align*}
We can optimize this bound with respect to $\beta$ and $\delta$, which gives $\beta^*$ depending only $z^*$, $\alpha$ and $\pi^*$. This concludes the proof of Lemma \ref{lem:in_the_set}. We now prove the different inequalities.\\
\emph{Proof of (\ref{eq:I11}).} For $x_1,x_2\in [-1,0)^2$, we have
\begin{align*}
&a(x_1,x_2) = \frac{1-\alpha}{2 |x_1|^{\alpha/2}|x_2|^{\alpha/2}}\\
&\times  \left| \sqrt{ \pi_{11}^* + \pi_{12}^* \frac{\mathbbm{1}_{|x_2-z^*|\in(0,1]} |x_2|^{\alpha}}{|x_2-z^*|^{\alpha}} + \pi_{22}^* \frac{ \mathbbm{1}_{|x_1-z^*|\in(0,1]} \mathbbm{1}_{|x_2-z^*|\in(0,1]} |x_1|^{\alpha} |x_2|^{\alpha} }{ |x_1-z^*|^{\alpha} |x_2-z^*|^{\alpha} } + \pi_{21}^* \frac{ \mathbbm{1}_{|x_1-z^*|\in(0,1]} |x_1|^{\alpha} }{ |x_1-z^*|^{\alpha} } } \right.\\
&- \left. \sqrt{ \pi_{11} + \pi_{12} \frac{ \mathbbm{1}_{|x_2-z|\in(0,1]} |x_2|^{\alpha} }{ |x_2-z|^{\alpha} } + \pi_{22} \frac{ \mathbbm{1}_{|x_1-z|\in(0,1]} \mathbbm{1}_{|x_2-z|\in(0,1]} |x_1|^{\alpha} |x_2|^{\alpha} }{ |x_1-z|^{\alpha} |x_2-z|^{\alpha} } + \pi_{21} \frac{ \mathbbm{1}_{|x_1-z|\in(0,1]} |x_1|^{\alpha} }{ |x_1-z|^{\alpha} } } \right|.
\end{align*}
We set $b=\min(z^*, z, 1)-1$. For $x_1,x_2\in[-1,b)^2$, we have
\[ a(x_1,x_2) = \frac{1-\alpha}{2|x_1|^{\alpha/2}|x_2|^{\alpha/2}} \left| \sqrt{ \pi^*_{11} } - \sqrt{ \pi_{11} } \right| \]
and
\[ \int_{[-1,b)^2} a(x_1,x_2)^2 dx_1 dx_2 \geq \frac{ \left[ 1-(-)^{1-\alpha} \right]^2 }{4} \left| \sqrt{ \pi^*_{11} } - \sqrt{ \pi_{11} } \right|^2. \]
Finally, with (\ref{eq:taylor}) we always have
\begin{align*}
\int_{[-1,b)^2} a(x_1,x_2)^2 dx_1 dx_2 &\geq \frac{ \left[ 1-(1-z\wedge z^*)_+^{1-\alpha} \right]^2 }{4} \left( \sqrt{ \pi^*_{11} } - \sqrt{ \pi_{11} } \right)^2\\
&\geq \frac{ (1-\alpha) \left( 1\wedge |z^*|/2 \right)^2 }{4^2} \left( \pi^*_{11} - \pi_{11} \right)^2.
\end{align*}
\emph{Proof of (\ref{eq:I22}).} We need to consider two different cases.
\begin{itemize}
 \item[$\bullet$] \emph{First case $z^*\geq z$.} For $x \in I_{22} = (z^*,z^*+1) \times (z^*, z^* + V|z-z^*|)$ with $V<\frac{1}{|z-z^*|}$, we have $
\frac{|x_2-z^*|}{|x_2-z|} \leq V$, $\frac{|x_2-z^*|}{|x_2|} \leq V$, $\frac{|x_1-z^*|}{|x_1-z|} \leq 1 - |z-z^*| \leq 1$ and $\frac{|x_1-z^*|}{|x_1|} \leq \frac{1}{1+z^*} \leq 1$. Therefore, for $x\in I_{22}$, we have
\[ a(x_1,x_2) \geq \frac{1-\alpha}{2 |x_1-z^*|^{\alpha/2} |x_2-z^*|^{\alpha/2}} \left( \sqrt{ \pi^*_{22} } - V^{\alpha/2} \right)_+. \]
For $V=(1-\alpha)^{2/\alpha} (\pi^*_{22})^{1/\alpha}<\frac{1}{|z-z^*|}$, we have
\begin{align*}
\int_{I_{22}} a^2(x_1,x_2) dx &= \frac{\left( \sqrt{ \pi^*_{22} } - V^{\alpha/2} \right)_+^2}{4} \left( V |z-z^*| \wedge 1 \right)^{1-\alpha} \\
&\geq \frac{ (\pi_{22}^*)^{1/\alpha} \alpha^2  (1-\alpha)^{2(1-\alpha)/\alpha} |z-z^*|^{1-\alpha} }{ 4 }.
\end{align*}
\item[$\bullet$] \emph{Second case $z^*<z$.} For $x\in \left(\frac{z^*}{2}\vee (z^*-1),z^*\right)\times (z^*,z^*+a|z-z^*|)$, $b\leq 1/2$ we have
\[ a(x_1,x_2) \geq \frac{1-\alpha}{2|x_1-z^*|^{\alpha}|x_2-z^*|^{\alpha}} \left( \sqrt{\pi^*_{22}} - \left(\frac{b}{1-b} \right)^{\alpha/2} \right)_+. \]
For $b=\left(\pi^*_{22}\right)^{1/\alpha} b'$ we have
\begin{align*}
&\int_{I_{22}} a^2(x_1,x_2) dx \geq \frac{\left( \sqrt{ \pi^*_{22} } - \left(\frac{b}{1-b}\right)^{\alpha/2} \right)_+^2}{4} \left( 1\wedge |z^*|/2 \right)^{1-\alpha} b^{1-\alpha} |z-z^*|^{1-\alpha}\\
&\geq \frac{ \left( \pi^*_{22} \right)^{1/\alpha} \left( 1 \wedge |z^*|/2 \right)^{1-\alpha} |z-z^*|^{1-\alpha} }{4} \left( 1 - \left( \frac{b'}{1-b'(\pi^*_{22})^{1/\alpha}}\right)^{\alpha/2} \right)_+^2  (b')^{1-\alpha}\\
&\geq \frac{ \left( \pi^*_{22} \right)^{1/\alpha} \left( 1 \wedge |z^*|/2 \right)^{1-\alpha} |z-z^*|^{1-\alpha} }{4} \frac{\alpha^2}{4}  \left( \frac{1-b'\left(1+(\pi^*_{22})^{1/\alpha}\right)}{1-b'(\pi^*_{22})^{1/\alpha}}\right)^2  (b')^{1-\alpha}.
\end{align*}
With $b'=\frac{1-\alpha}{(1-\alpha)(2\pi+1)+2}$ we have
\begin{align*}
\int_{I_{22}} a^2(x_1,x_2) dx &\geq \frac{ \alpha^2 \left( \pi^*_{22} \right)^{1/\alpha} \left( 1 \wedge |z^*|/2 \right)^{1-\alpha} |z-z^*|^{1-\alpha} }{4^2}\\
&\times \left( \frac{2 + (1-\alpha) (\pi^*_{22})^{1/\alpha} }{ 2 + (1-\alpha) \left(1+(\pi^*_{22})^{1/\alpha} \right)}\right)^2  \left(\frac{1-\alpha}{(1-\alpha)(2\pi+1)+2}\right)^{1-\alpha}\\
&\geq \frac{ \alpha^2 \left( \pi^*_{22} \right)^{1/\alpha} \left( 1 \wedge |z^*|/2 \right)^{1-\alpha} |z-z^*|^{1-\alpha} }{4^2}\\
&\times \left( \frac{3-\alpha }{ 2 + 2 (1-\alpha) }\right)^2  \left(\frac{1-\alpha}{5-3\alpha}\right)^{1-\alpha}\\
&= \frac{ \alpha^2 (3-\alpha)^2 }{4^3 \left( 2-\alpha \right)^2 }  \left(\frac{1-\alpha}{5-3\alpha}\right)^{1-\alpha} \left( \pi^*_{22} \right)^{1/\alpha} \left( 1 \wedge |z^*|/2 \right)^{1-\alpha} |z-z^*|^{1-\alpha}.
\end{align*}
\end{itemize}
Finally, we always have have
\[ \int_{I_{22}} a^2(x_1,x_2) dx \geq \frac{\alpha^2}{4^3} \left( \frac{3-\alpha}{2-\alpha} \right)^2 \left( \frac{1-\alpha}{5-3\alpha} \right)^{1-\alpha} (\pi^*_{22})^{1/\alpha} \left( 1 \wedge |z^*|/2 \right)^{1-\alpha}    |z-z^*|^{1-\alpha}.\]
\emph{Proof of (\ref{eq:I12}).} We prove it for $I_{12}$ assuming $z^*\leq z$. The proof is similar for $I_{21}$ and for $z\leq z^*$. For $b=0\wedge(z^*-1)$ and $0<c_-<c_+<1-|z-z^*|$, we set $I_{12} = (-1,b) \times (z + c_-,z^*+1)$. For $x_1,x_2\in I_{12}$, we have
\[ \frac{2 |x_1|^{\alpha/2} |x_2-z|^{\alpha/2}}{1-\alpha} a(x_1,x_2) = \left| \frac{ (\pi^*_{12}-\pi_{12}) + \pi^*_{12} \left( \frac{|x_2-z|^{\alpha}}{|x_2-z^*|^{\alpha}} -1 \right) + (\pi^*_{11}-\pi_{11}) \frac{|x_2-z|^{\alpha}\mathbbm{1}_{x_2\leq 1}}{|x_2|^{\alpha}} }{ \sqrt{ \pi^*_{12} \frac{|x_2-z|^{\alpha}}{|x_2-z^*|^{\alpha}} + \pi^*_{11} \frac{|x_2-z|^{\alpha}\mathbbm{1}_{x_2\leq 1}}{|x_2|^{\alpha}} } + \sqrt{ \pi^*_{12} + \pi_{11} \frac{|x_2-z|^{\alpha}\mathbbm{1}_{x_2\leq 1}}{|x_2|^{\alpha}} } } \right|.\]
We also have
\[ \mathbbm{1}_{x_2\leq 1} \frac{|x_2-z|}{|x_2|} \leq U(z,c_-,c_+) := \left\{ \begin{array}{ll}  \frac{c_+}{z+c_+} &\text{if } z+c_+\leq 1,\\
      1-z &\text{if } z+c_-< 1 < z+c_+,\\
      0 &\text{if } 1\leq z+c_-.
\end{array} \right.
\]
For $c_+=\mathbbm{1}_{z\geq \beta^*}(1-|z-z^*|) + \mathbbm{1}_{z<\beta^*} \frac{z(1-\beta^*)}{\beta^*}$ we have $U(z,c_-,c_+)\leq 1-\beta^*$. We also have
\begin{align*}
1-\frac{|x_2-z|^{\alpha}}{|x_2-z^*|^{\alpha}} &\leq 1 - \left( \frac{ c_- }{ c_-+|z-z^*| } \right)^{\alpha}\\
&\leq \alpha \frac{ \frac{ |z-z^*| }{ c_-+|z-z^*| } }{ \left( \frac{ c_- }{ c_-+|z-z^*| } \right)^{1-\alpha} } = \alpha \frac{|z-z^*|}{c_-} \left( \frac{c_-}{c_-+|z-z^*|} \right)^{\alpha}\\
&\leq  \frac{\alpha |z-z^*|}{c_-}.
\end{align*}
Therefore, with $c_-=c_+\delta$, $\delta\in(0,1)$, on $I_{12}$ we have
\[ \frac{2 |x_1|^{\alpha/2} |x_2-z|^{\alpha/2}}{1-\alpha} a(x_1,x_2) \geq \frac{ \left[ |\pi^*_{12}-\pi_{12}| - \frac{\alpha |z^*-z|}{b_-} - \left(\pi^*_{11}-\pi_{11}\right) (1-\beta^*)^{\alpha} \right]_+ }{2}.\]
If $|\pi^*_{12}-\pi_{12}| \geq 2 \left[ |\pi^*_{11}-\pi_{11}| (1-\beta^*)^{\alpha} + \alpha |z-z^*|/c_- \right]$ then
\[ \frac{2 |x_1|^{\alpha/2} |x_2-z|^{\alpha/2}}{1-\alpha} a(x_1,x_2) \geq \frac{ |\pi^*_{12}-\pi_{12}| }{4} \]
and
\begin{align*}
\int_{I_{12}} a^2(x_1,x_2) dx &\geq \frac{(\pi^*_{12}-\pi_{12})^2}{8^2} \left[ 1 - (1-z\wedge z^*)_+^{1-\alpha} \right] \left(c_+\right)^{1-\alpha} \left[ 1 - \delta^{1-\alpha}  \right]\\
&\geq (\pi^*_{12}-\pi_{12})^2 \frac{ (1-\alpha)^2 \left(1\wedge |z| \wedge |z^*| \right) \left(c_+\right)^{1-\alpha} \left( 1-\delta \right) }{ 8^2 }.
\end{align*}
Otherwise we have $|\pi^*_{12}-\pi_{12}| < 2 \left[ |\pi^*_{11}-\pi_{11}| (1-\beta^*)^{\alpha} + \alpha |z-z^*|/b_- \right]$.
\subsubsection{Proof of Lemma \ref{lem:out_the set}}
We need to go through numerous cases and subcases. Let $\beta^*$ be given in Lemma \ref{lem:in_the_set}. Without loss of generality we are going to assume that $z^*>0$.\\
\emph{Case 1}: $z\geq 0$ and $|z-z^*|\geq \beta^*$. Let $c$ be a positive constant.
\begin{itemize}
\item[$\bullet$] \emph{Subcase 1.1}: $z^*>z$ or ($z^*<z$ and $\pi_{22}\geq c^2 \pi^*_{22}$). For $x\in I = (z\vee z^*+\beta^*,z\vee z^*+1)^2$, we have
\[ a(x_1,x_2) = \frac{(1-\alpha) \left(\mathbbm{1}_{z>z^*} \pi_{22} + \mathbbm{1}_{z^*>z} \pi^*_{22}\right)}{2 |x_1-z\vee z^*|^{\alpha/2} |x_2-z\vee z^*|^{\alpha/2}},\]
and therefore
\begin{align*}
\int_I a^2(x_1,x_2) dx &= \frac{\mathbbm{1}_{z>z^*} \pi_{22} + \mathbbm{1}_{z^*>z} \pi^*_{22}}{4} \left( 1- (\beta^*)^{1-\alpha} \right)^2\\
&\geq \frac{c^2 \pi_{22}^* (1-\alpha)^2}{4} \left( 1-\beta^*\right)^2.
\end{align*}
\item[$\bullet$] \emph{Subcase 1.2}: $1\leq z^*<z$ and $\pi_{22}< c^2 \pi^*_{22}$. For $x\in \left(z^*,z^*+1\wedge (|z-z^*|/2) \right)^2$, we have
\[ \frac{|x_1-z^*|}{|x_1-z|} \leq \frac{1 \wedge |z-z^*|/2 }{ z-z^*-1\wedge |z-z^*|/2 } \leq 1. \]
We have
\begin{align*}
a(x_1,x_2) &\geq \frac{1-\alpha}{2 |x_1-z^*|^{\alpha/2} |x_2-z^*|^{\alpha/2}} \left( \sqrt{\pi^*_{22}} - \sqrt{\pi_{22}} \right)\\
&\geq \frac{(1-\alpha) \sqrt{\pi^*_{22}} }{2 |x_1-z^*|^{\alpha/2} |x_2-z^*|^{\alpha/2}} \left( 1 - c \right),
\end{align*}
and therefore
\begin{align*}
\int_I a^2(x_1,x_2)dx &\geq \frac{\pi^*_{22}}{4} \left( 1 - c \right)^2 \left( 1\wedge \frac{|z-z^*|}{2} \right)^{2(1-\alpha)}\\
&\geq \frac{\pi^*_{22} \left( 1 - c \right)^2}{2^{2(2-\alpha)}} \left( \beta^* \right)^{2(1-\alpha)}. 
\end{align*}
\item[$\bullet$] \emph{Subcase 1.3}: $z^*\in(0,1-\beta^*]$ and $z^*<z$. Let $b$ be in $(0,1)$. For $x\in I=\left(z^* - b z^*,z^* \right)^2$ we have
\begin{align*}
\frac{|x_1-z^*|}{|x_1|} &\leq \frac{ bz^* }{z^*-bz^*}=\frac{b}{1-b}\\
\frac{|x_1-z^*|}{|x_1-z|} &\leq \frac{ bz^* }{z-z^*+bz^*}\leq \frac{b\beta^*}{\beta^*+b\beta^*} \leq \frac{b}{1-b}. 
\end{align*}
It implies
\begin{align*}
a(x_1,x_2) &\geq \frac{1-\alpha}{2 |x_1-z^*|^{\alpha/2} |x_2-z^*|^{\alpha/2}} \left( \sqrt{\pi^*_{22}} - \left( \frac{b}{1-b} \right)^{\alpha} \right)_+,
\end{align*}
and for $b=b' \left(\pi^*_{22}\right)^{1/2\alpha}$ we get
\begin{align*}
\int_I a^2(x_1,x_2)dx &\geq \frac{ (z^*)^{2(1-\alpha)} \left(\pi^*_{22}\right)^{(1-\alpha)/\alpha} (b')^{2(1-\alpha)}}{4} \left( \sqrt{\pi^*_{22}} - \sqrt{\pi^*_{22}} \left( \frac{b'}{1-(\pi^*_{22})^{1/2\alpha} b'} \right)^{\alpha} \right)_+^2\\
&\geq \frac{ (z^*)^{2(1-\alpha)} \left(\pi^*_{22}\right)^{1/\alpha} (b')^{2(1-\alpha)}}{4} \alpha^2 \left( 1 - \frac{b'}{1-(\pi^*_{22})^{1/2\alpha} b'}  \right)_+^2.
\end{align*}
For $b' = \frac{1}{ 1 + 2 \left(\pi^*_{22}\right)^{1/2\alpha} + \frac{1}{1-\alpha} }$, we have
\begin{align*}
\int_I a^2(x_1,x_2)dx &\geq \frac{ (z^*)^{2(1-\alpha)} \left(\pi^*_{22}\right)^{1/\alpha} \alpha^2 }{ 4 \left( 1 + 2 \left(\pi^*_{22}\right)^{1/2\alpha} + \frac{1}{1-\alpha} \right)^{2(1-\alpha)} } \left( 1 - \frac{ 1 }{ 1 + (\pi^*_{22})^{1/2\alpha} + \frac{1}{1-\alpha} } \right)^2.
\end{align*}
\item[$\bullet$] \emph{Subcase 1.4}: $z^*<z$ and $z^*\in[1-\beta^*,1]$. Let $b$ be in $(0,1)$. For $x\in I=\left(z^*,z^* + b \beta^* \right)^2$ we have
\begin{align*}
\frac{|x_1-z^*|}{|x_1-z|} &\leq \frac{b \beta^*}{z-z^*-b\beta^*} \leq \frac{b}{1-b},\\
\frac{|x_1-z^*|}{|x_1|} &\leq \frac{ b \beta^* }{z^*+b\beta^*} \leq \frac{b}{1+b} \leq \frac{b}{1-b}.
\end{align*}
It implies
\[ a(x_1,x_2) \geq \frac{1-\alpha}{2 |x_1-z^*|^{\alpha/2} |x_2-z^*|^{\alpha/2}} \left( \sqrt{\pi^*_{22}} - \left( \frac{b}{1-b} \right)^{\alpha} \right)_+,\]
and for $b=b'\left(\pi^*_{22}\right)^{1/2\alpha}$ we get
\begin{align*}
\int_I a^2(x_1,x_2)dx &\geq \frac{(\beta^*)^{2(1-\alpha)} \left( \pi^*_{22} \right)^{(1-\alpha)/\alpha} (b')^{2(1-\alpha)} }{4} \pi^*_{22} \left( 1 - \left( \frac{b'}{1-b' (\pi^*_{22})^{1/2\alpha}} \right)^{\alpha} \right)^2_+\\
&\geq \frac{ (\beta^*)^{2(1-\alpha)} \left( \pi^*_{22} \right)^{1/\alpha} \alpha^2 }{4} (b')^{2(1-\alpha)} \left( 1 - \frac{b'}{1-b'(\pi^*_{22})^{1/2\alpha}} \right)^2_+.
\end{align*}
For $b' = \frac{1}{ 1 + 2 \left(\pi^*_{22}\right)^{1/2\alpha} + \frac{1}{1-\alpha} }$ we have
\begin{align*}
\int_I a^2(x_1,x_2)dx &\geq \frac{ (\beta^*)^{2(1-\alpha)} \left(\pi^*_{22}\right)^{1/\alpha} \alpha^2 }{ 4 \left( 1 + 2 \left(\pi^*_{22}\right)^{1/2\alpha} + \frac{1}{1-\alpha} \right)^{2(1-\alpha)} } \left( 1 - \frac{ 1 }{1 + (\pi^*_{22})^{1/2\alpha} + \frac{1}{1-\alpha}} \right)^2.
\end{align*}
\end{itemize}
We can optimize the subcases 1.1 and 1.2 with $c=\frac{ \left( \beta^*/2\right)^{2(1-\alpha)} }{ \left( \beta^*/2 \right)^{2(1-\alpha)} + (1-\alpha)(1-\beta^*) }$. Gathering the different results, there is a positive constant $C_1(z^*,\pi^*_{22},\alpha)$ such that $\int_{\mbb{R}^2} a(x_1,x_2) dx \geq C_1(\pi^*_{22},z^*,\alpha)$ for all $z$ satisfying $z\geq 0$ and $|z-z^*|\geq 1-\beta^*$.\\
\emph{Case 2}: $z<0$.
\begin{itemize}
 \item[$\bullet$] \emph{Subcase 2.1}: $z^*\leq 1$. Let $b$ be in $(0,1)$. For $x\in(z^*,z^*+b)^2$ we have $\frac{|x_1-z^*|}{|x-z|} \leq \frac{|x_1-z^*|}{|x_1|} \leq \frac{b}{z^*}$ and therefore
\[ a(x_1,x_2) \geq \frac{1-\alpha}{2 |x_1-z^*|^{\alpha/2} |x_2-z^*|^{\alpha/2}} \left[ \sqrt{\pi^*_{22}} - \left(\frac{b}{z^*}\right)^{\alpha} \right]_+.\]
We get $\int_{(z^*,z^*+b)^2} a^2(x_1,x_2) dx \geq \frac{ b^{2(1-\alpha)} \left[ \sqrt{\pi^*_{22}}-\left(\frac{b}{z^*}\right)^{\alpha} \right]^2_+ }{4}$. For $b=z^* (\pi^*_{22})^{1/2\alpha} (1-\alpha)^{1/\alpha}\leq 1$, we have
\[ \int_{(z^*,z^*+b)^2} a^2(x_1,x_2)dx \geq \frac{ \alpha^2 (1-\alpha)^{2(1-\alpha)/\alpha} (z^*)^{2(1-\alpha)} (\pi^*_{22})^{1/\alpha} }{4}.\]
\item[$\bullet$] \emph{Subcase 2.2}: $z^*>1$. For $x\in(z^*,z^*+1)^2$ we have
\[ a(x_1,x_2) = \frac{1-\alpha}{2} \sqrt{ \frac{\pi^*_{22}}{ |x_1-z^*|^{\alpha} |x_2-z^*|^{\alpha} } } .\]
Therefore we get $\int_{(z^*,z^*+1)^2} a^2(x_1,x_2)dx \geq \frac{\pi^*_{22}}{4}$.
\end{itemize}
Finally, we have
\[ \int_{\mbb{R}^2} a^2(x_1,x_2) dx \geq \frac{ \alpha^2 (1-\alpha)^{2(1-\alpha)/\alpha} (1\wedge z^*)^{2(1-\alpha)} (\pi^*_{22})^{1/\alpha} }{4}. \]
\emph{Case 3}: $|z-z^*|<\beta^*$ and $z\leq z^*/2$. Let $b$ be in $(0,1/|z-z^*|)$. For $x \in (z^*,z^*+b|z-z^*|)^2$ we have
\begin{align*}
\frac{|x_1-z^*|}{|x_1|} &\leq \frac{b|z-z^*|}{z^*+b|z-z^*|}\leq b\\ 
\frac{|x_1-z^*|}{|x_1-z|} &\leq \frac{b|z-z^*|}{b|z-z^*|+|z^*-z|} \leq b.
\end{align*}
Therefore we get
\[ a(x_1,x_2) \geq \frac{1-\alpha}{2 |x_1-z^*|^{\alpha/2} |x_2-z^*|^{\alpha/2}} \left[ \sqrt{\pi^*_{22}} - b^{\alpha} \right]_+.\]
We get
\[ \int_{(z^*,z^*+b|z-z^*|)^2} a^2(x_1,x_2)dx \geq \frac{b^{2(1-\alpha)} |z-z^*|^{2(1-\alpha)} \left[ \sqrt{\pi^*_{22}} - b^{\alpha} \right]^2_+ }{4}\]
and for $b = (\pi^*_{22})^{1/2\alpha} (1-\alpha)^{1/\alpha}\leq 1/|z-z^*|$ we have
\begin{align*}
\int_{(z^*,z^*+b|z-z^*|)^2} a^2(x_1,x_2)dx &\geq \frac{ |z-z^*|^{2(1-\alpha)} (1-\alpha)^{2(1-\alpha)/\alpha} \left(\pi_{22}^*\right)^{(1-\alpha)/\alpha} }{4}  \pi^*_{22} \alpha^2\\
&\geq \frac{ \alpha^2 \left( |z^*|/2 \right)^{2(1-\alpha)} (1-\alpha)^{2(1-\alpha)/\alpha} \left(\pi_{22}^*\right)^{1)/\alpha} }{4}.
\end{align*}
\subsubsection{Proof of Lemma \ref{lem:racor}.}
\begin{itemize}
\item[$\bullet$] For $\theta$ in $\Omega_2\cap\Omega_3$, we have
\begin{align*}
g(\theta) &\geq D_1 A_1^2 + D_{2,3} \left( A_2^2 + A_3^2 \right) + D_B B^{1-\alpha}\\
&\geq \min\left( D_1, D_{2,3}, D_B \right) \left[ A_1^2 + A_2^2 + A_3^2 + B^{1-\alpha} \right]. 
\end{align*}
\item[$\bullet$] For $\theta$ in $\Omega_2\cap\Omega_3^C$, we have
\[ g(\theta) \geq D_1 A_1^2 + D_{2,3} A_2^2 + D_B B^{1-\alpha}\]
and
\[ A_3^2 < \left(C_A A_1 + C_B B\right)^2 \leq 2 C_A^2 A_1^2 + 2 C_B^2 B^{1-\alpha}.\]
For $b=\frac{D_B}{1+2C_B^2}\bigwedge \frac{D_1}{1+2C_A^2}>0$ we have
\begin{align*}
g(\theta) &\geq D_{2,3} A_2^2 + \left( D_1 - b 2 C_A^2 \right) A_1^2 + D_{2,3} A_2^2 + (D_B-b2C_B^2) B^{1-\alpha} + b A_3^2\\
&\geq \min\left( \frac{D_B}{1+2C_B^2}, \frac{D_1}{1+2C_A^2}, D_{2,3}\right) \left[ A_1^2 + A_2^2 + A_3^2 + B^{1-\alpha} \right].
\end{align*}
\item[$\bullet$] For $\theta$ in $\Omega_2^C\cap\Omega_3^C$, we have
\[ g(\theta) \geq D_1 A_1^2 + D_B B^{1-\alpha}\]
and
\[ A_2^2 + A_3^2 < 2 \left(C_A A_1 + C_B B\right)^2 \leq 4 C_A^2 A_1^2 + 4 C_B^2 B^{1-\alpha}.\]
For $b=\frac{D_B}{1+4C_B^2}\bigwedge \frac{D_1}{1+4C_A^2}>0$ we have
\begin{align*}
g(\theta) &\geq D_{2,3} A_2^2 + \left( D_1 - b 4 C_A^2 \right) A_1^2 + (D_B-b 4 C_B^2) B^{1-\alpha} + b \left( A_2^2 + A_3^2 \right)\\
&\geq \min\left( \frac{D_B}{1+4C_B^2}, \frac{D_1}{1+4C_A^2} \right) \left[ A_1^2 + A_2^2 + A_3^2 + B^{1-\alpha} \right].
\end{align*}
\end{itemize}
Finally, we always have
\[ g(\theta) \geq \min\left( \frac{D_B}{1+4C_B^2}, \frac{D_1}{1+4C_A^2}, D_{2,3} \right) \left[ A_1^2 + A_2^2 + A_3^2 + B^{1-\alpha} \right]. \]
\subsubsection{Proof of Lemma \ref{lem:pi_to_w_and_q}}
We assume there is $w,w^*,q_{12},q_{12}^2,q_{21},q_{21}^*$ in $[0,1]$ such that
\[ \pi_{11}=w(1-q_{12}), \pi_{12}=w q_{12}, \pi_{21}=(1-w) q_{21}\]
and
\[ \pi_{11}^*=w^*(1-q^*_{12}), \pi^*_{12}=w^* q^*_{12}, \pi^*_{21}=(1-w^*) q^*_{21}.\]
\begin{itemize}
\item We have 
\begin{align}
(\pi_{11}-\pi^*_{11})^2 + (\pi_{12}-\pi^*_{12})^2&= w^2 \left[ 2 \left( q_{12} -\frac{1}{2} \right)^2 + \frac{1}{2} \right]\nonumber\\
&- 2ww^* \left[ 2\left( q^*_{12}-\frac{1}{2}\right) \left(q_{12}-\frac{1}{2} \right) + \frac{1}{2} \right]\nonumber\\
&+ (w^*)^2 \left[ 2 \left(q_{12}^*-\frac{1}{2}\right)^2 + \frac{1}{2} \right]\nonumber\\
&= \frac{1}{2} (w-w^*)^2 + 2 \left( w\left(q_{12}-\frac{1}{2}\right) - w^*\left(q^*_{12}-\frac{1}{2}\right) \right)^2\nonumber\\
&\geq \frac{1}{2} (w-w^*)^2.\label{eq:low_w}
\end{align}
Therefore, we also have
\begin{align}
&(\pi_{11}-\pi^*_{11})^2 + (\pi_{12}-\pi^*_{12})^2 + (\pi_{21} -\pi^*_{21})^2\nonumber\\
&\geq \frac{1}{2} (w-w^*)^2 + \left( (1-w) q_{21} - (1-w^*) q^*_{21} \right)^2\nonumber\\
&= (1-w)^2 \left[ \frac{1}{2} + q_{21}^2 \right] + (1-w^*)^2 \left[ \frac{1}{2} + (q^*_{21})^2 \right]\nonumber\\
&- (1-w) (1-w^*) \left[ 1 + 2q_{21} q^*_{21} \right]\nonumber\\
&= \left[ \frac{1}{2} + q_{21}^2 \right] \left( (1-w) - (1-w^*) \frac{1+2q_{21}q^*_{21}}{1+2q_{21}^2} \right)^2\nonumber\\
&+ (1-w^*)^2 \left[ \frac{1}{2} + (q^*_{21})^2 \right] - \left[ \frac{1}{2} + q_{21}^2 \right] (1-w^*)^2 \left( \frac{1+2q_{21}q^*_{21}}{1+2q^2_{21}} \right)^2\nonumber\\
&\geq \frac{(1-w^*)^2}{2\left( 1+2q_{21}^2 \right) } \left[ (1+2(q^*_{21})^2)(1+2q_{21}^2) - (1+2q_{21}q^*_{21})^2 \right]\nonumber\\
&= \frac{(1-w^*)^2}{1+2q_{21}^2} \left(q^*_{21}-q_{21}\right)^2\nonumber\\
&\geq \frac{(1-w^*)^2}{3} \left(q^*_{21}-q_{21}\right)^2.\label{eq:low_q21}
\end{align}
\item Similarly, we have
\begin{align}
(\pi_{11}-\pi^*_{11})^2 + (\pi_{12}-\pi^*_{12})^2 &= w^2 \left[ q_{12}^2 + (1-q_{12})^2 \right] + (w^*)^2 \left[ (q^*_{12})^2 + (1-q^*_{12})^2 \right]\nonumber\\
&- 2ww^* \left[ q_{12} q^*_{12} + (1-q_{12})(1-q^*_{12}) \right]\nonumber\\
&= \left[ q_{12}^2 + (1-q_{12})^2 \right] \left( w - w^* \frac{q_{12} q^*_{12} + (1-q_{12})(1-q^*_{12})}{q_{12}^2 + (1-q_{12})^2} \right)^2\nonumber\\
&+ (w^*)^2 \left[ (q^*_{12})^2 + (1-q^*_{12})^2 - \frac{\left(q_{12} q^*_{12} + (1-q_{12})(1-q^*_{12})\right)^2}{q_{12}^2 + (1-q_{12})^2} \right]\nonumber\\
&\geq \frac{(w^*)^2}{q_{12}^2 + (1-q_{12})^2} \left[ \left( (q^*_{12})^2 + (1-q^*_{12})^2 \right) \left( (q_{12})^2 + (1-q_{12})^2 \right)\right.\nonumber\\
&- \left. \left(q_{12} q^*_{12} + (1-q_{12})(1-q^*_{12})\right)^2 \right]\nonumber\\
&= \left(w^*\right)^2 \frac{\left( q_{12} - q^*_{12} \right)^2}{q_{12}^2 + (1-q_{12})^2}\nonumber\\
&\geq \left(w^*\right)^2 \left( q_{12} - q^*_{12} \right)^2.\label{eq:low_q12}
\end{align}
\end{itemize}
Finally, with (\ref{eq:low_w}),(\ref{eq:low_q21}) and (\ref{eq:low_q12}), we get
\begin{align*}
&(\pi_{11}-\pi^*_{11})^2 + (\pi_{12}-\pi^*_{12})^2 + (\pi_{21} -\pi^*_{21})^2\\
&\geq \max\left( \frac{1}{2} (w-w^*)^2, \frac{(1-w^*)^2}{3} \left(q^*_{21}-q_{21}\right)^2, \left(w^*\right)^2 \left( q_{12} - q^*_{12} \right)^2 \right).
\end{align*}
\section{Selection of the spacing parameter}
This section gathers the proofs of Theorem \ref{th:rho_selection_s}, \ref{th:selection_s_hmm}, Lemma \ref{lem:robust_beta} and Corollary \ref{coro:robust_hmm_selection_s}.
\subsection{Proof of Theorem \ref{th:rho_selection_s}}
\label{sec:proof_th_rho_selection_s}
We first need the following result.
\begin{lemme}
\label{lem:aux_beta_2}
Let $\mscr{M}$ be a finite set of probability distributions associated to the set of probability density functions $\mcal{M}$, with respect to the $\sigma$-finite measure $\mu$. Let $\hat{P}=\hat{P}\left(n,\mbf{X},\mcal{M}\right)$ be the $\rho$-estimator given by (\ref{eq:rho_estimateur}). For $t\in[n]$, there is an event $\Omega^*$ such that $\mbb{P}(\Omega^*)\geq 1-\lceil n/t\rceil\beta_t\left(\mbf{X}\right)$ and for all $\xi>0$, with probability at least $1-2|\mcal{M}|e^{-\xi}$, we have
\begin{align*}
\mathbbm{1}_{\Omega^*}\suml_{i=1}^n h^2\left(P_i,\hat{P} \right) &\leq \left( \frac{4a_0}{a_1} +1 \right) \inf_{Q\in\mscr{M}} \suml_{i=1}^n h^2 \left(P_i,Q\right)\\
&+ \frac{8}{3a_1} \left( \xi + 1.47 \right) \left[ 1 + \sqrt{1 + 18 t a_2^2 \alpha_0(t)} \right] + \frac{16.48}{a_1},
\end{align*}
with $\alpha_0(t)=\frac{32 \times 1.175 t a_2^2}{a_1^2} + \frac{8}{3 a_1}$, $a_0=4, a_1=3/8$ and $a_2^2=3\sqrt{2}$.
\end{lemme}
Consequently, we have
\begin{align*}
\mbb{E}\left[ \suml_{i=1}^n h^2\left(P_i,\hat{P} \right) \right] &\leq n \mbb{P}\left( \left(\Omega^*\right)^C \right) + \int_0^{\infty} \mbb{P}\left( \mathbbm{1}_{\Omega^*}\suml_{i=1}^n h^2\left(P_i,\hat{P} \right) \geq u \right) du\\
&\leq n \lceil n/t\rceil \beta_t\left(\mbf{X}\right) + \left( \frac{4a_0}{a_1} + 1 \right) \inf_{Q\in\mscr{M}}\suml_{i=1}^n h^2\left(P_i,\mscr{M}\right) + \frac{16.48}{a_1}\\
&+ \frac{8}{3a_1} \left(2.47 + \log(2|\mcal{M}|) \right) \left[ 1 + \sqrt{1 + 18 t a_2^2 \alpha_0(t)} \right].
\end{align*}
We apply this with $\mcal{M}=\widehat{\mcal{M}}_S\left(\mbf{X}^{(1)}\right)$ and conditionally on $\mbf{X}^{(1)}$. One can check that we have $\sqrt{ 1+18 t a_2^2 \alpha_0(t) } \leq 1 + 24 \frac{t a_2^2}{a_1} \sqrt{1.175} $. We get
\begin{align*}
\mbb{E}\left[ \suml_{i=1}^{n_2} h^2\left(P_i^{(2)},\hat{P}_{\hat{s}} \right) \big| \mbf{X}^{(1)} \right] &\leq c'_0 \inf_{s\in S} \suml_{i=1}^{n_2} h^2\left( P_i^{(2)}, \hat{P}_s\left(\mbf{X}^{(1)}\right) \right)\\
&+ c'_1 \left(2.47 + \log(2|S|) \right) \left[ 1 + 96\sqrt{2.35} t \right]\\
&+ c'_2 + n_2 \lceil n_2/t\rceil \beta_t\left(\mbf{X}^{(2)}\right),
\end{align*}
with $c'_0=\frac{4a_0}{a_1}+1=\frac{131}{3}$, $c'_1=\frac{2\times 8}{3 a_1}=\frac{128}{9}$ and $c'_2= \frac{16.48}{a_1} = \frac{131.84}{3}$. As $t$ can be any number in $[n_2]$ we can take the infimum with respect no $t$ in the upper bound. Let $\overline{P}$ be in $\mscr{P}_X$. We get
\begin{align*}
\mbb{E}\left[ h^2\left(\overline{P}, \hat{P}_{\hat{s}}\right) \right] &\leq \frac{2}{n_2} \mbb{E}\left[ \suml_{i=1}^{n_2} h^2\left(P_i^{(2)},\hat{P}_{\hat{s}} \right) \right] + \frac{2}{n_2} \suml_{i=1}^{n_2} h^2\left( P_i^{(2)},\overline{P} \right)\\
&\leq \frac{2}{n_2} \suml_{i=1}^{n_2} h^2\left( P_i^{(2)},\overline{P} \right) + \frac{2 c'_0}{n_2} \inf_{s\in S} \mbb{E}\left[ \suml_{i=1}^{n_2} h^2\left( P_i^{(2)},\hat{P}_s \right) \right]\\
&+ \inf_{t\in[n_2]} \left\{ \frac{c'_1}{n_2} \left(2.47 + \log(2|S|) \right) \left[ 1 + 96\sqrt{2.35} t \right] + 2 \lceil n_2/t\rceil \beta_t\left(\mbf{X}^{(2)}\right) \right\}\\
&+ \frac{2 c'_2}{n_2} .
\end{align*}
From (\ref{eq:th_1_comment}), for $s$ in $S$, we have
\begin{align*}
\frac{1}{n_2} \mbb{E}\left[ \suml_{i=1}^{n_2} h^2\left( P_i^{(2)},\hat{P}_s \right) \right] &\leq \frac{2}{n_2} \suml_{i=1}^{n_2} h^2\left( P_i^{(2)},\overline{P} \right) + \frac{4}{n_1} \suml_{i=1}^{n_1} h^2\left( P_i^{(1)},\overline{P} \right)\\
&+ \frac{4}{n_1} \mbb{E}\left[ \suml_{i=1}^{n_1} h^2\left( P_i^{(1)},\hat{P}_s \right) \right] \\
&\leq \frac{2}{n_2} \suml_{i=1}^{n_2} h^2\left( P_i^{(2)},\overline{P} \right) + \frac{4}{n_1} \suml_{i=1}^{n_1} h^2\left( P_i^{(1)},\overline{P} \right)\\
&+ \frac{4 c_0}{n_1} \inf_{Q\in\mscr{M}_s} \suml_{i=1}^{n_1} h^2\left( P_i^{(1)}, Q \right) + 4 c_1 \frac{(s+1)}{n_1} \left[ 17 + D_{n(s,1)}(\mscr{M}_s) \right]\\
&+  \frac{4 c_2}{n_1} \suml_{b=1}^{s+1} \mbf{K}\left(\mbf{P}^*_{s,b}||\mbf{P}^{ind}_{s,b}\right).
\end{align*}
We get
\begin{align*}
\mbb{E}\left[ h^2\left(\overline{P}, \hat{P}_{\hat{s}}\right) \right] &\leq \frac{2 + 4 c'_0}{n_2} \suml_{i=1}^{n_2} h^2\left( P_i^{(2)},\overline{P} \right) + \frac{8 c'_0}{n_1} \suml_{i=1}^{n_1} h^2\left( P_i^{(1)},\overline{P} \right)\\
&+ \inf_{t\in[n_2]} \left\{ \frac{c'_1}{n_2} \left(2.47 + \log(2|S|) \right) \left[ 1 + 96\sqrt{2.35} t \right] + 2 \lceil n_2/t\rceil \beta_t\left( \mbf{X}^{(2)} \right) \right\}\\
&+ \frac{2 c'_2}{n_2} + \frac{8 c'_0}{n_1} \inf_{s\in S} \bigg\{ c_0 \inf_{Q\in\mscr{M}_s} \suml_{i=1}^{n_1} h^2\left( P_i^{(1)}, Q \right)\\
&+ c_1 (s+1) \left[ D_{n(s,1)}(\mscr{M_s}) + 17 \right] +  c_2 \suml_{b=1}^{s+1} \mbf{K}\left(\mbf{P}^*_{s,b}||\mbf{P}^{ind}_{s,b}\right) \bigg\}.
\end{align*}
We also have
\[ \frac{1}{n_1} \inf_{Q\in\mscr{M}_s} \suml_{i=1}^{n_1} h^2\left( P_i^{(1)},Q \right) \leq 2 h^2(\overline{P},\mscr{M}_s) + \frac{2}{n_1} \suml_{i=1}^{n_1} h^2\left( P_i^{(1)},\overline{P} \right). \]
\subsubsection{Proof of Lemma \ref{lem:aux_beta_2}}
For $P_i=\mcal{L}(X_i),i=1,\dots,n,$ we write
\[ H^2_{Q,Q'} := \suml_{i=1}^n h^2\left( P_i,Q \right) + h^2\left( P_i,Q' \right).\]
\begin{lemme}
\label{lem:beta_aux_2}
Let $\delta>1$ and $\nu>0$ be such that
\[ e^{-\nu} + \suml_{j\geq 1} e^{-\delta^j \nu} \leq 1.\]
For $t$ in $\{1,\dots,n\}$, there is an event $\Omega^*$ satisfying $\mbb{P}(\Omega^*)\geq 1-\lceil n/t\rceil\beta_t$ such that for all $p$ in $\mcal{M}$ and all $\xi>0$,  we have
\[ \mbf{P}^*\left( \sup_{q\in\mcal{M}} \left\{ \left| \mbf{Z}_n(\mbf{X},p,q) \right| \mathbbm{1}_{\Omega^*} - \frac{a_1}{2} H^2_{P,Q} \right\} > \frac{2 (\upsilon + \xi)}{3} \left[ 1 + \sqrt{1 + 18 t a_2^2 \alpha} \right] \right) \leq 2|\mcal{M}| e^{-\xi},\]
with $\mbf{P}^*=\mcal{L}(\mbf{X})$ and $\alpha\geq\alpha_0(t) = \frac{ 32 t a_2^2 \delta}{a_1^2} + \frac{8}{3 a_1}$.
\end{lemme}
We take $\delta=1.175$ and $\upsilon=1.47$ as in \cite{baraudrevisited} Section A.1. Let $\xi>0$ and $p\in\mcal{M}$. On the event $\Omega^*$ defined by Lemma \ref{lem:beta_aux_2} and with Proposition 3 \cite{baraudrevisited}, we have
for all $q\in\mcal{M}$,
\begin{align*}
\mbf{T}_n\left(\mbf{X},p,q\right) &\leq \mbb{E} \mbf{T}_n\left(\mbf{X},p,q\right) + |\mbf{Z}\left(\mbf{X},p,q\right)|\\
&\leq \suml_{i=1}^n \left[ a_0  h^2\left(P_i,P\right) - a_1 h^2\left(P_i,Q\right) \right]\\
&+  \frac{a_1}{2} H^2_{P,Q} + \frac{2 (\xi+\upsilon) }{3} \left[ 1 + \sqrt{1 + 18 t a_2^2 \alpha_0(t) } \right]\\
&= \suml_{i=1}^n \left[ \left( a_0 + \frac{a_1}{2} \right) h^2\left(P_i,P\right) - \frac{a_1}{2} h^2\left(P_i,Q\right) \right]\\
&+ \frac{2}{3} (\xi+\upsilon) \left[ 1 + \sqrt{1 + 18 t a_2^2 \alpha_0(t)} \right].
\end{align*}
Then,
\begin{align*}
\mbf{\Upsilon}_n\left(\mbf{X},p\right) &= \sup_{q\in\mcal{M}} \mbf{T}_n\left(\mbf{X},p,q\right)\\
&\leq \left( a_0 + \frac{a_1}{2} \right) \suml_{i=1}^n h^2\left(\mbf{P}^{ind}_i,P \right)\\
&- \frac{a_1}{2} \inf_{Q\in\mscr{M}} \suml_{i=1}^n h^2\left( P_i,Q\right)\\
&+ \frac{2}{3} (\xi+\upsilon) \left[ 1 + \sqrt{1 + 18 t a_2^2 \alpha_0(t)} \right],
\end{align*}
and
\begin{align*}
\mbf{\Upsilon}_n\left(\mbf{X},q\right) &= \sup_{q'\in\mcal{M}} \mbf{T}_n\left(\mbf{X},q,p\right)\\
&\geq \mbf{T}_n\left(\mbf{X},q,p\right) = - \mbf{T}_n\left(\mbf{X},p,q\right)\\
&\geq -\left( a_0 + \frac{a_1}{2} \right) \suml_{i=1}^n h^2\left( P_i,P \right)  + \frac{a_1}{2} \suml_{i=1}^n h^2\left( P_i,Q \right)\\
&- \frac{2}{3} (\xi+\upsilon) \left[ 1 + \sqrt{1 + 18 t a_2^2 \alpha_0(t)} \right].
\end{align*}
Since $\mbf{\Upsilon}_n\left(\mbf{X},\hat{p}\right) < \mbf{\Upsilon}_n\left(\mbf{X},p\right) + 8.24$, we have
\begin{align*}
\frac{a_1}{2} \suml_{i=1}^n h^2\left( P_i, \hat{P} \right) &\leq 2 \left( a_0 + \frac{a_1}{2} \right) \suml_{i=1}^n h^2 \left( P_i,P \right) - \frac{a_1}{2} \inf_{Q\in\mscr{M}} \suml_{i=1}^n h^2 \left( P_i,\mscr{M} \right)\\
&+  \frac{4}{3} (\xi+\upsilon) \left[ 1 + \sqrt{1 + 18 t a_2^2 \alpha_0(t)} \right] + 8.24.
\end{align*}
Given that $\mscr{M}$ is finite we can take $P$ such that
\[ \inf_{Q\in\mscr{M}} \suml_{i=1}^n h^2\left( P_i,Q \right) = \suml_{i=1}^n h^2\left(P_i, P\right).\]
Hence we have
\begin{align*}
\suml_{i=1}^n h^2\left( P_i, \hat{P} \right) &\leq \left( \frac{4 a_0}{a_1} + 1 \right) \inf_{Q\in\mscr{M}} \suml_{i=1}^n h^2\left( P_i,Q \right)\\
&+  \frac{8}{3 a_1} (\xi+\upsilon) \left[ 1 + \sqrt{1 + 18 t a_2^2 \alpha_0(t)} \right] + \frac{16.48}{a_1}.
\end{align*}
\subsubsection{Proof of Lemma \ref{lem:beta_aux_2}}
\begin{lemme}
\label{lem:beta_aux_1}
For $t$ in $[n]$, there is an event $\Omega^*$ such that $\mbb{P}(\Omega^*)\geq 1-\lceil n/t\rceil\beta_t(\mbf{X})$ and
\begin{equation}
\forall q,q'\in\mcal{M},\forall x>0,\mbb{P}\left( \left|\mbf{Z}_n\left(\mbf{X},q,q'\right)\right| \mathbbm{1}_{\Omega^*} > \frac{2x}{3} \left[ 1 + \sqrt{1+\frac{ 18 t a_2^2 H^2_{Q,Q'}}{x}} \right] \right) \leq 2e^{-x}.
\end{equation}
\end{lemme}
Let $\xi>0$ and $\alpha >0$. We define $x_0=\upsilon+\xi$ and for $j\geq 0$, 
\begin{equation}
y_{j+1}^2=\delta y_j^2= \delta \alpha x_j.
\end{equation}
Let $q,q'$ be in $\mcal{M}$. We apply Lemma \ref{lem:beta_aux_1} according to the value of $H^2_{Q,Q'}$.
\begin{itemize}
\item If there is $j\geq 0$ such that $y_j^2\leq H^2_{Q,Q'} <y_{j+1}^2$, with probability at least $1-2e^{-x_j}$, we have
\begin{align*}
|\mbf{Z}_n(\mbf{X},q,q')|\mathbbm{1}_{\Omega^*} - \frac{a_1}{2} H^2_{Q,Q'} &\leq \frac{2 x_j}{3} \left[ 1 + \sqrt{1+\frac{ 18 t a_2^2 H^2_{Q,Q'}}{x_j}} \right] - \frac{a_1}{2} H^2_{q,q'}\\
&\leq \frac{2 x_j}{3} \left[ 1 + \sqrt{1+\frac{ 18 t a_2^2 y_{j+1}^2}{x_j}} \right] - \frac{a_1}{2} y_j^2\\
&\leq \frac{2 x_j}{3} \left[ 1 + \sqrt{1 + 18 t a_2^2 \delta\alpha } - \frac{ 3 a_1 \alpha}{4} \right]\\
&\leq 0,
\end{align*}
for
\begin{equation}
\alpha \geq \alpha_0(t) := \frac{32 \delta t a_2^2}{a_1} + \frac{8}{3a_1}.
\end{equation}
\item If $H^2_{Q,Q'} < y_0^2$, with probability at least $1-2e^{-x_0}$, we have
\begin{align*}
|\mbf{Z}_n(\mbf{X},q,q')| \mathbbm{1}_{\Omega^*} - \frac{a_1}{2} H^2_{Q,Q'} &\leq |\mbf{Z}_n(\mbf{X},q,q')|\mathbbm{1}_{\Omega^*}\\
&\leq \frac{2 x_0}{3} \left[ 1 + \sqrt{1 + 18 t a_2^2 \alpha } \right].
\end{align*}
\end{itemize}
Let $\overline{p}$ be in $\mcal{M}$. Finally, we have
\begin{align*}
&\mbb{P}\left( \sup_{q\in\mcal{M}} \left\{ \left| \mbf{Z}_n(\mbf{X},\overline{p},q) \right| \mathbbm{1}_{\Omega^*} - \frac{a_1}{2} H^2_{\overline{p},q} \right\} > \frac{2 x_0}{3} \left[ 1 + \sqrt{1 + 18 t a_2^2 \alpha} \right] \right)\\
&\leq \suml_{ \substack{q\in\mcal{M}: \\ H^2_{\overline{P},Q} < y_0^2} } \mbb{P}\left( |\mbf{Z}_n(\mbf{X},\overline{p},q)|\mathbbm{1}_{\Omega^*} - \frac{a_1}{2} H^2_{\overline{p},q} > \frac{2 x_0}{3} \left[ 1 + \sqrt{1 + 18 t a_2^2 \alpha} \right] \right)\\
&+ \suml_{j \geq 0} \suml_{\substack{q\in\mcal{M} :\\ y_j^2 \leq H^2_{\overline{p},q} < y_{j+1}^2}} \mbb{P}\left( \left|\mbf{Z}_n(\mbf{X},\overline{p},q)\right| \mathbbm{1}_{\Omega^*} - \frac{a_1}{2} H^2_{\overline{p},q} > 0 \right)\\
&\leq \suml_{\substack{q\in\mcal{M}:\\ H^2_{\overline{p},q} < y_0^2}} 2e^{-x_0} + \suml_{j\geq 0} \suml_{ \substack{ q\in\mcal{M}:\\ y_j^2 \leq H^2_{\overline{p},q} < y_{j+1}^2} } 2e^{-x_j}\\
&\leq 2|\mcal{M}| \left( e^{-x_0} + \suml_{j\geq 1} e^{-x_j}\right) = 2|\mcal{M}| \left( e^{-(\upsilon+\xi)} + \suml_{j\geq 1} e^{-\delta^j (\upsilon+\xi)}\right)\\
&\leq 2 |\mcal{M}| e^{-\xi} \left( e^{-\upsilon} + \suml_{j\geq 1} e^{-\delta^j \upsilon}\right) \leq 2|\mcal{M}|e^{-\xi}.
\end{align*}
\subsubsection{Proof of Lemma \ref{lem:beta_aux_1}}
We follow the proof of Sart \cite{sarttransition} (Proposition B.1). Let $t$ be a positive integer in $[n]$. Let $l$ be the smallest integer larger than $\slfrac{n}{2t}$. We derive from Berbee’s lemma and more precisely from Viennet [36] (page 484) that there exist $B^*_1,\dots,B^*_{2lt}$ such that
\begin{itemize}
\item For $i = 1,\dots,l$, the random vectors
\begin{equation}
B_{i,1} = \left( X_{2(i-1)t+1},\dots,X_{(2i-1)t} \right) \text{ and } B^*_{i,1} = \left( X^*_{2(i-1)t+1},\dots,X^*_{(2i-1)t} \right)
\end{equation}
have the same distribution, and so have the random vectors
\begin{equation}
B_{i,2} =\left( X_{(2i-1)t+1},\dots,X_{2it} \right)\text{ and }B^*_{i,2} = \left( X^*_{(2i-1)t+1},\dots,X^*_{2it} \right).
\end{equation}
\item The random vectors $B^*_{1,1},\dots,B*_{l,1}$ are independent. The random vectors $B^*_{1,2},\dots,B*_{l,2}$ are also independent.
\item The event
\[ \Omega^* = \bigcap_{1\leq j\leq l} \left\{ B_{j,1} = B^*_{j,1} \right\} \cap \left\{ B_{j,2} = B^*_{j,2} \right\} \]
satisfies $\mbb{P}\left((\Omega^*)^C\right) \leq 2l\beta_t\left( \mbf{X} \right)$.
\end{itemize}
Let $q,q'$ be in $\mcal{M}$. For simplicity, we write $Z_{q,q'}=\mbf{Z}(\mbf{B},q,q')$ and we define
\begin{align*}
Z^*_{q,q',1} &:= \suml_{i=1}^{l} \suml_{j=1}^t \left\{ \psi\left( \sqrt{\frac{q'}{q}\left(X^*_{2(i-1)t+j}\right)}\right) - \mbb{E} \left[ \psi\left( \sqrt{\frac{q'}{q}\left(X^*_{2(i-1)t+j}\right)}\right) \right] \right\} \mathbbm{1}_{2(i-1)t+j\leq n}\\
&= \suml_{i=1}^{l} \suml_{j=1}^t z^{q,q'}_{2(i-1)t+j} \mathbbm{1}_{2(i-1)t+j\leq n}
\end{align*}
and
\begin{align*}
Z^*_{q,q',2} &:= \suml_{i=1}^{l} \suml_{j=1}^t  \left\{ \psi\left( \sqrt{\frac{q'}{q}\left(X^*_{(2i-1)t+j}\right)}\right) - \mbb{E} \left[ \psi\left( \sqrt{\frac{q'}{q}\left(X^*_{(2i-1)t+j}\right)}\right) \right] \right\} \mathbbm{1}_{(2i-1)t+j\leq n}\\
&= \suml_{i=1}^{l} \suml_{j=1}^m z^{q,q'}_{(2i-1)t+j} \mathbbm{1}_{(2i-1)t+j\leq n}.
\end{align*}
Let $\xi$ be a positive real number. Since
\begin{equation}
|Z_{q,q'}|\mathbbm{1}_{\Omega^*} > \xi \Rightarrow |Z^*_{q,q',1}|\mathbbm{1}_{\Omega^*} > \xi/2 \text{ or } |Z^*_{q,q',2}|\mathbbm{1}_{\Omega^*} > \xi/2,
\end{equation}
we have
\begin{align*}
\mbb{P}\left( \left|Z_{q,q'}\right|\mathbbm{1}_{\Omega^*} > \xi \right) &\leq \mbb{P}\left(|Z^*_{q,q',1}|\mathbbm{1}_{\Omega^*} > \xi/2 \right) + \mbb{P}\left(|Z^*_{q,q',2}|\mathbbm{1}_{\Omega^*} > \xi/2 \right)\\
&\leq \mbb{P}\left( |Z^*_{q,q',1}| > \xi/2 \right) + \mbb{P}\left( |Z^*_{q,q',2}| > \xi/2 \right).
\end{align*}
One can notice that $Z^*_{q,q',1}$ and $Z^*_{q,q',2}$ are sums of $l$ independent variables. Therefore, we can use classic concentration inequalities. First, we can see that
\begin{align*}
V_{q,q',1} &= \suml_{i=1}^{l} \mbb{E}\left[  \left( \suml_{j=1}^t z_{q,q'}^{2(i-1)t+j} \mathbbm{1}_{2(i-1)t + j} \right)^2 \right]\\
&\leq \suml_{i=1}^{l} \suml_{j=1}^t t \mbb{E}\left[   \left( z_{q,q'}^{2(i-1)t+j}\right)^2 \mathbbm{1}_{2(i-1)t + j} \right]\\
&\leq t \suml_{i=1}^{n} \text{Var}\left(  \psi\left( \sqrt{ \frac{q'}{q}\left( X^*_i \right) } \right) \right)\\
&\leq t \suml_{i=1}^n a_2^2 \left[ h^2(P_i,Q) + h^2(P_i,Q') \right] = t a_2^2 H^2_{Q,Q'}.
\end{align*}
The last inequality comes from Proposition 3 in Baraud \& Birg\'e \cite{baraudrevisited} and $a_2^2=3\sqrt{2}$. Similarly  we have $V_{Q,Q',2} \leq t a_2^2 L_{Q,Q'}$. Therefore, Bennett's inequality (see Proposition 2.8 and inequality (2.16) in Massart \cite{Massart}) guarantees that for all $\xi>0$ we have
\begin{align*}
\mbb{P}\left(|Z_{q,q'}|\mathbbm{1}_{\Omega^*} > \xi \right) \leq 2\exp\left( -\frac{(\xi/2)^2}{2 (t a_2^2 H^2_{q,q'} + \xi/6)}\right).
\end{align*}
For $x>0$, we take $\xi=\frac{2x}{3} \left[ 1 + \sqrt{1+\frac{18 t a_2^2 H^2_{Q,Q'}}{x}} \right]$ and with probability less than or equal to $2e^{-x}$, we have
\begin{equation}
|Z_{q,q'}|\mathbbm{1}_{\Omega^*} > \frac{2x}{3} \left[ 1 + \sqrt{1+\frac{ 18 t a_2^2 H^2_{Q,Q'}}{x}} \right].
\end{equation}
\subsection{Proof of Lemma \ref{lem:robust_beta}}
\label{sec:proof_lem_robust_beta}
We have
\begin{align*}
\beta_t\left(\mbf{Y}\right) &= \sup_i \beta\left( \sigma(Y_1,\dots,Y_i); \sigma(Y_{i+t},\dots,Y_n)\right)\\
&= \sup_i d_{TV}\left( \mcal{L}\left(Y_1,\dots,Y_i\right)\otimes\mcal{L}\left(Y_{i+t},\dots,Y_n\right), \mcal{L}\left(Y_1,\dots,Y_i,Y_{i+t},\dots,Y_n\right) \right).
\end{align*}
We use the notation $X_a^b=(X_a,\dots,X_b)$ and similarly for $\mbf{E}$, $\mbf{Y}$ and $\mbf{Z}$. The triangle inequality implies
\begin{align*}
&d_{TV}\left( \mcal{L}\left(Y_1^i\right)\otimes\mcal{L}\left(Y_{i+t}^n\right), \mcal{L}\left(Y_1^n\right) \right)\\
&\leq \suml_{\mbf{e}\in\{0,1\}^n} \mbb{P}\left( \mbf{E}=\mbf{e} \right) d_{TV}\left( \mcal{L}( Y_1^i | E_1^i=e_1^i ) \otimes \mcal{L}(Y_{i+t}^n | E_{i+t}^N=e_{i+t}^n), \mcal{L}( Y_1^i,Y_{i+t}^n | E_1^i=e_1^i,E_{i+t}^n=e_{i+t}^n ) \right)\\
&= \suml_{\mbf{e}\in\{0,1\}^n} \mbb{P}\left( \mbf{E}=\mbf{e} \right) \beta\left( \sigma( (X_j)_{\substack{j\leq i,\\e_j=1}} ), \sigma( (X_j)_{\substack{j\geq i+k,\\e_j=1}} ) \right).
\end{align*}
We now need the following result to conclude.
\begin{lemme}
\label{lem:beta_formula}
For any random variables $A_1,A_2,B_1,B_2$, we have
\begin{equation*}
\beta\left(\sigma(A_1),\sigma(A_2)\right) \leq \beta\left( \sigma(A_1,B_1), \sigma(A_2,B_2) \right).
\end{equation*}
\end{lemme}
Combining the different inequalities above, we get
\begin{align*}
\beta_t\left(\mbf{Y}\right) &\leq \sup_i \beta\left( \sigma(Y_1^i); \sigma(Y_{i+t}^n) \right)\\
&= \sup_i \suml_{\mbf{e}\in\{0,1\}^n} \mbb{P}\left( \mbf{E}=\mbf{e} \right) \beta\left( \sigma( (X_j)_{\substack{j\leq i,\\e_j=1}} ), \sigma( (X_j)_{\substack{j\geq i+t,\\e_j=1}} ) \right)\\
&\leq \sup_i \suml_{\mbf{e}\in\{0,1\}^n} \mbb{P}\left( \mbf{E}=\mbf{e} \right) \beta\left( \sigma( (X_j)_{j\leq i} ), \sigma( (X_j)_{j\geq i+t} ) \right) = \beta_t\left(\mbf{X}\right).
\end{align*}
\subsubsection{Proof of Lemma \ref{lem:beta_formula}}
Let $\mu_1$, $\mu_2$, $\nu_1$ and $\nu_2$ be measures dominating respectively $\mcal{L}(A_1)$, $\mcal{L}(A_2)$, $\mcal{L}(B_1)$ and $\mcal{L}(B_2)$. We have
\begin{align*}
&\beta\left(\sigma(A_1),\sigma(A_2)\right)\\
&= \frac{1}{2} \int |p_A(a_1,a_2) - p_{A_1}(a_1) p_{A_2}(a_2) | \mu_1(da_1) \mu_2(da_2) \\
&= \frac{1}{2} \int | \int \left(p_{A,B}(a_1,b_1,a_2,b_2) - p_1(a_1,b_1) p_2(a_2,b_2) \right) \nu_1(db_1) \nu_2(db_2) | \mu_1(da_1) \mu_2(da_2)\\
&\leq \frac{1}{2} \int |p_{A,B}(a_1,b_1,a_2,b_2) - p_1(a_1,b_1) p_2(a_2,b_2) | \nu_1(db_1) \nu_2(db_2)  \mu_1(da_1) \mu_2(da_2\\
&= \beta\left( \sigma(A_1,B_1); \sigma(A_2,B_2) \right), 
\end{align*}
with $p_A = \frac{d\mcal{L}(A_1,A_2)}{d\mu_1\otimes\mu_2}$, $p_{A_1} = \frac{d\mcal{L}(A_1)}{d\mu_1}$, $p_{A_2} = \frac{d\mcal{L}(A_2)}{d\mu_2}$, $p_{A,B} = \frac{d\mcal{L}(A_1,B_1,A_2,B_2)}{d\mu_1\otimes \nu_1\otimes\mu_2\otimes \nu_2}$, $p_1 = \frac{d\mcal{L}(A_1,B_1)}{d\mu_1\otimes\nu_1}$  and $p_2 = \frac{d\mcal{L}(A_2,B_2)}{d\mu_2\otimes\mu_2}$.
\subsection{Proof of Theorem \ref{th:selection_s_hmm}}
\label{sec:proof_th_selection_s_hmm}
From (\ref{eq:approximation_delta}) we have
\begin{align*}
h^2\left( \overline{P},\mscr{M}_s\right) &\leq 2 L\epsilon^2 + 2 L(K-1)\delta(s) + 2 h^2\left( \overline{P}, \overline{\mscr{M}} \right)\\
&\leq 2 L\epsilon^2 + 2 h^2\left( \overline{P}, \mscr{M} \right) + 2 (s+1) L \frac{\overline{V}}{n_1}.
\end{align*}
From Proposition \ref{prop:rho_dimension_hmm} we have $D_{n_1(s,1)}\left(\mscr{M}_{s}\right) \leq C L \overline{V} \log n_1$, for a constant $C$. For $S$ defined by (\ref{eq:S}), we have
\[ |S| = 2 + \lfloor \log_{\tau}( \lfloor (n_1-2)/2 \rfloor ) \rfloor \leq 2 + \frac{\log n_1}{\log \tau} \leq C \log n_1, \]
for some positive constant $C$. Theorem \ref{th:rho_selection_s} allows to obtain (\ref{eq:th_selection_s_hmm}).\par
The following result is proven in Section \ref{sec:proof_lem_hyp_exponential_ergodicity_s}.
\begin{lemme}
\label{lem:hyp_exponential_ergodicity_s}
Under Assumption \ref{hyp:hmm_ergodic_s}, there exist positive constants $r(Q^*),C(Q^*)>0$ such that
\begin{itemize}
\item for all $j\in[2]$ and all $i\in[n_j]$, we have
\begin{equation}
\label{eq:hyp_s_1}
h^2\left(P_i^{(j)},P^*\right)\leq C(Q^*) e^{-r(Q^*) i},
\end{equation}
\item for all $t\in[n_2]$, we have
\begin{equation}
\label{eq:hyp_s_2}
\beta_t\left(\mbf{X}^{(2)}\right) \leq C(Q^*) e^{- r(Q^*) t/2},
\end{equation}
\item for all $s\geq L-1$, all $b$ in $[s+1]$,
\begin{equation}
\label{eq:hyp_s_3} 
\mbf{K}\left(\mbf{P}^*_{s,b}||\mbf{P}^{ind}_{s,b}\right) \leq n(s,b) C(Q^*) e^{-r(Q^*) s}.
\end{equation}
\end{itemize}
\end{lemme}
From (\ref{eq:hyp_s_1}) we have
\[ \suml_{i=1}^{n_1} h^2\left(P^{(1)}_i,P^*\right), \suml_{i=1}^{n_1} h^2\left(P^{(1)}_i,P^*\right) \leq \frac{C(Q^*)}{e^{r(Q^*)}-1}. \]
For $t=n_2\wedge \left\lceil 4 r(Q^*)^{-1} \log n_2 \right\rceil$, with (\ref{eq:hyp_s_2}) we have
\begin{align*}
\lceil n_2/t \rceil \beta_t\left( \mbf{X}^{(2)} \right) &\leq \begin{cases}
1 \text{  for  } n_2\leq r(Q^*)^{-1} 4 \log n_2,\\
C(Q^*) n_2^{-1} \text{  otherwise},
      \end{cases}\\
      &\leq n_2^{-1} \left( C(Q^*) \vee r(Q^*)^{-1} 4\log n_2 \right).
\end{align*}
We have the following
\begin{align*}
\left\lceil \frac{\log \log n_1 - \log r(Q^*)}{\log \tau} \right\rceil > \left\lfloor \frac{\log \left\lfloor \frac{n_1-2}{2} \right\rfloor}{\log\tau} \right\rfloor &\Rightarrow \frac{\log \log n_1 - \log r(Q^*)}{\log\tau} > \frac{\log \left\lfloor \frac{n_1-2}{2} \right\rfloor}{\log\tau} - 1\\
&\Rightarrow  \tau r(Q^*)^{-1} \log n_1 \geq \left\lfloor \frac{n_1-2}{2} \right\rfloor\\
&\Rightarrow  2 \frac{2+\tau r(Q^*)^{-1}\log n_1}{n_1} \geq 1.
\end{align*}
For $s=\lceil \tau^j \rceil$ with $j=\left\lceil \frac{\log \log n_1 - \log r(Q^*)}{\log\tau} \right\rceil \wedge \left\lfloor \frac{\log \left\lfloor \frac{n_1-2}{2} \right\rfloor}{\log\tau} \right\rfloor$, we have
\[ s \leq \tau^{\frac{\log \log n_1 - \log r(Q^*)}{\log \tau}+1} + 1 = 1 + \tau r(Q^*)^{-1} \log n_1, \]
and inequality (\ref{eq:hyp_s_3}) gives
\begin{align*}
\suml_{b=1}^{s+1} \mbf{K}\left(\mbf{P}^*_{s,b}||\mbf{P}^{ind}_{s,b}\right) &\leq C(Q^*) n_1 e^{-r(Q^*)s}\\
&\leq C(Q^*) n_1 \left( 2 \frac{2+\tau r(Q^*)^{-1} \log n_1}{n_1} \vee  \frac{1}{n_1} \right) = 2 C(Q^*) (2+\tau r(Q^*)^{-1} \log n_1).
\end{align*}
These last inequalities give (\ref{eq:th_selection_s_hmm_2}).
\subsubsection{Proof of Lemma \ref{lem:hyp_exponential_ergodicity_s}}
\label{sec:proof_lem_hyp_exponential_ergodicity_s}
We just have to follow the proof of Lemma \ref{lem:k_s_hmm}. We already have (\ref{eq:hyp_s_1}) and (\ref{eq:hyp_s_3}). The inequality (\ref{eq:hyp_s_2}) can be deduced from the inequality
\[ d_{TV}\left(Q^t_{k,\cdot},\pi\right) \leq C e^{-rt}, \]
and from the definition of $\beta_t$.
\subsection{Proof of Corollary \ref{coro:robust_hmm_selection_s}}
\label{sec:proof_coro_robust_hmm_selection_s}
We have
\[ \mbb{P}\left( X^{(j)}_i = \left( \overline{Y}^{(j)}_i,\ldots,\overline{Y}^{(j)}_{i+L-1} \right) \right) \geq \mbb{P}\left( E^{(j)}_i = \dots = E^{(j)}_{i+L-1} = 1 \right) = p^{(j)}_i p^{(j)}_{i+1} \dots p^{(j)}_{i+L-1}, \]
and with the convexity of the squared Hellinger distance
\begin{align*}
h^2\left( P^{(j)}_i, P^* \right) &\leq p^{(j)}_i p^{(j)}_{i+1} \dots p^{(j)}_{i+L-1} h^2\left( \overline{P}^{(j)}_i,P^* \right) + \left( 1 - p^{(j)}_i p^{(j)}_{i+1} \dots p^{(j)}_{i+L-1} \right)\\
&\leq h^2\left( \overline{P}^{(j)}_i, P^* \right) + \left( 1 - p^{(j)}_i \right) + \dots + \left( 1 - p^{(j)}_{i+L-1} \right),
\end{align*}
where $\overline{P}^{(j)}_i = \mcal{L}\left( \overline{Y}^{(j)}_i,\dots,\overline{Y}^{(j)}_{i+L-1} \right)$. One can check that $n\geq 1+N/2$ with our conditions on $L$. With Theorem \ref{th:selection_s_hmm}, Lemma \ref{lem:robust_beta} and Lemma \ref{lem:hyp_exponential_ergodicity_s} we have
\begin{align*}
C \mbb{E} \left[ h^2\left( P^*, \hat{P}_s \right) \right] &\leq h^2\left( P^*, \mscr{M} \right) + \frac{C(Q^*)}{n_1(e^{r(Q^*)}-1)} + \frac{C(Q^*)}{n_2(e^{r(Q^*)}-1)}\\
&+ L \epsilon^2 + \frac{L}{N_1} \suml_{i=1}^{N_1} \left(1-p^{(1)}_i\right) + \frac{L}{N_2} \suml_{i=1}^{N_2} \left(1-p^{(2)}_i\right)\\
&+ \inf\limits_{t\in[n_2]}\left\{ \frac{t\log \log n_1}{n_2} + \lceil n_2/t\rceil C(Q^*) e^{-r(Q^*)t/2} \right\}\\
&+ \inf\limits_{s\in S}\left\{ (s+1) L \overline{V} \frac{\log n_1}{n_1} + e^{-r(Q^*)s} \right\},
\end{align*}
for some positive constant $C$ and $s\geq L-1$. We can control the last terms with reasonable choices of $t$ and $s$ following the proof of Theorem \ref{th:selection_s_hmm}.
\end{document}